\documentclass[11pt]{article}
\usepackage[utf8]{inputenc}
\usepackage{amsmath}
\usepackage{amsfonts,enumerate}
\usepackage{color}
\usepackage{amssymb}
\usepackage{mathrsfs}
\usepackage{esint}
\usepackage[colorinlistoftodos]{todonotes}
\usepackage[colorlinks=true]{hyperref}
\hypersetup{
	colorlinks=blue,%
	citecolor=red,%
	filecolor=black,%
	linkcolor=blue,%
	urlcolor=blue
}
\usepackage[margin=1.3in]{geometry}

\usepackage[amsmath,thmmarks,hyperref]{ntheorem}
{
	\theoremstyle{nonumberplain}
	\theoremheaderfont{\bfseries}
	\theorembodyfont{\normalfont}
	\theoremsymbol{\mbox{$\Box$}}
	\newtheorem{pf}{Proof.}
}

\numberwithin{equation}{section}
\def\R{\mathbb{R}}
\def\S{\mathbb{S}}
\def\N{\mathbb{N}}
\def\B{\mathbb{B}}

\def \no{\nonumber}

\def\ve{\varepsilon}
\newcommand{\ud}{\mathrm{d}}
\newcommand{\Sn}{\mathbb{S}^n}
\newcommand{\Rn}{\mathbb{R}^n}
\def\pa{\partial}
\newtheorem*{conjecture}{Conjecture.}
\newtheorem{thm}{Theorem}[section]
\newtheorem{definition}{Definition}[section]

\newtheorem{lem}{Lemma}[section]
\newtheorem{rem}{Remark}[section]

\newtheorem{prop}{Proposition}[section]

\newtheorem*{convention}{\textbf{Convention.}}

\newtheorem{thm A}{Theorem A}
\newtheorem{cor}{Corollary}[section]
\usepackage{upgreek}
\usepackage{graphicx}

\makeatletter

\newdimen\bibspace
\renewenvironment{thebibliography}[1]{%
	\section*{\refname 
		\@mkboth{\MakeUppercase\refname}{\MakeUppercase\refname}}%
	\list{\@biblabel{\@arabic\c@enumiv}}%
	{\settowidth\labelwidth{\@biblabel{#1}}%
		\leftmargin\labelwidth
		\advance\leftmargin\labelsep
		\itemsep\bibspace
		\parsep\z@skip     %
		\@openbib@code
		\usecounter{enumiv}%
		\let\p@enumiv\@empty
		\renewcommand\theenumiv{\@arabic\c@enumiv}}%
	\sloppy\clubpenalty4000\widowpenalty4000%
	\sfcode`\.\@m}
{\def\@noitemerr
	{\@latex@warning{Empty `thebibliography' environment}}%
	\endlist}

\makeatother

\makeatletter

\allowdisplaybreaks

\title{article}
\begin{document}
	\title{\bf Classification of  solutions to \\the $Q$-flat and constant $T$-curvature equation \\on the half-space and ball} 
	\author{Xuezhang Chen\thanks{X. Chen: xuezhangchen@nju.edu.cn. Both authors are partially supported by NSFC (No.12271244). }~
		and Shihong Zhang\thanks{S. Zhang: dg21210019@smail.nju.edu.cn.}\\
		{\small $^{\ast}$$^{\dag}$School of Mathematics \& IMS, Nanjing University, Nanjing 210093, P. R. China }}
	\date{}

	\maketitle

	{\noindent\small{\bf Abstract}:
	For conformal boundary operators associated with the Paneitz operator, we introduce a rigorous definition of the biharmonic Poisson kernel consisting of a pair of kernel functions and derive its explicit representation formula. With this powerful tool, we establish classification theorems of nonnegative solutions to the $Q$-flat and constant $T$-curvature equations on $\R_+^{n+1}$ and $\B^{n+1}$.
	}

		\medskip 
		
		{{\bf $\mathbf{2020}$ MSC:} 35J40, 53C18 (35J91, 53C24)}
		
		\medskip 
		{\small{\bf Keywords:}
			 $Q$-curvature, conformal boundary operator, $T$-curvature, biharmonic Poisson kernel.}
		
		{\footnotesize \tableofcontents}	
		
		\section{Introduction}

		In a seminal work \cite{Caffarelli&Gidas&Spruck}, Caffarelli-Gidas-Spruck   established  a Liouville-type theorem of the second-order critical semilinear elliptic equation: For $n\geq 2$, every positive $C^2$ solution to 
	       \begin{align*}
			-\Delta u=(n+1)(n-1)u^{\frac{n+3}{n-1}}, \qquad\quad \mathrm{~~in~~} \R^{n+1}
		\end{align*}
		is of the form
		$$u(X)=\left(\frac{\ve}{\ve^2+|X-X_0|^2}\right)^{\frac{n-1}{2}}, \qquad \mathrm{~~for~~some~~} \ve \in \R_+, X_0 \in \R^{n+1}.$$
	Under an additional hypothesis $u(X)=O(|X|^{1-n})$ as $|X| \to \infty$, the same result was also proved earlier by Obata \cite{Obata} and Gidas-Ni-Nirenberg \cite{GNN}. Escobar \cite{Escobar1} used  the Obata-type argument to characterize minimizers of a sharp Sobolev trace inequality on the unit ball  $\B^{n+1}$, via a conformal map between $\R_+^{n+1}:=\{X=(x,t)\in \Rn\times \R_+\}$ and $\B^{n+1}$, which is equivalent to seeking positive smooth solutions to
		\begin{align}\label{PDE:Li-Zhu_Liouville}
			\begin{cases}
			\displaystyle \Delta u=0,&\mathrm{~~in~~}\quad \R_{+}^{n+1},\\
			\displaystyle -\partial_tu(x,0)=(n-1)u^{\frac{n+1}{n-1}}\qquad\qquad&\mathrm{~~on~~}\quad  \partial\R^{n+1}_{+},
			\end{cases}
			\end{align}
			under a decay hypothesis that $u(X)=O(|X|^{1-n})$ as $|X| \to \infty$.
Y. Y. Li-Zhu \cite{Li-Zhu} (see also, Y. Y. Li-Zhang \cite{Li-Zhang}) employed  the moving sphere method to prove the following Liouville-type theorem without the decay hypothesis: Every positive  $C^2$ solution to \eqref{PDE:Li-Zhu_Liouville} is of the form 
		\begin{align*}
			u(x,t)=\left(\frac{\ve}{(\ve+t)^2+|x-x_0|^2}\right)^{\frac{n-1}{2}} , \qquad \mathrm{~~for~~some~~} \ve \in \R_+, x_0 \in \Rn.
		\end{align*}
 Such Liouville-type theorems of  critical  semilinear  second-order elliptic equations on $\R_+^{n+1}$ are referred to Escobar \cite{Escobar2}, Chipot-Shafrir-Fila \cite{CSF},  Y. Y. Li-Zhang \cite{Li-Zhang} and references therein.

Sun-Xiong \cite{Sun-Xiong} classified nonnegative smooth solutions of the following conformally invariant polyharmonic boundary value problem (BVP): For $m \geq 2,  n+1>2m$ and $1<p\leq \frac{n+2m-1}{n+1-2m}$,
	\begin{align}\label{bdry_value_problem:Xiong-Sun}
		\begin{cases}
	\displaystyle \Delta^m u=0\qquad\qquad\qquad\qquad\, &\mathrm{in}\qquad \R_{+}^{n+1},\\
	\displaystyle \partial_t\Delta^ku(x,0)=0\qquad\,\,\,\,\,\,\,&\mathrm{on}\quad\,\,\,\, \partial\R^{n+1}_{+},\\
	\displaystyle  (-1)^m\partial_t\Delta^{m-1}u(x,0)=u^{p}\qquad\,\,\,\,\,\,\,&\mathrm{on}\quad\,\,\,\, \partial\R^{n+1}_{+},
	\end{cases}
	\end{align}
	where $0 \leq k \leq m-2$. When $m$ is even, assume in addition that $u(X)=o(|X|^{2m-1})$ as $|X| \to \infty$. Then for the critical Sobolev exponent $p=\frac{n+2m-1}{n-(2m-1)}$,
	\begin{align*}
		u(x,t)=&c(n,m)\int_{\Rn}\frac{t^{2m-1}}{(|x-y|^2+t^2)^{\frac{n+2m-1}{2}}}\left(\frac{\lambda}{\lambda^2+|y|^2}\right)^{\frac{n-2m+1}{2}} \ud y+\sum_{k=1}^{m-1}t^{2k}P_{2m-2-2k}(x),
	\end{align*}	
	where  $x_0 \in \Rn, \lambda \in \R_+$, $c(n,m)$ is a normalized constant and $P_{2m-2-2k}(x)$ denotes a polynomial of degree $\leq 2m-2-2k$ satisfying $\liminf_{|x| \to \infty}P_{2m-2-2k}(x)\geq 0$; for $1<p<\frac{n+2m-1}{n+1-2m}$,
	\begin{align*}
		u(x,t)=\sum_{k=1}^{m-1}t^{2k}P_{2m-2-2k}(x),
	\end{align*}
	where $P_{2m-2-2k}$ is a nonnegative polynomial of degree $\leq 2m-2-2k$.
	One of crucial steps in \cite{Sun-Xiong}  is  to improve the integrability of a Kevin transformation $u^\ast$ for $(-\Delta)^m$ through an iteration method, and to show that the singularity of $u^\ast$ can be eventually removed at $t=0$. Next the integral equation of $u^\ast$ on $\pa \R_+^{n+1}$ is solved by the classification theorems due to Y. Y. Li \cite{Li} and Chen-Li-Ou \cite{Chen&Li&Ou}, then   the classification of  solutions to the BVP \eqref{bdry_value_problem:Xiong-Sun} follows by the the standard procedure.
	
	Conformal boundary operators $\mathscr{B}_k^3$ for $k \in \{1,2,3\}$ associated with the Paneitz operator  have attracted a lot of attention, see \cite{GS,BFGM,CLMNY,Chen-Zhang} etc. On a smooth four-manifold with boundary, Chang-Qing \cite{Chang-Qing} discovered a third-order conformal boundary operator $\mathscr{B}_3^3$ and its associated curvature $T_3^3$. To generalize the Chang-Qing boundary operator  to an $(n+1)$-manifold $(X^{n+1},g)$ with boundary $\pa X:=M^n$ for all conformal boundary operators $\mathscr{B}_k^3, k \in \{1,2,3\}$,  Branson-Gover \cite{Branson-Gover} initiated the Dirichlet-to-Neumann programme for all the higher order GJMS operators, and introduced the first- and second-order  conformal  boundary operators $\mathscr{B}_k^3$ associated with the Paneitz operator. Grant \cite{Grant} and Stafford \cite{Stafford} made suitable adjustments of the approach in \cite{Branson-Gover} to conduct  $\mathscr{B}_3^3$. Readers are referred  to  Juhl \cite{Juhl}, Gover-Peterson \cite{Gover-Peterson} for other treatments.  The Paneitz operator $P_4^g$ and $\mathscr{B}_k^3$ enjoy the property that if $\tilde g=e^{2\tau} g$, then for all $\Psi \in C^\infty(\overline{X^{n+1}})$,
\begin{equation}\label{conf_transformation:Paneitz}
\tilde P_4(\Psi)=e^{-\frac{n+5}{2}\tau}P_4^g(e^{\frac{n-3}{2}\tau} \Psi)
\end{equation}
and
\begin{equation}\label{conf_transformation:bdry_operators}
 \tilde{\mathscr{B}}_k^3 (\Psi)=e^{-\frac{n+2k-3}{2}\tau} \mathscr{B}_k^3(e^{\frac{n-3}{2}\tau} \Psi).
\end{equation}
 For future reference, adopting the notations as in \cite{Chang-Qing} and \cite{Gover-Peterson}, we call $(T_k)_g:=\frac{2}{n-3}\mathscr{B}_k^3(1), k \in \{1,2,3\}$ the $T$-curvatures on $M$, as the higher order generalization of the boundary mean curvature $h_g=\frac{2}{n-3}\mathscr{B}_1^3(1)$. Readers are referred to Case \cite{Case}  or the authors \cite{Chen-Zhang} for explicit formulae of $P_4^g$ and $\mathscr{B}_k^3$ on manifolds with boundary.

We are concerned with two models: $\R_+^{n+1}$ and $\B^{n+1}$. On  $(\R^{n+1}_{+}, \partial \R^{n+1}_{+}, |\ud X|^2)$, the explicit formulae of conformal boundary operators are  $\mathscr{B}^3_0(u)=u$ and
\begin{align*}
	\mathscr{B}^3_1(u)=-\partial_tu;\qquad
	\mathscr{B}^3_2(u)=\partial_t^2u-\overline{\Delta}u;\qquad
	\mathscr{B}^3_3(u)=&\partial_t\Delta u+2\overline{\Delta}\partial_tu
\end{align*}
for all $u \in C^\infty(\overline{\R_+^{n+1}})$. On $(\B^{n+1}, \Sn, |\ud \xi|^2)$, these conformal boundary operators are $\mathscr{B}^3_0(U)=U$ and
	\begin{align*}
		\mathscr{B}^3_1(U)=&\frac{\partial U}{\partial r}+\frac{n-3}{2}U,\\
			\mathscr{B}^3_2(U)=&\frac{\partial^2 U}{\partial r^2}-\Delta_{\S^n}U+(n-2)\frac{\partial U}{\partial r}+\frac{(n-3)(n-1)}{2}U,\\
			\mathscr{B}^3_3(U)=&-\frac{\partial \Delta U}{\partial r}-2\Delta_{\S^n}\frac{\partial U}{\partial r}-\frac{n-3}{2}\frac{\partial^2 U}{\partial r^2}-\frac{3n-5}{2}\Delta_{\S^n}U\\
			&+\frac{n-3}{2}\frac{\partial U}{\partial r}+\frac{(n^2-1)(n-3)}{4}U
	\end{align*}
	for all $U \in C^\infty(\overline{\B^{n+1}})$, where $r=|\xi|$ for $\xi \in \overline{\B^{n+1}}$;
 the corresponding $T$-curvatures $(T_i^3)_{\Sn}$ are constant and set $\mathbb{T}^3_i=\frac{n-3}{2}(T_i^3)_{\Sn}$ for $ i \in \{1,2,3\}$, explicitly,
\begin{align*}
	\mathbb{T}^3_1=\frac{n-3}{2}, \qquad \mathbb{T}^3_2=\frac{(n-1)(n-3)}{2} \qquad \mathrm{and}\qquad \mathbb{T}^3_3=\frac{(n^2-1)(n-3)}{4}
\end{align*}
will be used throughout the paper. 
\begin{convention} 
	 Let $\Delta$ and $\overline{\Delta}$ denote Euclidean Laplacians on $\R^{n+1}$ and $\R^n$, respectively.  Let $p_k^\ast=\frac{n+2k-3}{n-3}$  for $k \in \{1,2,3\}$, as $p_k^\ast+1$ is the critical Sobolev exponent associated to $\mathscr{B}^3_k$. 
			\end{convention}	

We are interested in the classification of nonnegative solutions to  the $Q$-flat and constant $T$-curvature type equations on $\R_+^{n+1}$: Let $n\geq 4$  and $p_i,p_j>0$ for $i,j \in \{1,2,3\}, i<j$, consider nonnegative smooth solutions to
\begin{align}\label{Intro main eqn}
		\begin{cases}
	\displaystyle ~~~ \Delta^2 u=0  &\mathrm{~~in~~}\quad~~ \R_+^{n+1},\\
	\displaystyle \mathscr{B}^3_i(u)=\mathbb{T}^3_i u^{p_i}\qquad \qquad&\mathrm{~~on~~}\quad \partial\R^{n+1}_{+},\\
	\displaystyle \mathscr{B}_j^3(u)=\mathbb{T}^3_ju^{p_j} &\mathrm{~~on~~}\quad \partial\R^{n+1}_{+},
	\end{cases}
\end{align} 
under hypothesis that $u(x,0)$ has some suitable decay at infinity.
The nonlinear interaction among $\mathscr{B}^3_i$ and $\mathscr{B}^3_j$ brings us new obstacles, compared with vanishing  of  lower order derivative boundary conditions in \cite{Sun-Xiong}. So, it has great difficulty in applying the approach in \cite{Sun-Xiong}  to \eqref{Intro main eqn}, and peculiar phenomena should occur simultaneously. 
 
   Although the Poisson kernel of a higher order elliptic BVP might appear elsewhere in the literature,  a rigorous definition of the \emph{biharmonic Poisson kernel} to the fourth-order elliptic BVP \eqref{Intro main eqn} is still in demand.
 \begin{definition}[Biharmonic Poisson kernel] 
For $i,j \in \{0,1,2,3\}, i<j$, a pair of kernel functions $P_i^{3} \oplus P_j^{3}(X)$ for $X \in \R_+^{n+1}$ is called a biharmonic Poisson kernel associated to a conformal  boundary operator pair $(\mathscr{B}_i^3,\mathscr{B}_j^3)$, if the following Poisson integral
\begin{align*}
v(X)=&[P_i^3 \ast f_i+P_j^3\ast f_j](X)\\
=&\int_{\pa\R_+^{n+1}}P_i^{3}(x-y,t)f_i(y) \ud y+\int_{\pa\R_+^{n+1}}P_j^{3}(x-y,t)f_j(y) \ud y
\end{align*}
makes sense, provided  $f_i, f_j \in C^2(\Rn)$ have suitable decays at infinity, and gives a classical solution to
\begin{align}\label{prob:bdry_value-Possion_kernel}
		\begin{cases}
	\displaystyle ~~~ \Delta^2 v=0  &\mathrm{~~in~~}\quad~~ \R_{+}^{n+1},\\
	\displaystyle \mathscr{B}^3_i(v)=f_i\qquad \qquad&\mathrm{~~on~~}\quad \partial\R^{n+1}_{+},\\
	\displaystyle \mathscr{B}_j^3(v)=f_j \qquad\,\,\,\,\,\,\,&\mathrm{~~on~~}\quad \partial\R^{n+1}_{+}.
	\end{cases}
\end{align} 
\end{definition}

We are now ready to give the explicit formulae of  biharmonic Poisson kernels.

\begin{thm} \label{Thm:Poisson kernel_half-space}
For $n \geq 2$ and $X=(x,t) \in \R_+^{n+1}$, we let\footnote{ A heuristic observation allows us to find these explicit kernel functions $P_k^3, k \in \{0,1,2,3\}$.}
\begin{align*}
P_0^3(X)=&\frac{2(n+1)}{|\S^n|}\frac{t^3}{(t^2+|x|^2)^{\frac{n+3}{2}}};\\
P_1^3(X)=&-\frac{2}{|\S^n|}\frac{t^2}{(t^2+|x|^2)^{\frac{n+1}{2}}};\\
P_2^3(X)=&-\frac{1}{(n-1)|\S^n|}\frac{t}{(t^2+|x|^2)^{\frac{n-1}{2}}};
\end{align*}
and
\begin{align*}
P_3^3(X)=\frac{1}{(n-1)(n-3)|\S^n|}\frac{1}{(t^2+|x|^2)^{\frac{n-3}{2}}}, \quad &\qquad \mathrm{~~when~~} n\neq 3;\\
P_3^3(X)=-\frac{1}{4|\S^3|}\log (t^2+|x|^2)+C, \quad C \in \R, &\qquad \mathrm{~~when~~} n= 3.
\end{align*}
Suppose (i) $f_0 \in C^2(\Rn)\cap L^\infty(\Rn)$ ; (ii) $ k \in \{1,2,3\}$, $f_k \in C^2(\Rn)$ and $f_k=O(|x|^{-k-\delta_k})$ as $|x| \to \infty$ for some $\delta_k \in \R_+$. Then for $i,j \in \{0,1,2,3\}, i+j \neq 3$,  the explicit formula of Poisson kernel corresponding to  problem \eqref{prob:bdry_value-Possion_kernel} is $P_i^{3} \oplus P_j^{3}$. 
\end{thm}
\begin{rem}
It is direct to verify that other than $(\mathscr{B}_0^3,\mathscr{B}_3^3)$ and $(\mathscr{B}_1^3,\mathscr{B}_2^3)$, the BVP $(\Delta^2; \mathscr{B}_i^3,\mathscr{B}_j^3)$ for problem \eqref{prob:bdry_value-Possion_kernel}  satisfies the \emph{Lopatinski-Shapiro condition} (e.g., Branson-Gover \cite[Section 6]{Branson-Gover}; or \cite[p.311]{Case}) or \emph{complementing condition} (e.g.,  Agmon-Douglis-Nirenberg \cite[Section 10]{ADN} or \cite[Section 2.3]{GGS}). 
\end{rem}

	 Let
	 \begin{align*}
	\Gamma(X-Y)=\Gamma(|X-Y|)=\begin{cases}
	\displaystyle \frac{1}{2(n-1)(n-3)|\S^n|}|X-Y|^{3-n} \qquad &\mathrm{~~when~~} n \neq 3,\\
	\displaystyle -\frac{1}{4|\S^3|}\log|X-Y|   &\mathrm{~~when~~} n=3.
	\end{cases}
	\end{align*}
	denote the fundamental solution of $\Delta^2$ in $\R^{n+1}$. 
	
		\begin{rem}
	As a traditional way of the classical Poisson kernel, conformal boundary operators $\mathscr{B}_k^3$ indeed act as a bridge between the biharmonic Poisson kernel and fundamental solution: Let $n \geq 2$, for every fixed $X=(x,t) \in \R_+^{n+1}$ and  all $Y=(y,0) \in \pa \R_+^{n+1}$, a direct computation yields
	\begin{align}\label{Poisson_Fund-sol}
	\begin{split}
\mathscr{B}^3_0(\Gamma(X-Y))=\frac{1}{2}P^3_3(x-y,t);&\qquad \mathscr{B}^3_1(\Gamma(X-Y))=\frac{1}{2}P^3_2(x-y,t);\\
	\mathscr{B}^3_2(\Gamma(X-Y))=-\frac{1}{2}P^3_1(x-y,t);&\qquad \mathscr{B}^3_3(\Gamma(X-Y))=-\frac{1}{2}P^3_0(x-y,t).
	\end{split}
\end{align}

	\end{rem}

Consider 
\begin{align}\label{BDVP:Green_fcn}
		\begin{cases}
	\displaystyle \quad \Delta^2 v=f  &\mathrm{~~in~~}\quad~~ \R_{+}^{n+1},\\
	\displaystyle \mathscr{B}^3_i(v)=0\qquad \qquad&\mathrm{~~on~~}\quad \partial\R^{n+1}_{+},\\
	\displaystyle \mathscr{B}_j^3(v)=0 \qquad\,\,\,\,\,\,\,&\mathrm{~~on~~}\quad \partial\R^{n+1}_{+}.
	\end{cases}
\end{align} 

\begin{definition}[Biharmonic Green function]
Let $n \geq 3$, $i,j \in \{0,1,2,3\},  i+j \neq 3$ and $f \in C^\infty(\R_+^{n+1})$. A function $G^{(i,j)}:(X,Y)\in \overline{\R_+^{n+1}}\times \overline{\R_+^{n+1}} \mapsto \R\cup \{\infty\}$ for \eqref{BDVP:Green_fcn}
	is called a Green function, if $G^{(i,j)}(X,Y)=\Gamma(X-Y)+H^{(i,j)}(X,Y)$, where $H^{(i,j)}$ is a  solution of
	\begin{align}\label{reg_part:Green-fcn}
\begin{cases}
\displaystyle ~~~ \Delta^2 H^{(i,j)}=0 \qquad &\mathrm{in}\qquad~ \R_+^{n+1},\\
\displaystyle \mathscr{B}^3_i(H^{(i,j)})= -\mathscr{B}^3_i(\Gamma) \qquad &\mathrm{on}\qquad \pa \R_+^{n+1},\\
\displaystyle \mathscr{B}^3_j(H^{(i,j)})=-\mathscr{B}^3_j(\Gamma)\qquad &\mathrm{on}\qquad \pa \R_+^{n+1}.
\end{cases}
\end{align}
	\end{definition}
	\begin{rem}
	As a direct consequence of Theorem \ref{Thm:Poisson kernel_half-space} and \eqref{Poisson_Fund-sol}, problem \eqref{reg_part:Green-fcn} admits a unique classical solution in $\R_+^{n+1}$ for $n \geq 3$.
	\end{rem}	
	
       More importantly, we can derive the explicit formulae of $H^{(i,j)}$ with the help of  the Boggio formula in Boggio \cite{Boggio} or Gazzola-Grunau-Sweers \cite[Lemma 2.27]{GGS}. 
       
       For clarity, we distinguish these representation formulae into two cases:  noncritical and critical dimensions.
	\begin{thm}\label{Green-fcn:noncritical}
	For $n \geq 2$ and $n \neq 3$, the biharmonic Green functions for $\R_+^{n+1}$ are
	\begin{align*}
	G^{(0,2)}(X,Y)=&\Gamma(X-Y)-\Gamma(\bar{X}-Y);\\
	G^{(1,3)}(X,Y)=&\Gamma(X-Y)+\Gamma(\bar{X}-Y);\\
	G^{(0,1)}(X,Y)=&\Gamma(X-Y)+\frac{1}{2}\Gamma(|\bar X-Y| )\left[(n-3)\frac{|X-Y|^2}{ |\bar X-Y|^2}-(n-1)\right];\\
	G^{(2,3)}(X,Y)=&\Gamma(X-Y)-\frac{1}{2}\Gamma(|\bar X-Y| )\left[(n-3)\frac{|X-Y|^2}{ |\bar X-Y|^2}-(n-1)\right];
	\end{align*}
	where $\bar X=(x,-t)$ for $X \in \R_+^{n+1}$.
	Moreover, there holds
$$0\leq G^{(0,1)}(X,Y)\leq G^{(0,2)}(X,Y)\leq G^{(1,3)}(X,Y)\leq G^{(2,3)}(X,Y), \quad \forall~ X,Y \in \overline{\R_+^{n+1}},X \neq Y,$$
with equality if and only if either $X \in \Rn$ or $Y \in \Rn$.
	\end{thm}
	
	\begin{thm}\label{Thm:Green fcn_half-space_dim-4}
	The representation formulae of Green functions for $\R_+^4$ are
	\begin{align*}
		G^{(0,1)}(X,Y)=&\Gamma(X-Y)-\Gamma(\bar X-Y)+\frac{1}{8|\S^3|}\left(\frac{|X-Y|^2}{|\bar X-Y|^2}-1\right);\\
		G^{(2,3)}(X,Y)=&\Gamma(X-Y)+\Gamma(\bar X-Y)-\frac{1}{8|\S^3|}\left(\frac{|X-Y|^2}{|\bar X-Y|^2}-1\right)+C;\\
			G^{(0,2)}(X,Y)=&\Gamma(X-Y)-\Gamma(\bar X-Y);\\
			G^{(1,3)}(X,Y)=&\Gamma(X-Y)+\Gamma(\bar X-Y)+C;
	\end{align*}
	for $X,Y \in \overline{\R_+^4}, X \neq Y$ and $C \in \R$.
\end{thm}

	Since $\R_+^{n+1}$ and $\B^{n+1}$ are conformally equivalent, we extend to introduce the notions of the biharmonic Poisson kernel $\bar P_k^3, k \in \{0,1,2,3\}$ and Green function $\bar G^{(i,j)}, i,j \in \{0,1,2,3\},i<j, i+j \neq 3$ for $\B^{n+1}$, and derive  their explicit formulae as well. 
	
	Until now, we have succeeded in obtaining a fourth-order analogue of the classical Poisson kernel and Green function for elliptic equations.
	
	For $\R_+^{n+1}$ and $\B^{n+1}$, we assemble the biharmonic Poisson kernel part in Section \ref{Sect3} and the biharmonic Green function part in Section \ref{Sect4}.

\medskip

To continue,  we introduce a \emph{`geometric bubble'}  by 
$$U_{x_0,\ve}(x,t)=\left(\frac{2\ve}{(\ve+t)^2+|x-x_0|^2}\right)^{\frac{n-3}{2}}, \qquad (x,t) \in \R_+^{n+1},$$
and its restriction to the boundary by letting $U_{x_0,\ve}(x)=U_{x_0,\ve}(x,0)$, where $\ve \in \R_+, x_0 \in \Rn$.

For $(\mathscr{B}^3_1, \mathscr{B}^3_3)$,  under appropriate decay rate of the initial data, the intrinsic  $\mathscr{B}^3_3$ equation with critical/subcritical growth enables us to determine the explicit  formula for $u$ on $\pa \R_+^{n+1}$, see Lemma \ref{Boundary data lem}.  With the aid of the basics of $U_{x_0,\ve}$ and biharmonic Poisson integral, we establish that every nonnegative smooth solution is a sum of a bubble and a polynomial. In sharp contrast to second-order elliptic equations, fix $p_3=p_3^\ast$, a \emph{striking} feature is that a nontrivial bubble (in general, without geometric meaning)  exists even for $p_1 \neq p_1^\ast$, and coincides with  $U_{x_0,\ve}$  for $p_1=p_1^\ast$; whereas for $p_1<p_1^\ast$, the solution is trivial.
	
		\begin{thm}\label{Main thm 1}
			For $n\geq 4$ and $p_1,p_3>0$, let $u$ be a nonnegative smooth solution of 
			\begin{align*}
			\begin{cases}
			\displaystyle ~~~\Delta^2 u=0  &\mathrm{~~in~~}\qquad~ \R_{+}^{n+1},\\
			\displaystyle \mathscr{B}_1^3(u)=\mathbb{T}^3_1u^{p_1}&\mathrm{~~on~~}\qquad \partial\R^{n+1}_{+},\\
			\displaystyle \mathscr{B}_3^3(u)=\mathbb{T}^3_3u^{p_3}\qquad\qquad &\mathrm{~~on~~}\qquad \partial\R^{n+1}_{+}.
			\end{cases}
			\end{align*}
			Assume that $u(x,0)=O(|x|^{-c})$ as $|x| \to \infty$ for some $c>\max\{3/p_3, 1/p_1\}$. Then we have
			\begin{enumerate}[(1)]
				\item	If $p_3=p_3^\ast=\frac{n+3}{n-3}$ and $p_1>\frac{1}{n-3}$, then either $u=0$ or there exist some $\ve \in \R_+, x_0\in \R^n$ and nonnegative constant $c_2$ such that
				\begin{align*}
				u(x,t)=M_1(x,t)+c_2t^2,
				\end{align*}
				where
				\begin{align*}
				M_1(x,t)=\left[\mathbb{T}_1^3 P_1^3 \ast U_{x_0,\ve}^{p_1}+\mathbb{T}_3^3 P_3^3 \ast U_{x_0,\ve}^{p_3^\ast}\right](x,t).
				\end{align*}
				In particular, when $p_1=p_1^\ast=\frac{n-1}{n-3}$, then  
								$$u(x,t)=U_{x_0,\ve}(x,t)+c_2 t^2.$$
				
				\item If $0<p_3<p_3^\ast$, then $u(x,t)=c_2t^2$ for every nonnegative constant $c_2$.	
			\end{enumerate}
		\end{thm}
		
		Although the  BVP $(\Delta^2; \mathscr{B}^3_1, \mathscr{B}^3_2)$ on $\R_+^{n+1}$ is not well-posed,  instead we take the derivative $\pa_t$ of  biharmonic Poisson integral of type $(\mathscr{B}^3_0, \mathscr{B}^3_2)$ to obtain the integral equation of $u^{p_1}(x,0)$, and the classification theorem by Y. Y. Li \cite{Li} and Chen-Li-Ou \cite{Chen&Li&Ou}) gives the explicit formula of $u(x,0)$. Our classification result  relies not on $p_1$ and $p_2$ but on the condition whether  $p_2/p_1$ is critical or subcritical Sobolev exponents.

			\begin{thm}\label{Main thm 2}
			For $n\geq 4$ and $p_1,p_2>0$, let $u$ be a nonnegative smooth solution of 
			\begin{align*}
			\begin{cases}
			\displaystyle ~~~\Delta^2 u=0, &\mathrm{~~in~~}\qquad~ \R_{+}^{n+1},\\
			\displaystyle \mathscr{B}^3_1(u)=\mathbb{T}^3_1u^{p_1}&\mathrm{~~on~~}\qquad\partial\R^{n+1}_{+},\\
			\displaystyle \mathscr{B}^3_2(u)=\mathbb{T}^3_2u^{p_2}\qquad\qquad&\mathrm{~~on~~}\qquad \partial\R^{n+1}_{+}.
			\end{cases}
			\end{align*}
			Assume that $u(x,0)=O(|x|^{-c})$ as $|x| \to \infty$ for some $c>2/p_2$. Then we have
				\begin{enumerate}[(1)]
				\item	If $\frac{p_2}{p_1}=\frac{n+1}{n-1}$ and $p_1>0$, then either $u=0$ or there exist some $\ve \in \R_+, x_0\in \Rn$ and nonnegative constant $c_3$ such that 
				\begin{align*}
				u(x,t)=M_2(x,t)+c_3t^3,
				\end{align*}
				where
				\begin{align*}
				M_2(x,t)=\Big[P_0^3 \ast U_{x_0,\ve}^{\frac{p_1^\ast}{p_1}}+\mathbb{T}_2^3 P_2^3 \ast U_{x_0,\ve}^{p_2^\ast}\Big](x,t).
				\end{align*}
				In particular, when $p_1=p_1^\ast=\frac{n-1}{n-3}$, then 			
								\begin{align*}
				u(x,t)=U_{x_0,\ve}(x,t)+c_3t^3.
				\end{align*}
				
				\item If $\frac{p_2}{p_1}<\frac{n+1}{n-1}$ , then $u(x,t)=c_3t^3$ for every nonnegative constant $c_3$.
			\end{enumerate}
		\end{thm}

		For $(\mathscr{B}^3_2, \mathscr{B}^3_3)$, the analysis of the polynomial part becomes more challenging, and requires  the techniques developed in Theorems \ref{Main thm 1} and \ref{Main thm 2}.
	
	       \begin{thm}\label{Main thm 3}
		For $n\geq 4$ and $p_2,p_3>0$, let $u$ be a nonnegative smooth solution of 
		\begin{align*}
		\begin{cases}
		\displaystyle ~~~\Delta^2 u=0  &\mathrm{~~in~~}\qquad~ \R_{+}^{n+1},\\
		\displaystyle \mathscr{B}^3_2(u)=\mathbb{T}^3_2u^{p_2} &\mathrm{~~on~~}\qquad \partial\R^{n+1}_{+},\\
		\displaystyle \mathscr{B}_3^3(u)=\mathbb{T}^3_3u^{p_3}\qquad\qquad &\mathrm{~~on~~}\qquad \partial\R^{n+1}_{+}.
		\end{cases}
		\end{align*}
		Assume that $u(x,0)=O(|x|^{-c})$  as $|X| \to \infty$ for some  $c>\max\{2/p_2,3/p_3\}$. Then we have
		\begin{enumerate}[(1)]
				\item	If $p_3=p_3^\ast=\frac{n+3}{n-3}$ and $p_2>\frac{2}{n-3}$, then either $u=0$ or there exist some $\ve \in \R_+, x_0\in \R^n$ and  nonnegative constant $c_1$ such that 
			\begin{align*}
			u(x,t)=M_3(x,t)+c_1t,
			\end{align*}
			where
			\begin{align*}
			M_3(x,t)=\left[P_0^3 \ast U_{x_0,\ve}+\mathbb{T}_2^3 P_2^3 \ast U_{x_0,\ve}^{p_2}\right](x,t).
			\end{align*}
			In particular, when $p_2=p_2^\ast=\frac{n+1}{n-3}$, then 
			\begin{align*}
			u(x,t)=U_{x_0,\ve}(x,t)+c_1t.
			\end{align*}
			
			\item If $0<p_3<\frac{n+3}{n-3}$, then $u(x,t)=c_1t$ for every nonnegative constant $c_1$.
		\end{enumerate}
	\end{thm}
	
In conclusion, $p_3$  is  the dominant exponent and the \emph{dominant impacts} of $\mathscr{B}^3_k$ are different, precisely, $\mathscr{B}^3_1\preceq \mathscr{B}^3_2\preceq \mathscr{B}^3_3$.  Especially, for  $p_i=p_i^\ast, i \in\{1,2,3\}$, the bubble happens to be the same  $U_{x_0,\ve}$; meanwhile, the corresponding polynomial varies according to various circumstances.  Theorems \ref{Main thm 1}-\ref{Main thm 3}  give a complete picture of Liouville-type theorems to the study of the  $Q$-flat and constant $T$-curvature problem on compact manifolds with boundary in the future.

In Section \ref{Sect2} we study the nontrivial solutions  of a homogeneous biharmonic BVP to \eqref{prob:bdry_value-Possion_kernel} (i.e., $f_i=f_j=0$) under suitable decay assumption of the initial data, and also describe asymptotic behavior of two typical singular integrals originating from the biharmonic Poisson integrals for $\R_+^{n+1}$. In Section \ref{Sect5} we present the proof of classification Theorems \ref{Main thm 1}-\ref{Main thm 3}. 

	\medskip

		Theorems \ref{Main thm 1}-\ref{Main thm 3} enable us to classify nonnegative solutions to the $Q$-flat and constant $T$-curvature equations  on $\B^{n+1}$ with one boundary point singularity.

	\begin{thm}\label{Thm:bdry_singularity}
	Suppose $n\geq 4$ and $U \in C^4(\overline{\B^{n+1}}\backslash \{-\mathbf{e}_{n+1}\})$ is a nonnegative solution to
	\begin{align}\label{bdry_value_problem_singular pt}
	\begin{cases}
	\displaystyle ~~~\Delta^2 U=0  &\mathrm{~~in~~}\qquad \B^{n+1},\\
	\displaystyle \mathscr{B}^3_i(U)=\mathbb{T}^3_i U^{\frac{n+2i-3}{n-3}} &\mathrm{~~on~~}\qquad \Sn\backslash\{-\mathbf{e}_{n+1}\},\\
	\displaystyle \mathscr{B}^3_j(U)=\mathbb{T}^3_j U^{\frac{n+2j-3}{n-3}}\qquad\qquad &\mathrm{~~on~~}\qquad\Sn\backslash\{-\mathbf{e}_{n+1}\},
	\end{cases}
	\end{align}
	where $i,j \in \{1,2,3\},i<j$. Assume that there exist $\alpha\in (0,\frac{n(n-3)}{n+3})$ and $C \in \R_+$ such that 
	\begin{align}\label{assump:isolated_singularity}
	\sup_{\xi\in\S^{n}\backslash\{-\mathbf{e}_{n+1}\}}U(\xi)|\xi+\mathbf{e}_{n+1}|^{\alpha}\leq C .
	\end{align}
	Then there exist some $\xi_{0}\in \B^{n+1}$ and nonnegative constant $\bar c$ such that 
	\begin{align*}
	U(\xi)=\left(\frac{1-|\xi_0|^2}{|\xi|^2 |\xi_0|^2-2 \xi_0 \cdot \xi+1}\right)^{\frac{n-3}{2}}
	+\bar c |\xi+\mathbf{e}_{n+1}|^{3-n}\left(\frac{1-|\xi|^2}{|\xi+\mathbf{e}_{n+1}|}\right)^{6-i-j}.
	\end{align*}
	\end{thm}
	
	\begin{rem}\label{Intro rem 1}
		If $i=1,j=2$, then the assumption of $\alpha$ can be relaxed to $\alpha\in (0,\frac{(n-1)(n-3)}{n+1})$, which ensures that  $n-3-\alpha>\frac{2}{p_2^\ast}$ by virtue of Theorem \ref{Main thm 2}. 
	\end{rem}	
	
	We would digress for the moment to comment on  the decay assumption of initial data  that $u(x,0)=O(|x|^{-c})$ as $|x| \to \infty$ in Theorems \ref{Main thm 1}-\ref{Main thm 3}. For critical exponents $p_i^\ast$, such assumptions are to ensure that involved biharmonic Poisson integrals are well-defined. Consequently, this leads to $c>\frac{3(n-3)}{n+3}$ in Theorems \ref{Main thm 1} ,\ref{Main thm 3} and $c>\frac{2(n-3)}{n+1}$ in Theorem \ref{Main thm 2}. As will show in the proof of Theorem \ref{Thm:bdry_singularity},  we indeed have $c=n-3-\alpha$, so  \eqref{assump:isolated_singularity}  gives   $\alpha<\frac{n(n-3)}{n+1}$ and $\alpha<\frac{(n-1)(n-3)}{n+1}$, respectively. Besides this, if $\mathscr{B}^3_3(U)=\frac{n-3}{2}\tilde T^3_3 U^{\frac{n+3}{n-3}}\in L^1(\S^n)$ for a conformal metric $\tilde g=U^{4/(n-3)}|\ud \xi|^2$, then near  $-\mathbf{e}_{n+1}$ we have
		\begin{align*}
		U^{\frac{n+3}{n-3}}(\xi)\lesssim |\xi+\mathbf{e}_{n+1}|^{-\frac{\alpha(n+3)}{n-3}}\in L^1(\S^n),
		\end{align*}
		which requires  $0<\alpha<\frac{n(n-3)}{n+3}$. Intuitively, $\mathscr{B}^3_3(U)\in L^1(\S^n)$ means $U\in W^{3,1}(\S^n)\hookrightarrow W^{2,\frac{n}{n-1}}(\S^n)$. This roughly yields  $\mathscr{B}^3_2(U)\in L^{\frac{n}{n-1}}(\S^n)$, and also
		\begin{align*}
		\left(U^{\frac{n+1}{n-3}}(\xi)\right)^{\frac{n}{n-1}}\lesssim|\xi+\mathbf{e}_{n+1}|^{\frac{\alpha n(n+1)}{(n-1)(n-3)}}\in L^1(\S^n),
		\end{align*}
		which requires $\alpha<\frac{(n-1)(n-3)}{n+1}$. 
	
	\begin{cor}\label{Geom_Rigidity}
	Suppose $n\geq 4$ and $U$ is a positive smooth solution in $\overline{\B^{n+1}}$ to
	\begin{align*}\label{Introduction geometry equ}
	\begin{cases}
	\displaystyle ~~~\Delta^2 U=0 &\mathrm{~~in~~}\qquad \B^{n+1},\\
	\displaystyle \mathscr{B}^3_i(U)=\mathbb{T}^3_i U^{\frac{n+2i-3}{n-3}} &\mathrm{~~on~~}\qquad \Sn,\\
	\displaystyle \mathscr{B}^3_j(U)=\mathbb{T}^3_j U^{\frac{n+2j-3}{n-3}} \qquad \qquad &\mathrm{~~on~~}\qquad \Sn,
	\end{cases}
	\end{align*}
	where $i,j \in \{1,2,3\},i<j$. Then
	$$U(\xi)=\left(\frac{1-|\xi_0|^2}{|\xi|^2 |\xi_0|^2-2 \xi_0 \cdot \xi+1}\right)^{\frac{n-3}{2}}$$
	for some $\xi_0 \in \B^{n+1}$.
	\end{cor}
	
	A geometric interpretation of Corollary \ref{Geom_Rigidity} is that, \emph{Let $g$ be a $Q$-flat conformal metric on $(\B^{n+1},|\ud \xi|^2)$ for $n \geq 4$. Assume that  $\big((T^3_i)_g, (T^3_j)_g\big)=\big((T^3_i)_{\Sn},(T^3_j)_{\Sn}\big)$ for $i,j\in \{1,2,3\},i<j$. Then there exists a conformal diffeomorphism $\Psi$ on $\overline{\B^{n+1}}$ such that $g= \Psi^\ast (|\ud \xi|^2)$.}
	
	\medskip

		In Section \ref{Sect6}, we give the proof of Theorem \ref{Thm:bdry_singularity}. As applications we establish two sharp geometric inequalities over certain $Q$-flat conformal class of $(\B^{n+1},|\ud \xi|^2)$ for $n\geq 4$, including a sharp fourth-order  isoperimetric  inequality.
		
		In Section \ref{Sect7}, we raise a rigidity conjecture for a geodesic ball in $\S^{n+1}$, which is a natural extension of Corollary \ref{Geom_Rigidity}, and could be reduced to the classification of positive solutions to the constant $Q$-curvature and constant $T$-curvature equation in $\B^{n+1}$. We prove a uniqueness result of positive solutions to a  fourth-order geometric ODE in Appendix \ref{Append:B}. This allows us to answer this conjecture affirmatively for radial solutions. In Appendix \ref{Append:A}, we give an alternative proof of Proposition \ref{Prop:ex-intrinsic_GJMS}.
		
		\bigskip

\noindent{\bf Acknowledgments:} The authors thank Professor  Qianqiao Guo for bringing  \cite{Gluck} into their sight.

					\section{Preliminaries}	\label{Sect2}

		 This section is devote to building some technical lemmas and  giving asymptotic behavior of two prototypes of singular integrals, which are two variants of biharmonic Poisson integrals  for $\R_+^{n+1}$.
	
	\begin{lem}\label{singular integral estimate lem 1}
		Suppose $f\in L_{\mathrm{loc}}^{1}(\R^n)$ with  $f(x)=O(|x|^{-a})$ as $|x| \to \infty$ for some $a \in \R_+$. Assume that $0<\alpha<n/2$ and $a+2\alpha>n$.  We define
		\begin{align*}
			v(x,t)=\int_{\R^n}\frac{f(y)}{(t^2+|x-y|^2)^{\alpha}} \ud y.
		\end{align*} 
		Then  there holds 
		\begin{align*}
		\left|v(X)\right|\lesssim& \begin{cases}
	|X|^{n-2\alpha-\min\{a,n\}}  \quad &\mathrm{if~~} a \neq n\\
	|X|^{-2\alpha} \log |X| &\mathrm{if~~} a=n
	\end{cases} \qquad \quad \mathrm{~~as~~} \quad |X| \to \infty.
		\end{align*}
		
	\end{lem}
	\begin{pf}
			We decompose $\R^n=A_1\cup A_2\cup A_3$ with 
		\begin{align*}
		A_1=&\left\{y; |y|<|x|/2\right\},\qquad \qquad \qquad\qquad  A_2=\left\{y; |x-y|<|x|/2\right\},\\
		A_3=&\left\{y; |y|\geq |x|/2,|x-y|\geq |x|/2\right\}.
		\end{align*}
		
			\medskip
    \noindent\underline{\emph{Case 1.}}  $|X| \gg 1$ and $t>\frac{|x|}{2}$.
    \medskip
    
    There hold $t\gg 1$ and $|X| \sim t$ due to
		\begin{align*}
		t^2 \leq 	|X|^2=t^2+|x|^2\leq 5 t^2.
		\end{align*}
		
		In $A_1$,  we have $|x-y|\sim |x|$ and $ t^{2}+|x-y|^2\sim t^2$, hence
	\begin{align*}
	|v(x,t)|\lesssim t^{-2\alpha}\int_{A_1}|f(y)|\ud y	\lesssim \begin{cases}
	|X|^{n-2\alpha-\min\{a,n\}} &\mathrm{if~~} a \neq n\\
	|X|^{-2\alpha} \log |X| &\mathrm{if~~} a=n
	\end{cases}.
	\end{align*}
	
	In $A_2$, we have $|x-y|<|x|/2<t$, $|y|\sim |x|$ and  $ t^{2}+|x-y|^2\sim t^2$, then 
	\begin{align*}
	|v(x,t)|\lesssim t^{-2\alpha}\int_{A_2}|f(y)|\ud y\lesssim& \begin{cases}
	|X|^{n-2\alpha-\min\{a,n\}} &\mathrm{if~~} a \neq n\\
	|X|^{-2\alpha} \log |X| &\mathrm{if~~} a=n
	\end{cases}.
	\end{align*}
	
	In $A_3$, we have $|x-y|\leq |y|+|x|\leq3|y|$ and $|y|\leq|x-y|+|x|\leq 3|x-y|$, whence $|x-y|\sim|y|$ . We further decompose $A_3=A_{31}\cup A_{32}$ with
	\begin{align*}
	A_{31}=A_3\cap \{y; |y|>t\}, \qquad\qquad A_{32}=A_3\cap \{y;|y|\leq t\}.
	\end{align*}
	In $A_{31}$, we have $t^2+|x-y|^2\sim |y|^2$ and 
	\begin{align*}
	|v(x,t)|\lesssim&\int_{A_{31}}\frac{|f(y)|}{|y|^{2\alpha}}\ud y
	\lesssim \int_{|y|>t}|y|^{-a-2\alpha}\ud y\lesssim t^{n-a-2\alpha}\lesssim |X|^{n-a-2\alpha}.
	\end{align*}
	In $A_{32}$, we have $t^2+|x-y|^2\sim t^2$ and
	\begin{align*}
	|v(x,t)|\lesssim t^{-2\alpha}\int_{A_{32}}|f(y)|\ud y\lesssim& \begin{cases}
	t^{n-2\alpha-\min\{a,n\}} &\mathrm{if~~} a \neq n\\
	t^{-2\alpha} \log t &\mathrm{if~~} a=n
	\end{cases}\\
	\lesssim& \begin{cases}
	|X|^{n-2\alpha-\min\{a,n\}}  &\mathrm{if~~} a \neq n\\
	|X|^{-2\alpha} \log |X| &\mathrm{if~~} a=n
	\end{cases}.
	\end{align*}
	
	Hence we combine the estimates above to obtain
	\begin{align*}
	\left|v(x,t)\right|\lesssim  \begin{cases}
	|X|^{n-2\alpha-\min\{a,n\}} \qquad &\mathrm{if~~} a \neq n\\
	|X|^{-2\alpha} \log |X| &\mathrm{if~~} a=n
	\end{cases}.
	\end{align*}
\medskip
    \noindent\underline{\emph{Case 2.}} $|X|\gg 1$ and  $t\leq \frac{|x|}{2}$. 
    \medskip
    
      We have $|x|\gg 1$ and $|X|\sim |x|$ due to
   \begin{align*}
		|x|^2	\leq 	|X|^2=t^2+|x|^2\leq \frac{5}{4} |x|^2.
		\end{align*}
    
	In $A_1$,  we have $\frac{|x|}{2}\leq|x-y|\leq \frac{3|x|}{2}$ and  $t^2+|x-y|^2\sim |x|^2$, then
	\begin{align*}\label{dim 3 singular integral estimate a_1}
	|v(x,t)|\lesssim |x|^{-2\alpha}\int_{A_1}|f(y)|\ud y\lesssim&	\begin{cases}
	|x|^{-2\alpha}|x|^{n-\min\{a,n\}} &\mathrm{if~~} a \neq n\\
	|x|^{-2\alpha} \log (|x|+1) &\mathrm{if~~} a=n
	\end{cases}\no\\
	\lesssim& \begin{cases}
	|X|^{n-2\alpha-\min\{a,n\}} &\mathrm{if~~} a \neq n\\
	|X|^{-2\alpha} \log |X| &\mathrm{if~~} a=n
	\end{cases}.
	\end{align*}
	
	In $A_2$, we have $|x|/2<|y|<3|x|/2$. Further decompose $A_2=A_{21}\cup A_{22}$ with
	\begin{align*}
	A_{21}=A_2\cap \{y; |x-y|>t\} \qquad \mathrm{and}\qquad A_{22}=A_2\cap \{y; |x-y|\leq t\}.
	\end{align*}
	For $A_{21}$,  we have $t^2+|x-y|^2\sim |x-y|^2$ and 
	\begin{align*}
	|v(x,t)|\lesssim \int_{A_{21}}\frac{|f(y)|}{|x-y|^{2\alpha}}\ud y\lesssim |x|^{-a}\int_{|x-y|<|x|/2}\frac{1}{|x-y|^{2\alpha}}\ud y\lesssim |X|^{n-2\alpha-a}.
	\end{align*}
	For $A_{22}$,  $t^2+|x-y|^2\sim t^2$, we know
	\begin{align*}
	|v(x,t)|
	\lesssim t^{-2\alpha}\int_{A_{22}}|f(y)|\ud y
	\lesssim t^{-2\alpha}|x|^{-a}\int_{|x-y|\leq t}\ud y\lesssim |X|^{n-2\alpha-a}.
	\end{align*}
	For $A_3$, we know $|y|/3<|x-y|<3|y|$, $|y|>|x|/2>t$, whence, $t^2+|x-y|^2\sim |y|^2$. Then,
	\begin{align*}
	|v(x,t)|
	\lesssim\int_{A_{3}}\frac{|f(y)|}{|y|^{2\alpha}}\ud y
	\lesssim\int_{|y|>|x|/2}|y|^{-(2\alpha+a)}\ud y\lesssim |X|^{n-2\alpha-a}.
	\end{align*}
	
	Hence, putting these facts above together we have
	\begin{align*}
	\left|v(x,t)\right|\lesssim  \begin{cases}
	|X|^{n-2\alpha-\min\{a,n\}} \qquad &\mathrm{if~~} a \neq n\\
	|X|^{-2\alpha} \log |X| &\mathrm{if~~} a=n
	\end{cases}.
	\end{align*}
	
	\end{pf}

\begin{lem}\label{singular integral estimate lem 2}

Suppose $f\in L_{\mathrm{loc}}^{1}(\R^n)$ satisfies that  $f(x)=O(|x|^{-b})$ as $|x| \to \infty$ for some $b \in \R_+$. Let
		\begin{align*}
			w(x,t)=\int_{\R^n}\frac{t^{\beta}  f(y)}{(t^2+|x-y|^2)^{\frac{n+\beta}{2}}} \ud y, \quad \beta \in \R_+.
		\end{align*} 
		Then  there holds 
	\begin{align*}
	\left|w(x,t)\right|\lesssim& \begin{cases}
	|X|^{-\min\{b,n\}} \qquad &\mathrm{if~~} b \neq n\\
	|X|^{-n} \log |X| &\mathrm{if~~} b=n
	\end{cases} \qquad \quad \mathrm{~~as~~} \quad |X| \to \infty.
	\end{align*}
	
\end{lem}

Since the same trick works for the above second type singular integral, we omit the proof.

\medskip

	 For $r \in \R_+$, let 
	$B_r^+=B_r(0)\cap \R^{n+1}_{+}$ and $\partial^{+}B_r^+=\partial B_r^+\cap\R^{n+1}_{+}$. Let $P_k(x)$ denote a polynomial with degree $\leq k$, where $k \in \N$. Given $u \in C^{\infty}(\overline{\R^{n+1}_+})$, we introduce a Kelvin transform for $\Delta^2$ by
		$$u^\ast(X)=|X|^{3-n}u\left(X|X|^{-2}\right), \qquad X \in \overline{\R_+^{n+1}}\backslash \{0\}.$$

	For our purpose, we need to a refinement of Sun-Xiong's results \cite[Proposition 3.1 and Theorem 3.2]{Sun-Xiong} for homogeneous biharmonic BVPs on $\R_+^{n+1}$.
	\begin{lem}\label{Uniqueness lem}
		Let $u\in C^{\infty}(\overline{\R^{n+1}_+})$ be a solution of
			\begin{align}\label{homo bdry value problem I}
		\begin{cases}
		\displaystyle ~\Delta^2 u=0 &\mathrm{~~in~~}\quad~ \R_{+}^{n+1},\\
		\displaystyle 	~~~ \partial_{t}u=0 &\mathrm{~~on~~}\quad \partial\R^{n+1}_{+},\\
		\displaystyle\partial_{t} \Delta u=0\qquad \qquad &\mathrm{~~on~~}\quad \partial\R^{n+1}_{+}.
		\end{cases}
		\end{align}
		Suppose    (i) $\lim_{|x| \to \infty} u(x,0)=0$; (ii) there exist some $\delta, C \in \R_+$ such that
		\begin{align}\label{Uniqueness thm condition}
u(X)\geq -C|X|^{-\delta}  \qquad \mathrm{as} \quad |X| \to \infty.
		\end{align}
	Then, $u(x,t)=c_2t^2$ for some nonnegative constant $c_2$.
	\end{lem}
	
For our purpose,  the nonnegativity assumption of solutions in Sun-Xiong \cite{Sun-Xiong}  is redundant. 
	\begin{pf}
	The proof consists of three steps.

			\medskip
    \noindent\underline{\emph{Step 1.}} $u^\ast \in L^1(B_1^+)$.
                          \medskip
		
		First, rewrite the boundary conditions as $\mathscr{B}_1^3 u=0, \mathscr{B}_3^3 u=0$. Then the conformal covariance of $\mathscr{B}_1^3$ and $\mathscr{B}_3^3$ yields that $u^\ast$ satisfies  problem \eqref{homo bdry value problem I} in $\overline{\R_+^{n+1}}\backslash \{0\}$. 
				
 		Define
		\begin{align*}
			\eta_{\ve}(t)=\begin{cases}
			\displaystyle \frac{1}{24}(t-\ve)^4 \qquad&\mathrm{for}\qquad t\geq \ve,\\
			\displaystyle 0 \qquad&\mathrm{for}\qquad 0\leq t< \ve.
			\end{cases}
		\end{align*}
		Notice that $\eta_{\ve} \in C^{3,1}(\R^{n+1}_{+})$ and $u^{*}$ is biharmonic on $\overline{\R^{n+1}_+}\backslash\{0\}$. Then integrating by parts gives
		\begin{align*}
			&\int_{B_1^+}u^{*}(X)\Delta^2\eta_{\ve} \ud X\\
						=&\int_{\partial^{+}B_1^+}\left(\Delta u^{*}(X)\frac{\partial}{\partial r}\eta_{\ve}-\frac{\partial}{\partial r}\Delta u^{*}(X)\eta_{\ve}+u^{*}(X)\frac{\partial}{\partial r}\Delta\eta_{\ve}-\frac{\partial}{\partial r}u^{*}(X)\Delta\eta_{\ve}\right)\ud \sigma.
		\end{align*}
		Letting $\ve\to 0$ we obtain 
		\begin{align}\label{Uniqueness lem formula a}
			\int_{B_1^+}u^{*}(X)\ud X\leq C.
		\end{align}
		
		Write $u^{*}(X)=u_{+}^{*}(X)-u_{-}^{*}(X)$. It follows from the assumption \eqref{Uniqueness thm condition} that $u_{-}^{*}(X)\leq C|X|^{\delta+3-n}$. Using \eqref{Uniqueness lem formula a} we obtain
		 \begin{align*}
			\int_{B_1^+}u_{+}^{*}(X)\ud X\leq C+\int_{B_1^+}u_{-}^{*}(X)\ud X\leq C
		\end{align*}
		and thus $u^\ast \in L^1(B_1^+)$.
		
		\medskip
    \noindent\underline{\emph{Step 2.}} $u(X)=o(|X|^4)$ as $|X| \to \infty$. 
                \medskip
		
		For $r \in \R_+$, we set $v_r(X)=u^{*}(rX)$, which is biharmonic in $\overline{\R^{n+1}_+}\backslash\{0\}$. Now we extend $v_r$ by an even reflection to $\R^{n+1}\backslash \{0\}$  and apply the local estimates for biharmonic functions (see e.g.,  \cite[Proposition 4]{Martinazzi}) to obtain
		\begin{align*}
			\|v_r\|_{L^{\infty}(B^{+}_{3/4}\backslash B^{+}_{1/2})}\lesssim \|v_r\|_{L^{1}(B^{+}_{5/4}\backslash B^{+}_{1/4})}\lesssim r^{-(n+1)}\|u^{*}\|_{L^{1}(B^{+}_{5r/4}\backslash B^{+}_{r/4})}.
		\end{align*}
		Since $u^\ast \in L^1(B_1^+)$, we may take $r$ sufficiently small such that
		\begin{align*}
			\|u^{*}\|_{L^{\infty}(B^{+}_{3r/4}\backslash B^{+}_{r/2})}=o(r^{-(n+1)}).
		\end{align*}
This in turn implies 
	\begin{align*}
		u(X)=o(|X|^{4}) \qquad \mathrm{~~as~~} |X| \to \infty.
	\end{align*}	
	
	\medskip
    \noindent\underline{\emph{Step 3.}} Completion of the proof.  
	\medskip	
		
		Extend $u$ by an even reflection
		\begin{align*}
		\tilde{u}(x,t)=\begin{cases}
		\displaystyle u(x,t)\qquad&\mathrm{for}\qquad t\geq 0,\\
		\displaystyle u(x,-t)\qquad&\mathrm{for}\qquad t<0.
		\end{cases}
		\end{align*}
		Notice that $\partial_{t} \tilde{u}(x,0)=0$ and $$\partial_{t}^3u(x,0)=\partial_{t}\Delta u(x,0)-\partial_{t}\overline{\Delta} u(x,0)=0,$$
		so $\tilde{u} \in C^4(\R^{n+1})$ and is biharmonic in $\R^{n+1}$. By \cite[Proposition 4]{Martinazzi} we obtain 
		\begin{align*}
		\|\nabla^{4} \tilde{u}\|_{L^{\infty}(B_R)}\leq \frac{C}{R^{4}}\frac{1}{|B_{2R}|}\int_{B_{2R}}|\tilde{u}| \ud X\leq \frac{o(R^{4})}{R^{4}}\to 0 \qquad\mathrm{as}\qquad R\to +\infty.
		\end{align*}
		Thus we have
		\begin{align*}
		u(x,t)=\sum_{i=0}^{3}t^{i}P_{3-i}(x).
		\end{align*} 
	The homogeneous boundary conditions on $\pa \R_+^{n+1}$ give
		\begin{align*}
			0=\partial_tu(x,0)=\sum_{i=1}^{3}i t^{i-1}P_{3-i}(x)=P_2(x)
		\end{align*}
		and
		\begin{align*}
		0=\partial^3_tu(x,0)=\sum_{i=1}^{3}i(i-1)(i-2)t^{i-3}P_{3-i}(x)=6P_0(x).
		\end{align*}
		This yields
		\begin{align*}
		u(x,t) =P_3(x)+t^2P_1(x).
		\end{align*}
		Since $\liminf_{|X|\to+\infty}u(X)\geq 0$ by virtue of \eqref{Uniqueness thm condition} and $u(x,0)=o(1)$ as $|x| \to \infty$, we conclude that $P_3(x)=0$ and  $P_1(x)=c_2$ for some nonnegative constant $c_2$.
	Hence, the desired assertion follows.
	\end{pf}

\begin{rem}\label{Rem:Uniqueness}
	If we relax the above  \eqref{Uniqueness thm condition} to  the assumption that $u(X)\geq -C|X|^{\delta}$ for some  constant $0<\delta<4$, which ensures that $u^\ast_- \in L^1(B_1^+)$, then the above argument yields
	\begin{align*}
		u(x,t)=t^2P_1(x).
	\end{align*}
\end{rem}

	\begin{lem}\label{Uniqueness lem 2}
		Let $u\in C^{\infty}(\overline{\R^{n+1}_+})$ be a solution of
		\begin{align}\label{homo bdry value problem II}
		\begin{cases}
		\displaystyle \Delta^2 u=0 &\mathrm{~~in~~}\qquad~ \R_{+}^{n+1},\\
		\displaystyle \quad~ u=0 &\mathrm{~~on~~}\qquad \partial\R^{n+1}_{+},\\
		\displaystyle ~~ \Delta u=0\qquad \qquad&\mathrm{~~on~~}\qquad\partial\R^{n+1}_{+}.
		\end{cases}
		\end{align}
		Suppose  there exist some $\delta, C \in \R_+$ such that
		\begin{align}\label{Cond:Uniqueness_weak}
u(X)\geq -C|X|^{-\delta}  \quad \mathrm{as} \quad |X| \to \infty.
		\end{align}
		Then, 
		\begin{align*}
		u(x,t)=tP_2(x)+c_3t^3,
		\end{align*}
		where $c_3$ is a nonnegative constant and
			\begin{align*}
		\liminf_{|x|\to\infty} P_2(x)\geq 0.
		\end{align*}
	\end{lem}
	\begin{pf}
		Notice that  these boundary conditions are equivalent to $\mathscr{B}_0^3 u=0, \mathscr{B}_3^3 u=0$ on $\pa \R_+^{n+1}$. By conformal covariance,  $u^\ast$  satisfies  problem \eqref{homo bdry value problem II} on $\overline{\R^{n+1}_+}\backslash\{0\}$. Following the same lines as before, we have $u^\ast \in L^1(B_1^+)$ and $u(x,t)=o(|X|^4)$ as $|X| \to \infty$. Finally, we extend $u$ to $\R^{n+1}$ by an odd reflection
		\begin{align*}
		\tilde{u}(x,t)=\begin{cases}
		\displaystyle u(x,t)\qquad&\mathrm{~~for~~}\quad t\geq 0,\\
		\displaystyle -u(x,-t)\qquad&\mathrm{~~for~~}\quad t<0.
		\end{cases}
		\end{align*}
		A similar argument as in Lemma \ref{Uniqueness lem} combined with the boundary conditions shows
		\begin{align*}
			u(x,t)=tP_2(x)+c_3t^3.
		\end{align*}
		This together with the decay assumption at infinity yields the desired assertion.
	\end{pf}

	\begin{rem}\label{Uniqueness lem 2 rem 2}
		If we relax  the above assumption \eqref{Cond:Uniqueness_weak} to $u(X)\geq -C|X|^{\delta}$ as $|X| \to \infty$ for some $0<\delta<3$, then the same result in Lemma \ref{Uniqueness lem 2} is true.
			\end{rem}

	\section{Biharmonic Poisson kernel}\label{Sect3}

	For a smooth compact Riemannian manifold $(X^{n+1}, g)$ with boundary $M^n=\pa X$ and $n\geq 4$, under the assumption that a certain conformal invariant $\lambda_1(P_4^g)>0$ (see \cite{Case}),  the existence and uniqueness of smooth solutions to  the following BVP was studied in \cite[Proposition 1.3]{Case}: For  $\bar f_0,\bar f_1 \in C^\infty(M)$,
	\begin{align*}
	\begin{cases}
	\displaystyle ~P_4^g U=0 &\mathrm{~~in~~}\qquad~ X,\\
	\displaystyle \mathscr{B}^3_0 U=\bar f_0 &\mathrm{~~on~~}\qquad M,\\
	\displaystyle \mathscr{B}^3_1U=\bar f_1 \qquad\qquad &\mathrm{~~on~~}\qquad  M.
	\end{cases}
	\end{align*} 
	
		Our purpose is to derive explicit representation formulae of biharmonic Poisson kernels for $\R_+^{n+1}$ and $\B^{n+1}$.
On  $(\B^{n+1}, \Sn, |\ud \xi|^2)$, for each proper conformal boundary operator pair, such a biharmonic BVP is expected to be better understood,  serving as a first step to the prescribing boundary curvature problem on manifolds with boundary. 
	
	\subsection{The half-space}

	We start with the classical Poisson kernel for $-\Delta$ in $\R_+^{n+1}$:
$$P(x,t)=\frac{2}{|\Sn|}\frac{t}{(t^2+|x|^2)^{\frac{n+1}{2}}}.$$
Moreover, $P$ enjoys the property that let $f \in C(\Rn) \cap L^\infty(\Rn)$, then 
\emph{(a)} the Poisson integral $P\ast f \in C^\infty(\R_+^{n+1})\cap L^\infty(\R_+^{n+1})$ and $\Delta (P\ast f)=0$ in $\R_+^{n+1}$; \emph{(b)} for each fixed $x \in \Rn$,  $\lim_{t\to 0} P \ast f(x,t)=f(x)$.

\begin{lem}\label{lem:singular_integral_II}
	Suppose $f\in C^2(\R^n)\cap L^\infty(\Rn)$ and $\alpha>1$. Then for every fixed $x \in \Rn$, we have
	\begin{align*}
		\lim_{t\searrow 0}\partial_t\int_{\R^n}\frac{t^{\alpha}}{\left(t^2+|x-y|^2\right)^{\frac{n+\alpha}{2}}}f(y)\ud y=0.
	\end{align*}
 \end{lem}
 
 It is not hard to see that the above result is still true provided both $f$ and $|\nabla f|$  are bounded. This elementary lemma will be used later.

\begin{pf}
	A direct calculation yields
	\begin{align*}
		&\partial_t\int_{\R^n}\frac{t^{\alpha}}{\left(t^2+|x-y|^2\right)^{\frac{n+\alpha}{2}}}f(y)\ud y\\
		=&\int_{\R^n}\left[\frac{\alpha t^{\alpha-1}}{{\left(t^2+|x-y|^2\right)^{\frac{n+\alpha}{2}}}}-\frac{(n+\alpha) t^{\alpha+1}}{{\left(t^2+|x-y|^2\right)^{\frac{n+\alpha+2}{2}}}}\right]f(y)\ud y\\
		=&\frac{\alpha}{t}\int_{\R^n}\frac{f(x-ty)}{(1+|y|^2)^{\frac{n+\alpha}{2}}}\ud y-\frac{n+\alpha}{t}\int_{\R^n}\frac{f(x-ty)}{(1+|y|^2)^{\frac{n+\alpha+2}{2}}}\ud y.
	\end{align*}
	
	For the first term, with some $\sigma\in (0,1)$ we have
	\begin{align*}
		\int_{\R^n}\frac{f(x-ty)}{(1+|y|^2)^{\frac{n+\alpha}{2}}}\ud y
		=\int_{|y|> t^{-\sigma}}
		+\int_{|y|\leq t^{-\sigma}}\frac{f(x-ty)}{(1+|y|^2)^{\frac{n+\alpha}{2}}}\ud y:=I_1+I_2.
				\end{align*}
				For $I_1$, it is not hard to see that
				\begin{align*}
					|I_1|\lesssim \int_{|y|> t^{-\sigma}}	\frac{1}{(1+|y|^2)^{\frac{n+\alpha}{2}}}\ud y\lesssim t^{\sigma\alpha}.
				\end{align*}
				For $I_2$, as $t|y|<t^{1-\sigma}$ for all $|y|<t^{-\sigma}$,  the Taylor expansion together with symmetry enables us to estimate
				\begin{align*}
					I_2=&\int_{|y|\leq t^{-\sigma}}
						\frac{f(x)-t\nabla f(x)\cdot y+O(t^2|y|^2)}{(1+|y|^2)^{\frac{n+\alpha}{2}}}\ud y\\
						=&f(x)\int_{\R^n}\frac{1}{(1+|y|^2)^{\frac{n+\alpha}{2}}}\ud y-f(x)\int_{|y|>t^{-\sigma}}\frac{1}{(1+|y|^2)^{\frac{n+\alpha}{2}}}\ud y+\int_{|y|\leq t^{-\sigma}}\frac{O(t^2|y|^2)}{(1+|y|^2)^{\frac{n+\alpha}{2}}}\ud y.
				\end{align*}
				As $\alpha>1$, we now take $\sigma \in\left(\frac{1}{\alpha}, \frac{2}{\alpha}\right)$ for $\alpha \geq 2$ and $\sigma \in (\frac{1}{\alpha}, 1)$ for $\alpha<2$.
				Then we have
				\begin{align*}
					I_2=\frac{|\S^{n-1}|}{2}B(\frac{n}{2},\frac{\alpha}{2})f(x)+ O(t^{\sigma \alpha}).
				\end{align*}
				Thus, we obtain
				\begin{align*}
				\int_{\R^n}\frac{f(x-ty)}{(1+|y|^2)^{\frac{n+\alpha}{2}}}\ud y=\frac{|\S^{n-1}|}{2}B(\frac{n}{2},\frac{\alpha}{2})f(x)+ O(t^{\sigma \alpha}).
				\end{align*}
				
				For the second term, with the same $\sigma$ as above a similar argument yields
				$$\int_{\R^n}\frac{f(x-ty)}{(1+|y|^2)^{\frac{n+\alpha+2}{2}}}\ud y=\frac{|\S^{n-1}|}{2}B(\frac{n}{2},\frac{\alpha+2}{2}) f(x)+O(t^{\sigma(\alpha+2)})+O(t^2).$$
				
				Therefore, with the choice of $\sigma$ we combine these estimates to conclude that for $t \ll 1$,
				\begin{align*}
						&\partial_t\int_{\R^n}\frac{t^{\alpha}}{\left(t^2+|x-y|^2\right)^{\frac{n+\alpha}{2}}}f(y)\ud y\\
						=&\frac{|\S^{n-1}|}{2t}\left[\alpha B(\frac{n}{2},\frac{\alpha}{2})-(n+\alpha)B(\frac{n}{2},\frac{\alpha+2}{2})\right]f(x)+ O(t^{\sigma \alpha-1}) \\
						=&\ O(t^{\sigma \alpha-1}),
				\end{align*}
from which the desired conclusion follows.		
\end{pf}

\begin{rem}
The result in Proposition \ref{lem:singular_integral_II} is sharp. Otherwise, if $\alpha=1$ and let $f$ be a nonnegative function in $\Rn$ with compact support, then the classical Poisson integral
$$u(x,t)=P\ast f(x,t)$$
solves $\Delta u=0$ in $\R_+^{n+1}$ and $u=f(x)$ on $\pa \R_+^{n+1}$. However, for each $x \notin \mathrm{supp}(f)$, the Hopf boundary point lemma yields $\pa_t u>0$ on $\pa \R_+^{n+1}$.
\end{rem}

	\noindent \textbf{Proof of Theorem \ref{Thm:Poisson kernel_half-space}.} It suffices to show that fix $i,k \in \{0,1,2,3\}$, let
	$$ [(\mathscr{B}_i^3 P_k^3)\ast f_k](x,t)=\int_{\Rn}(\mathscr{B}_i^3 P_k^3)(x-y,t) f_k(y) \ud y,$$
	then we need to show
\begin{enumerate}[(i)]
\item $\Delta^2 P_i^3=0$ in $\R_+^{n+1}$.

\item For each pair $(i,k)$ other than $(0,2),(1,3)$, and fixed $x \in \Rn$ there hold 
\begin{align*}
\lim_{t\to 0} [(\mathscr{B}_i^3 P_i^3)\ast f_i](x,t)=f_i(x); \qquad \lim_{t\to 0} [(\mathscr{B}_i^3 P_k^3)\ast f_k](x,t)=0;\\ 
\lim_{t\to 0} [(\mathscr{B}_k^3 P_k^3)\ast f_k](x,t)=f_k(x); \qquad \lim_{t\to 0} [(\mathscr{B}_k^3 P_i^3)\ast f_i](x,t)=0.
\end{align*}
\end{enumerate}

For (i), notice that
\begin{align*}
\Delta^2 |X|^{3-n}=0, \qquad \pa_t |X|^{3-n}=(3-n)|X|^{1-n}t,
\end{align*}
this yields that $P_3^3$ is biharmonic in $\R_+^{n+1}$, so is $P_2^3=\pa_t P_3^3$. Observe that
\begin{align*}
\Delta |X|^{1-n}=0, \qquad \pa_t |X|^{1-n}=(1-n)|X|^{-1-n} t,
\end{align*}
then
\begin{align*}
\Delta^2 ( t \cdot t |X|^{-1-n})=\Delta\left(2 \pa_t (t |X|^{-1-n})\right)=\frac{2}{1-n}\pa_t^2 \Delta |X|^{1-n}=0,
\end{align*}
from which we know that $P_1^3$ is biharmonic in $\R_+^{n+1}$. Finally, we use 
\begin{equation}\label{fourth_Poisson_kernel}
\pa_t ( t^2 |X|^{-1-n})=2t |X|^{-1-n}-(1+n)t^3 |X|^{-3-n}
\end{equation}
to conclude that $\Delta^2 P_0^3=0$ in  $\R_+^{n+1}$.

\medskip

		In the following, we prove assertion (ii) case by case.
		
		\medskip
 \noindent\underline{\emph{Case 1.}} $(\mathscr{B}^3_0,\mathscr{B}^3_1)$-boundary condition.
\medskip

The kernel $P_0^3$ shares the common property as $P$ that
\begin{equation}\label{measure:Dirac_another}
\lim_{t\to 0}P_0^3\ast f_0(x,t)=f_0(x).
\end{equation}
This is due to
\begin{align*}
			&\lim_{t\to 0}\frac{2(n+1)}{|\S^n|}\int_{\R^n}\frac{t^3f_0(y)}{(t^2+|x-y|^2)^{\frac{n+3}{2}}} \ud y\\
			=&\lim_{t\to 0}\frac{2(n+1)}{|\S^n|}\int_{\R^n}\frac{f_0(x-ty)}{(1+|y|^2)^{\frac{n+3}{2}}} \ud y=f_0(x),
		\end{align*}
		where the last identity follows from
                  \begin{align*}
			\frac{2(n+1)}{|\S^n|}\int_{\R^n}\frac{1}{(1+|y|^2)^{\frac{n+3}{2}}}\ud y=&\frac{2(n+1)|\S^{n-1}|}{|\S^n|}\int_{0}^{+\infty}\frac{r^{n-1}}{(1+r^2)^{\frac{n+3}{2}}}\ud r\\
			=&\frac{(n+1)|\S^{n-1}|}{|\S^n|}B\left(\frac{n}{2},\frac{3}{2}\right)=1.
		\end{align*}
		Here, we emphasize $f_0 \in C(\Rn)\cap L^\infty(\Rn)$ is enough to ensure \eqref{measure:Dirac_another}, same as $P$. 
		
		Notice that
		$$\lim_{t \to 0}P_1^3 \ast f_1(x,t)=-\lim_{t \to 0} t [P\ast f_1](x,t)=0.$$
		
		By Lemma \ref{lem:singular_integral_II} we have
$$\lim_{t \to 0}[(\mathscr{B}_1^3 P_0^3) \ast f_0](x,t)=-\frac{2(n+1)}{|\Sn|}\lim_{t \to 0} \pa_t\int_{\R^n}\frac{t^{3}}{\left(t^2+|x-y|^2\right)^{\frac{n+3}{2}}}f_0(y)\ud y=0.$$

Again by \eqref{fourth_Poisson_kernel} we obtain
\begin{align*}
\lim_{t \to 0}[(\mathscr{B}_1^3 P_1^3) \ast f_1](x,t)=-\lim_{t \to 0}[(\pa_t P_1^3) \ast f_1](x,t)=\lim_{t \to 0}[(2P-P_0^3)\ast f_1](x,t)=f_1(x).
\end{align*}

\medskip
	 \noindent\underline{\emph{Case 2.}}  $(\mathscr{B}^3_0,\mathscr{B}^3_2)$-boundary condition.
	 \medskip
	 
	 Rewrite
	 $$\mathscr{B}^3_2=\partial_t^2-\overline{\Delta}=\Delta-2 \overline{\Delta}.$$
	 Recall that
	 $$P_0^3(X)=\frac{2(n+1)}{|\Sn|} t^3|X|^{-3-n},\quad P_2^3(X)=-\frac{1}{(n-1)|\Sn|}t |X|^{1-n}.$$
	Also \eqref{measure:Dirac_another} yields
	 $$\lim_{t \to 0} (P_0^3 \ast f_0)(x,t)=f_0(x).$$
	Also we have
	\begin{align}\label{limit:zeroth_second_type}
	\lim_{t \to 0} (P_2^3 \ast f_2)(x,t)=&-\frac{1}{(n-1)|\Sn|}\int_{\Rn}\frac{t}{(t^2+|x-y|^2)^{\frac{n-1}{2}}} f_2(y) \ud y\no\\
	=&-\frac{1}{(n-1)|\Sn|}\int_{\Rn}\frac{t^2 f_2(x-ty)}{(1+|y|^2)^{\frac{n-1}{2}}}  \ud y=0,
	\end{align}
	where the last identity follows from the following estimate: For each fixed $x$, we choose $M=M(x)>2|x|$ sufficiently large such that $f_2(x-ty)=O(|x-ty|^{-2-\delta_2})$ for all $|y|>M/t$, then
		\begin{align*}
			\left|\int_{\R^n}\frac{t^2f_2(x-ty)}{(1+|y|^2)^{\frac{n-1}{2}}}\ud y\right|\lesssim& \int_{|y|<\frac{M}{t}}\frac{t^2}{(1+|y|^2)^{\frac{n-1}{2}}}\ud y+ \int_{|y|>\frac{M}{t}}\frac{t^2}{|y|^{n-1}|x-ty|^{2+\delta_2}}\ud y\nonumber\\
			\lesssim& t+t^{-\delta_2}\int_{|y|>\frac{M}{t}}\frac{1}{|y|^{n+1+\delta_2}}\ud y\lesssim t.
		\end{align*}

	Taking $\Delta$ to Equation \eqref{fourth_Poisson_kernel} to show
	\begin{align*}
	(n+1)\Delta(t^3 |X|^{-3-n})=-\pa_t\Delta ( t \cdot t |X|^{-1-n})=-2 \pa_t^2 (t |X|^{-1-n})=2\overline{\Delta} (t |X|^{-1-n}).
	\end{align*}
	Then we obtain
	$$\mathscr{B}^3_2 P_0^3 =2\overline{\Delta}(P-P_0^3).$$
	This implies
	$$\lim_{t \to 0}[(\mathscr{B}^3_2  P_0^3) \ast f_0](x,t)=2 \lim_{t\to 0}\overline{\Delta}[(P-P_0^3)\ast f_0](x,t)=0.$$
	
	Next we combine
	$$\Delta (t |X|^{1-n})=2\pa_t |X|^{1-n}=2(1-n)t |X|^{-1-n}$$
	and 
	\begin{align}\label{eq:elementary}
	-2\overline{\Delta} (t|X|^{1-n})=2t\pa_t^2(|X|^{1-n})=2(1-n)[t|X|^{-1-n}-(n+1)t^3 |X|^{-3-n}]
	\end{align}
	to show
	$$\mathscr{B}^3_2 P_2^3=P+(P-P_0^3),$$
	from which we have
	\begin{equation}\label{second_third}
	\lim_{t \to 0}[(\mathscr{B}^3_2  P_2^3) \ast f_2](x,t)=f_2(x).
	\end{equation}

\medskip
 \noindent\underline{\emph{Case 3.}} Type $(\mathscr{B}^3_1,\mathscr{B}^3_3)$ boundary condition.
\medskip

 A direct computation yields
\begin{align*}
				\mathscr{B}^3_1P^3_1=\partial_t\left(\frac{2}{|\S^n|}t^2|X|^{-1-n}\right)=&\frac{2}{|\S^n|}\left(\frac{2t}{(t^2+|x|^2)^{\frac{n+1}{2}}}-\frac{(n+1)t^3}{(t^2+|x|^2)^{\frac{n+3}{2}}}\right)\\
				:=&2P-P_0^3,
			\end{align*}
			 		Then, it is not hard to see that $\lim_{t\to 0} [(\mathscr{B}^3_1P^3_1)\ast f_1](x,t)=f_1(x)$.
		
		Since
		\begin{align*}
				\mathscr{B}^3_1P^3_3=\frac{1}{(n-1)|\Sn|}t |X|^{1-n},
		\end{align*}
 similar to \eqref{limit:zeroth_second_type} we have
		\begin{align*}
		\lim_{t\to 0}[(\mathscr{B}_1^3 P_3^3)\ast f_3](x,t)=&\lim_{t\to 0}\int_{\Rn}\frac{1}{(n-1)|\Sn|}\frac{t f_3(y)}{(t^2+|x-y|^2)^{\frac{n-1}{2}}} \ud y\\
		=&\lim_{t\to 0}\frac{1}{(n-1)|\Sn|}\int_{\Rn}\frac{t^2 f_3(x-ty)}{(1+|y|^2)^{\frac{n-1}{2}}} \ud y=0.
		\end{align*}		
		
		Recall that
		$$\mathscr{B}^3_3=\partial_t\Delta +2\overline{\Delta}\partial_t$$
		and
		$$ P_1^3(X)=-\frac{2}{|\Sn|}t^2 |X|^{-1-n},\quad P_3^3(X)=\frac{1}{(n-1)(n-3)|\Sn|} |X|^{3-n}.$$
	A direct computation yields
	\begin{align*}
	\pa_t \Delta (t \cdot t |X|^{-1-n})=&2\pa_t^2(t|X|^{-1-n});\\
	2\overline{\Delta}\partial_t (t \cdot t |X|^{-1-n})=&2\overline{\Delta}\left[2t |X|^{-1-n}-(1+n)t^3 |X|^{-3-n}\right].
	\end{align*}	
	Since $t |X|^{-1-n}$ is harmonic in $\R_+^{n+1}$, we have
	$$\mathscr{B}^3_3  P_1^3=-2\overline{\Delta}(P-P_0^3).$$
	This directly implies
	\begin{align*}
	\lim_{t \to 0}[(\mathscr{B}^3_3  P_1^3) \ast f_1](x,t)=-2 \lim_{t\to 0}\overline{\Delta} [(P-P_0^3) \ast f_1](x,t)=0.
	\end{align*}
	
	A direct computation yields
	\begin{align*}
	\pa_t \Delta |X|^{3-n}=&(3-n)\Delta(t |X|^{1-n})=2(n-1)(n-3) t |X|^{-1-n},\\
	2\overline{\Delta}\partial_t |X|^{3-n}=&2(3-n) \overline{\Delta} (t |X|^{1-n}).
	\end{align*}
	These together with \eqref{eq:elementary} give
	$$\mathscr{B}^3_3  P_3^3=P-\overline{\Delta}(P-P_0^3).$$
	This directly implies
 	\begin{align}\label{third_third}
	\lim_{t \to 0}[(\mathscr{B}^3_3  P_3^3) \ast f_3](x,t)=&f_3(x)-\overline{\Delta} \lim_{t\to 0} [(P-P_0^3)\ast f_3](x,t)=f_3(x).
	\end{align}
	
	\medskip
 \noindent\underline{\emph{Case 4.}} $(\mathscr{B}^3_2,\mathscr{B}^3_3)$-boundary condition.
\medskip

Since we already have \eqref{second_third} and \eqref{third_third}, it remains to show
\begin{align*}
\lim_{t \to 0}[(\mathscr{B}^3_2 P_3^3)\ast f_3](x,t)=0, \qquad \lim_{t \to 0}[(\mathscr{B}^3_3 P_2^3)\ast f_2](x,t)=0.
\end{align*}

To this end, notice that
\begin{align*}
\mathscr{B}^3_2 |X|^{3-n}=&\Delta |X|^{3-n}-2\overline{\Delta} |X|^{3-n}\\
=&2(3-n)|X|^{1-n}-2(3-n)[(n-1)|X|^{-1-n}t^2+|X|^{1-n}]\\
=&2(n-1)(n-3)t^2|X|^{-1-n},
\end{align*}
whence,
$$\mathscr{B}^3_2 P_3^3=\frac{2}{|\Sn|}t^2 |X|^{-1-n}.$$
This implies
$$\lim_{t \to 0}[(\mathscr{B}^3_2 P_3^3)\ast f_3](x,t)=\lim_{t \to 0} [t P\ast f_3(x,t)]=0.$$

Since 
\begin{align*}
\mathscr{B}^3_3(t |X|^{1-n})=&\pa_t \Delta (t |X|^{1-n})+2\overline{\Delta}\partial_t(t |X|^{1-n})\\
=&2\pa_t^2 |X|^{1-n}+2\overline{\Delta}|X|^{1-n}+2t\overline{\Delta} \pa_t |X|^{1-n}\\
=&2(1-n)t^2 \overline{\Delta}|X|^{-1-n},
\end{align*}
we have	
$$\mathscr{B}^3_3 P_2^3=\frac{2}{|\Sn|}\overline{\Delta} (t^2 |X|^{-1-n}).$$
This gives
$$\lim_{t \to 0}[(\mathscr{B}^3_3 P_2^3)\ast f_2](x,t)=\overline{\Delta}\lim_{t \to 0} [t P\ast f_2(x,t)]=0.$$ 

\medskip

For $n=3$, we only need to point out several significant differences. Observe that
	\begin{align*}
		\mathscr{B}^3_3 P^3_3(X)=(\partial_t\Delta+2\overline{\Delta}\partial_t)P^3_3(X)=\frac{8}{|\S^3|}\frac{t^3}{(t^2+|x|^2)^3}=P^3_0(X),
	\end{align*}
	so, $\lim_{t\to 0} [(\mathscr{B}^3_3P^3_3)\ast f_3](x,t)=f_3(x)$. Also, we can apply
	\begin{align*}
		\mathscr{B}^3_1P^3_3(X)=&\frac{t}{2|\S^3|(t^2+|x|^2)}=-P_2^3(X),\\
		 \mathscr{B}^3_2P^3_3(X)=&\frac{2}{|\S^3|}\frac{t^2}{(t^2+|x|^2)^2}=-P_1^3(X)
	\end{align*}
to show
	\begin{align*}
		\lim_{t\to 0} [(\mathscr{B}^3_1P^3_3)\ast f_3](x,t)=0 \qquad \mathrm{and}\qquad \lim_{t\to 0} [(\mathscr{B}^3_2P^3_3)\ast f_3](x,t)=0.
	\end{align*}
\hfill $\Box$

\subsection{The unit ball}

Our purpose is to give explicit representation formulae of \emph{biharmonic Poisson kernels} for  $\B^{n+1}$ with the help of conformal equivalence between $\B^{n+1}$ and $\R_+^{n+1}$.

\medskip

Let $F: (\R^{n+1}_{+}, |\ud X|^2) \to (\B^{n+1},|\ud \xi|^2) $ 
\begin{equation}\label{conf-map:ball_half-space}
	\xi=F(X)=-\mathbf{e}_{n+1}+\frac{2(X+\mathbf{e}_{n+1})}{|X+\mathbf{e}_{n+1}|^2}
	\end{equation}
	be a conformal map with the property that
	$$F^\ast(|\ud \xi|^2)=\left(\frac{2}{|X+\mathbf{e}_{n+1}|^2}\right)^2 |\ud X|^2:= U_0(X)^{\frac{4}{n-3}} |\ud X|^2.$$
	Clearly, $F$ is an inversion of the sphere centered at $-\mathbf{e}_{n+1}$ of radius $\sqrt{2}$ and  its inverse $F^{-1}=F$.

	The following  elementary identity will be frequently used: 
	\begin{equation}\label{formula:dist-fcns_two_models}
	2|\xi+\mathbf{e}_{n+1}|^{-1} |\eta+\mathbf{e}_{n+1}|^{-1}|\xi-\eta|= |F^{-1}(\xi)-F^{-1}(\eta)|, \quad \forall~\xi,\eta \in \overline{\B^{n+1}}\setminus\{-\mathbf{e}_{n+1}\}.
	\end{equation}
	This follows by a straightforward calculation.
	
	For $U \in C^\infty(\overline{\B^{n+1}})$ we let 
	\begin{align}\label{fcns:ball_half-space}
	u= U_0 (U\circ F).
	\end{align}
	Clearly, $u(X)=O^{(k)}(|X|^{3-n})$ as $|X| \to \infty$, where $O^{(k)}(|X|^{3-n})$ means the $k$-th derivatives of $u$ belong to $O(|X|^{3-n-k})$ for $k \in \N$. Moreover, by conformal covariance of the Paneitz operator, we have
	\begin{align*}
	\Delta^2 u =P_4^{|\ud X|^2}(u)=P_4^{|\ud X|^2}( U_0 U\circ F)=&U_0^{\frac{n+5}{n-3}}P_4^{U_0^{\frac{4}{n-3}}|\ud X|^2}(U\circ F)\\
	=&U_0^{\frac{n+5}{n-3}}(P_4^{|\ud \xi|^2} U)\circ F=U_0^{\frac{n+5}{n-3}}(\Delta^2 U)\circ F.
	\end{align*}
	Similarly,  by \eqref{conf_transformation:bdry_operators} we have
	\begin{equation}\label{conf_change_bdry_operators}
	\left(\mathscr{B}_k^3\right)_{|\ud X|^2} (u)=U_0^{\frac{n+2k-3}{n-3}}\left[\left(\mathscr{B}_k^3\right)_{|\ud \xi|^2}U\right] \circ F \qquad \mathrm{~~for~~}\quad   k\in \{0,1,2,3\}.
	\end{equation}

\begin{thm}\label{Thm:Poisson kernel_ball}
		Suppose $n \geq 2$ and $n\neq 3$, and $\bar f_i\in C^{\infty}(\S^n)$ for $i \in \{0,1,2,3\}$. For $\xi \in \B^{n+1}$ and $\eta \in \Sn$, we let
		\begin{align*}
		\bar P_0^3(\xi,\eta)=&\frac{n+1}{4|\S^n|} \frac{\left(1-|\xi|^2\right)^3}{|\xi-\eta|^{n+3}};\\
		\bar P_1^3(\xi,\eta)=&-\frac{1}{2|\S^n|}\frac{(1-|\xi|^2)^2}{|\xi-\eta|^{n+1}};\\
		\bar P_2^3(\xi,\eta)=&-\frac{1}{2(n-1)|\S^n|}\frac{1-|\xi|^2}{|\xi-\eta|^{n-1}};\\
		\bar P_3^3(\xi,\eta)=&\frac{1}{(n-1)(n-3)|\S^n|} \frac{1}{|\xi-\eta|^{n-3}}.
		\end{align*}
		Then for $i, j \in \{0,1,2,3\},i<j, i+j\neq 3$, the solution of the following  BVP
		\begin{align}\label{Biharmonic_bdry_value_problem_on_balls}
			\begin{cases}
		\displaystyle\Delta^2 U=0   &\mathrm{~~in~~}\quad \B^{n+1},\\
		\displaystyle \mathscr{B}^3_i U=\bar f_i  &\mathrm{~~on~~}\quad \Sn,\\
		\displaystyle \mathscr{B}^3_j U=\bar f_j\qquad \qquad&\mathrm{~~on~~}\quad \Sn,
		\end{cases}
		\end{align}
		is unique and has the form of
		$$U(\xi)=\int_{\Sn} \bar P_i^3(\xi,\eta) \bar f_i(\eta) \ud V_{\S^n}(\eta)+\int_{\Sn} \bar P_j^3(\xi,\eta) \bar f_j(\eta) \ud V_{\S^n}(\eta).$$
	\end{thm}
\begin{pf}
We first handle  $(\mathscr{B}^3_0,\mathscr{B}^3_2)$ type boundary condition.
	
	Let $u$ define as in \eqref{fcns:ball_half-space}, equivalently,
	$$U(\xi)=\left(\frac{2}{|\xi+\mathbf{e}_{n+1}|^2}\right)^{\frac{n-3}{2}} u\circ F^{-1}(\xi).$$

	By virtue of \eqref{conf_change_bdry_operators} we set
	\begin{align*}
	 f_0(x)=\left(\frac{2}{1+|x|^2}\right)^{\frac{n-3}{2}}\bar f_0\circ F(x,0) \quad \mathrm{and}\quad	f_2(x)=\left(\frac{2}{1+|x|^2}\right)^{\frac{n+1}{2}}\bar f_2\circ F(x,0).
	\end{align*}
	Then for $|X|,|x|\gg1$,  we know
	\begin{align*}
	u(X)=O(|X|^{3-n}), \quad  f_0(x)=O(|x|^{3-n}),\quad f_2(x) =O(|x|^{-1-n}),
	\end{align*}
	 and $u$ solves 
	\begin{align}\label{Sec 4.2 boundary equ a}
	\begin{cases}
	\displaystyle \qquad~ \Delta^2 u=0 &\mathrm{~~in~~}\qquad~ \R_{+}^{n+1},\\
	\displaystyle  \qquad \quad~~	u= f_0 &\mathrm{~~on~~}\qquad \partial\R^{n+1}_{+},\\
	\displaystyle \partial^2_t u-\overline{\Delta}u= f_2 \qquad\quad &\mathrm{~~on~~}\qquad \partial\R^{n+1}_{+}.
	\end{cases}
	\end{align}
	
	As a direct consequence of Remark \ref{Uniqueness lem 2 rem 2},  the above problem \eqref{Sec 4.2 boundary equ a} admits a unique solution. This together with Theorem \ref{Thm:Poisson kernel_half-space} gives
		\begin{align*}
	u(X)=\frac{2(n+1)}{|\S^n|}\int_{\R^n}\frac{t^3 f_0(y)}{(t^2+|x-y|^2)^{\frac{n+3}{2}}}\ud y
	-\frac{1}{(n-1)|\S^n|}\int_{\R^n}\frac{t f_2(y)}{(t^2+|x-y|^2)^{\frac{n-1}{2}}} \ud y.
	\end{align*}
	
	If we let
	\begin{align*}
	\xi=&F(X)=\left(\frac{2x}{(1+t)^2+|x|^2},\frac{1-t^2-|x|^2}{(1+t)^2+|x|^2}\right),\\
	\eta=&F(y,0)=\left(\frac{2y}{1+|y|^2},\frac{1-|y|^2}{1+|y|^2}\right),
	\end{align*}
	then 
	$$t\circ F^{-1}=\frac{1-|\xi|^2}{|\xi+\mathbf{e}_{n+1}|^2}$$
	and
	$$|\xi+\mathbf{e}_{n+1}|^2=\frac{4}{(1+t)^2+|x|^2},  \quad |\eta+\mathbf{e}_{n+1}|^2=\frac{4}{1+|y|^2}.$$
	Moreover, it follows from \eqref{formula:dist-fcns_two_models} that
	\begin{equation*}
	2|\xi+\mathbf{e}_{n+1}|^{-1} |\eta+\mathbf{e}_{n+1}|^{-1}|\xi-\eta|= |(x,t)-(y,0)|.
	\end{equation*}
	
	Therefore, putting these facts together we obtain
		\begin{align*}
	U(\xi)=&\left(\frac{2}{|\xi+\mathbf{e}_{n+1}|^2}\right)^{\frac{n-3}{2}}\frac{2(n+1)}{|\S^n|}\int_{\Sn}\frac{(t\circ F^{-1})^3 f_0\circ F^{-1}(\eta) \left(\frac{2}{|\eta+\mathbf{e}_{n+1}|^2}\right)^n}{((t\circ F^{-1})^2+|(x-y)\circ F^{-1}|^2)^{\frac{n+3}{2}}}\ud y\\
	&-\left(\frac{2}{|\xi+\mathbf{e}_{n+1}|^2}\right)^{\frac{n-3}{2}}\frac{1}{(n-1)|\S^n|}\int_{\Sn}\frac{t\circ F^{-1} f_2\circ F^{-1}(\eta) \left(\frac{2}{|\eta+\mathbf{e}_{n+1}|^2}\right)^n}{((t\circ F^{-1})^2+|(x-y)\circ F^{-1}|^2)^{\frac{n-1}{2}}} \ud y\\
	=&\frac{(n+1)\left(1-|\xi|^2\right)^3}{4|\S^n|}\int_{\S^n}\frac{\bar f_0(\eta)}{|\xi-\eta|^{n+3}}\ud V_{\S^n}(\eta)-\frac{(1-|\xi|^2)}{2(n-1)|\S^n|}\int_{\S^n}  \frac{\bar f_2(\eta)}{|\xi-\eta|^{n-1}} \ud V_{\S^n}(\eta).
	\end{align*}
	
	\medskip
	
	As indicated in the preceding discussion, for other proper conformal boundary operator pairs it remains to show the uniqueness. 
	
	The uniqueness for $(\mathscr{B}^3_0,\mathscr{B}^3_1)$ directly follows by an integration by parts. For $(\mathscr{B}^3_1,\mathscr{B}^3_3)$, the uniqueness follows from a similar argument above combined with Lemma \ref{Uniqueness lem} or Remark \ref{Rem:Uniqueness}. 
	
	The remaining boundary condition is  $(\mathscr{B}^3_2,\mathscr{B}^3_3)$, the uniqueness becomes a little subtle. Suppose $U$ is a smooth solution of
	\begin{align*}
	\begin{cases}
	\displaystyle\Delta^2 U=0 \qquad \qquad &\mathrm{in}\qquad \B^{n+1},\\
	\displaystyle \mathscr{B}^3_2 U=0  &\mathrm{on}\qquad \S^n,\\
	\displaystyle \mathscr{B}^3_3 U=0  &\mathrm{on}\qquad \S^n,
	\end{cases}
	\end{align*}
	and let $u$ define as in \eqref{fcns:ball_half-space}, then $u$ solves 
	\begin{align*}
\begin{cases}
\displaystyle \qquad \qquad~ \Delta^2 u=0 &\mathrm{~~in~~}\qquad~~ \R_{+}^{n+1},\\
\displaystyle 	\qquad \partial_t^2 u-\overline{\Delta} u=0 &\mathrm{~~on~~}\qquad \partial\R^{n+1}_{+},\\
\displaystyle \partial_t\Delta u+2\overline{\Delta}\partial_t u=0\qquad\quad&\mathrm{~~on~~}\qquad \partial\R^{n+1}_{+},
\end{cases}
\end{align*}
	with $u(X)=O(|X|^{3-n})$ as $|X| \to \infty$.
	
By Lemma \ref{Uniqueness lem} and Theorem \ref{Thm:Poisson kernel_half-space}, the unique solution for $(\mathscr{B}^3_1,\mathscr{B}^3_3)$ is
\begin{align*}
	u(x,t)=\frac{2}{|\S^n|}\int_{\R^n}\frac{t^2\partial_t u(y, 0)}{(t^2+|x-y|^2)^{\frac{n+1}{2}}} \ud y.
\end{align*}
Since $\partial_t u(y,0)=O(|y|^{1-n})$ as $|y| \to \infty$ by definition \eqref{fcns:ball_half-space} of $u$, we have
\begin{align*}
\lim_{t\to 0}\frac{2t}{|\S^n|}\int_{\R^n}\frac{t\partial_t u(y, 0)}{(t^2+|x-y|^2)^{\frac{n+1}{2}}} \ud y=0.
\end{align*}
So, by the uniqueness for  $(\mathscr{B}^3_0,\mathscr{B}^3_2)$ due to Lemma \ref{Uniqueness lem 2} we obtain $u= 0$ and $U=0$ as well. 
\end{pf}

\begin{rem}
	It is notable that restricting $\xi$ to $\Sn$, $2\bar P_3^3(\xi,\eta)$ coincides with the Green function for the fractional GJMS operator $P_3$ of order three on $\Sn$ (cf. the first author and Shi \cite[Theorem 1(3)]{Chen-Shi}).  
	\end{rem}

\begin{lem}\label{Lem:nonuniqueness}
	The following homogeneous biharmonic  BVPs:
	\begin{align}\label{Sec 4 bad boundary type 1}
	\begin{cases}
	\displaystyle ~~~\Delta^2 U=0  &\mathrm{in}\qquad \B^{n+1}\\
	\displaystyle \mathscr{B}^3_0(U)=0  &\mathrm{on}\qquad \S^n\\
	\displaystyle \mathscr{B}^3_3(U)=0\qquad \quad &\mathrm{on}\qquad \S^n
	\end{cases}
	\end{align}
	and
	\begin{align}\label{Sec 4 bad boundary type 2}
	\begin{cases}
	\displaystyle ~~~\Delta^2 U=0  &\mathrm{in}\qquad \B^{n+1}\\
	\displaystyle \mathscr{B}^3_1(U)=0  &\mathrm{on}\qquad \S^n\\
	\displaystyle \mathscr{B}^3_2(U)=0\qquad \quad &\mathrm{on}\qquad \S^n
	\end{cases}
	\end{align}
	have nontrivial smooth solutions.
\end{lem}
\begin{pf}
	For every $\varphi\in C_{\mathrm{c}}^{\infty}(\R^n)$, we introduce
	\begin{align*}
	u_\varphi(x,t)=-\frac{2}{|\S^n|}\int_{\R^n}\frac{t^2\varphi(y)}{(t^2+|x-y|^2)^{\frac{n+1}{2}}}\ud y
	\end{align*}
	for the boundary condition $(\mathscr{B}^3_0,\mathscr{B}^3_3)$;   and
	\begin{align*}
	u_\varphi (x,t)=\frac{1}{(n-1)(n-3)|\Sn|}\int_{\R^n}\frac{\varphi(y)}{(t^2+|x-y|^2)^{\frac{n-3}{2}}}\ud y
	\end{align*}
	for the boundary condition $(\mathscr{B}^3_1,\mathscr{B}^3_2)$.

	As in the proof of Theorem \ref{Thm:Poisson kernel_half-space}, we know that $u=u_\varphi$ satisfies both
	\begin{align*}
	\begin{cases}
	\displaystyle \Delta^2 u=0\qquad\qquad\qquad\qquad\, &\mathrm{in}\qquad \R_{+}^{n+1}\\
	\displaystyle  \quad~	u= 0\qquad\qquad \qquad\qquad\,\,&\mathrm{on}\quad\,\,\,\, \partial\R^{n+1}_{+}\\
	\displaystyle \partial_t\Delta u+2\overline{\Delta}\partial_t u=0\qquad\,\,\,\,\,\,\,&\mathrm{on}\quad\,\,\,\, \partial\R^{n+1}_{+}
	\end{cases}
	\end{align*}
	and
	\begin{align*}
	\begin{cases}
	\displaystyle \Delta^2 u=0\qquad\qquad\qquad\qquad\, &\mathrm{in}\qquad \R_{+}^{n+1}\\
	\displaystyle 	~~\partial_t u= 0\qquad\qquad \qquad\qquad\,\,&\mathrm{on}\quad\,\,\,\, \partial\R^{n+1}_{+}\\
	\displaystyle \partial^2_t u-\overline{\Delta} u=0\qquad\,\,\,\,\,\,\,&\mathrm{on}\quad\,\,\,\, \partial\R^{n+1}_{+}
	\end{cases},
	\end{align*}
	respectively. Finally, we may obtain nontrivial smooth solutions of \eqref{Sec 4 bad boundary type 1} and \eqref{Sec 4 bad boundary type 2} via 
	\begin{align*}
	U_\varphi(\xi)=\left(\frac{2}{|\xi+\mathbf{e}_{n+1}|^2}\right)^{\frac{n-3}{2}}u_\varphi\circ F^{-1}(\xi).
	\end{align*}
\end{pf}

This also reflects the loss of \emph{`the complementing condition'} for $(\mathscr{B}^3_0,\mathscr{B}^3_3)$ and $(\mathscr{B}^3_1,\mathscr{B}^3_2)$. Due to this reason, problem \eqref{Sec 4 bad boundary type 1} or \eqref{Sec 4 bad boundary type 2} fails to have the uniqueness.

\medskip

Motivated by Shi and the first author \cite[Theorem 1(1)]{Chen-Shi} we are able to derive explicit formulae of \emph{biharmonic Poisson kernel}s for $\R_+^4$ and $\B^{4}$. 

\begin{thm}\label{Thm:biharmonic_Poisson_kernerl_dim-4}
			Suppose $\bar f_i\in C^{\infty}(\S^3)$ for $i \in \{0,1,2,3\}$.  Then the biharmonic Poisson kernels $\bar P_i^3$ for $\B^4$ have the same representation formulae as in Theorem \ref{Thm:Poisson kernel_ball}  except for
			\begin{align*}
				\bar P_3^3(\xi,\eta)=-\frac{1}{2|\S^3|}\log|\xi-\eta|+C, \qquad C \in \R
		\end{align*}
			 where $\xi \in \B^{4}, \eta \in \S^3$. Moreover, we have
			\begin{enumerate}[(i)]
			\item $i=1,2$ and $j=3$, the solution of \eqref{Biharmonic_bdry_value_problem_on_balls} is
		\begin{align*}
&U(\xi)-\frac{1}{|\S^3|}\int_{\S^3}U   \ud V_{{\S^3}}\\
=&\int_{\S^3} \bar{P}^3_i(\xi-\eta)\bar{f}_i(\eta) \ud V_{{\S^3}}(\eta)-\frac{1}{2|\S^3|}\int_{\S^3}\log |\xi-\eta|\bar{f}_3(\eta)\ud V_{\S^3}(\eta), \quad \forall~\xi \in \overline{\B^4}.
\end{align*}

\item $(i,j)=(0,1) \mathrm{~or~} (0,2)$, the unique solution of \eqref{Biharmonic_bdry_value_problem_on_balls} is
		\begin{align*}
	U(\xi)=\int_{\S^3} \bar{P}^3_i(\xi-\eta)\bar{f}_i(\eta) + \bar{P}^3_j(\xi-\eta)\bar{f}_j(\eta)\ud V_{{\S^3}}(\eta), \quad \forall~\xi \in \overline{\B^4}.
\end{align*}	
\end{enumerate}
			
			\end{thm}
	
		\medskip
			
\begin{pf} Suppose $U$ is a smooth solution of \eqref{Biharmonic_bdry_value_problem_on_balls}. For $n=3$, we introduce 
	\begin{align*}
	u(x,t)=U\circ F(x,t)+\log \frac{2}{(1+t)^2+|x|^2},  \qquad \hat U_0(x,t):=\log \frac{2}{(1+t)^2+|x|^2}
\end{align*}
and
\begin{align*}
f_0(x)=&\bar f_0\circ F(x,0)+\hat U(x,0),\\
	f_k(x)=&\left[\bar{f}_k\circ F(x,0)+(T_k^3)_{\S^3}\right] e^{k\hat U_0(x,0)}, k \in \{1,2,3\}.
\end{align*}
By the conformal covariance of $P_4^g$ and $\mathscr{B}_k^3$ we have
\begin{align*}
e^{2u} |\ud X|^2=F^\ast(e^{2U} |\ud \xi|^2), \quad e^{2 \hat U_0} |\ud X|^2=F^\ast( |\ud \xi|^2)
\end{align*}
and
\begin{align*}
f_k=(\mathscr{B}_k^3)_{|\ud X|^2}(u)=&(\mathscr{B}_k^3)_{|\ud X|^2}(U\circ F)+(\mathscr{B}_k^3)_{|\ud X|^2}(\hat U_0)\\
\overset{\eqref{conf_transformation:bdry_operators}}{=}&e^{k\hat U_0}(\mathscr{B}_k^3)_{e^{2 \hat U_0}|\ud X|^2}(U \circ F)+(T_k^3)_{\S^3}e^{3\hat U_0}\\
=&e^{k\hat U_0}(\mathscr{B}_k^3)_{F^\ast(|\ud \xi|^2)}(U \circ F)+(T_k^3)_{\S^3}e^{k\hat U_0}\\
=&e^{k\hat U_0}[(\mathscr{B}_3^3)_{|\ud \xi|^2}(U)]\circ F+(T_k^3)_{\S^3}e^{k\hat U_0}\\
=&e^{k\hat U_0}(\bar f_k \circ F+(T_k^3)_{\S^3}).
\end{align*}
In particular, $(T_1^3)_{\S^3}=1,(T_2^3)_{\S^3}=2, (T_3^3)_{\S^3}=4$.

Suppose $u$ is a solution of \eqref{prob:bdry_value-Possion_kernel}. Using the biharmonic Poisson kernel for $\R_+^4$ we have
\begin{align*}
u(x,t)=&\int_{\pa \R_+^4}[P_i^3(x-y,t)f_i(y)+P_j^3(x-y,t)f_j(y)] \ud y\\
=&\int_{\pa \R_+^4}[P_i^3(x-y,t) e^{i \hat U_0(y,0)}f_i\circ F(y)+P_j^3(x-y,t) e^{j \hat U_0(y,0)}f_j\circ F(y)]\ud y+I(x,t),
\end{align*}
where
\begin{align*}
I(x,t)=\begin{cases}
\displaystyle \int_{\pa \R_+^4}[P_0^3(x-y,t)\hat U_0(y,0)+P_j^3(x-y,t)e^{j \hat U_0(y,0)}(T_j^3)_{\S^3}]\ud y,\\
\displaystyle \int_{\pa \R_+^4}[P_i^3(x-y,t) e^{i \hat U_0(y,0)}f_i\circ F(y)+P_j^3(x-y,t) e^{j \hat U_0(y,0)}f_j\circ F(y)]\ud y~~ \mathrm{~~if~~} i \neq 0.
\end{cases}
\end{align*}
Then it is not hard to see that $I=\hat U_0$ when $(i,j)=(0,1),(0,2)$; $I=\hat U_0+C, C \in \R$ when $(i,j)=(1,3),(2,3)$.

Hence, using the variable changes as in the proof of Theorem \ref{Thm:Poisson kernel_ball} we obtain
\begin{align*}
U(\xi)=\int_{\S^3} \bar{P}^3_i(\xi-\eta)\bar{f}_i(\eta) + \bar{P}^3_j(\xi-\eta)\bar{f}_j(\eta)\ud V_{{\S^3}}(\eta)
\end{align*}
when $(i,j)=(0,1),(0,2)$; and
\begin{equation}\label{Poisson_kernel:dim_four}
U(\xi)=\int_{\S^3} \bar{P}^3_i(\xi-\eta)\bar{f}_i(\eta) \ud V_{{\S^3}}(\eta)-\frac{1}{2|\S^3|}\int_{\S^3}\log |\xi-\eta|\bar{f}_3(\eta)\ud V_{\S^3}(\eta)+C,
\end{equation}
 when $(i,j)=(1,3),(2,3)$.

It remains to determine the constant $C$. On $\S^3$, a compatible condition of $\bar f_3$ arises from the fact that $\mathrm{ker} P_3=\R$ by virtue of \cite[Theorem 1(1)]{Chen-Shi}, more explicitly, 
\begin{equation}\label{compatible-cond:P_3}
\int_{\S^3} \bar f_3 \ud V_{\S^3}=\int_{\S^3}(-\frac{\partial \Delta U}{\partial r}-2\Delta_{\S^3}\frac{\partial U}{\partial r}-2\Delta_{\S^3}U) \ud V_{\S^3}=0.
\end{equation}
Restricting $\xi$ to $\S^3$ in \eqref{Poisson_kernel:dim_four} we find
\begin{align*}
	U(\xi)=-\frac{1}{2|\S^3|}\int_{\S^3}\log |\xi-\eta|\bar{f}_3(\eta)\ud V_{\S^3}+C.
\end{align*}
Integrating the above equation over $\S^3$  gives
\begin{align*}
\int_{\S^3}U \ud V_{\S^3}=-\frac{1}{2|\S^3|}\int_{\S^3}\bar{f}_3(\eta)\ud V_{\S^3}(\eta) \int_{\S^3}\log |\xi-\eta| \ud V_{\S^3}(\xi)+C|\S^3|=C |\S^3|,
\end{align*}
where we have used \eqref{compatible-cond:P_3} and the fact that $\forall~ \eta \in \S^3$, $\int_{\S^3}\log |\xi-\eta| \ud V_{\S^3}(\xi)$ is constant. This follows that
\begin{align*}
	C=\frac{1}{|\S^3|}\int_{\S^3}  U \ud V_{{\S^3}}.
\end{align*}
Then the desired assertion follows by substituting the above constant into \eqref{Poisson_kernel:dim_four}.
\end{pf}

\medskip

Finally, we point out that the proof of Theorems \ref{Thm:Poisson kernel_half-space} and \ref{Thm:Poisson kernel_ball} is still valid for $n=2$. However, we shall not go deep into this, since it has been well studied by the authors \cite{Chen-Zhang}.

\medskip

For a Poincar\'e-Einstein manfold $(X^{n+1},g_+)$ with conformal infinity $(\pa X=M,[h])$, and assume $(n^2-1)/4, (n^2-9)/4 \notin \sigma_{\mathrm{pp}}(-\Delta_{g_+})$, with the factorization of $P_4^{g_+}$ and explicit formulae of extrinsic GJMS boundary operators $\mathscr{B}_k^3$, Case \cite{Case} used the scattering operator to find fractional GJMS operators $P_1^g,P_3^g$ from  directly computing $(\mathscr{B}_k^3)_g u$, where $P_4^g u=0$ in $X$ for $g=r^2 g_+$ with $r$ being a geodesic defining function. As for the Poincar\'e ball $(\B^{n+1}, g_+=(\frac{2}{1-|x|^2})^2 |\ud x|^2)$, the spectrum $\sigma(-\Delta_{g_+})=[n^2/4,+\infty)$ and  \cite[Theorem 1.4]{Case} is applicable.

\begin{prop}\label{Prop:ex-intrinsic_GJMS}
			Suppose $U$ is a biharmonic function in $\B^{n+1}$ for $n\geq 2$, then there hold
			\begin{align*}
				\mathscr{B}^3_3U=2P_3 U \qquad \mathrm{and}\qquad \mathscr{B}^3_2\left(U\right)=2P_1(\mathscr{B}^3_1U) \qquad \mathrm{~~on~~} \Sn,
			\end{align*}
			where $P_1,P_3$ are the first- and third-order GJMS operators with respect to the round metric, see \cite{Case,Chen-Shi}.
			\end{prop}
			
			An alternative proof of Proposition \ref{Prop:ex-intrinsic_GJMS} using the \emph{biharmonic Poisson kernel} is available in Appendix \ref{Append:A}.

			\section{Biharmonic Green function}\label{Sect4}
		For $n\geq 2$, let $X=(x,t) \in \R_+^{n+1},Y=(y,s) \in \R_+^{n+1}$. We define  $\bar{X}=(x,-t)$, and $\xi=F(X),\eta=F(Y)$ with $F$ being the conformal map as in \eqref{conf-map:ball_half-space}. 	
		
		We first introduce a \emph{biharmonic Green function} on $\B^{n+1}$ with each of four proper conformal boundary operator pairs. Consider an nonhomogeneous  biharmonic BVP with vanishing boundary conditions
	\begin{align}\label{eqn:def_Green_fcn}
	\begin{cases}
	\Delta^2  U=\bar f \qquad\qquad &\mathrm{~~in~~}\quad \B^{n+1},\\
	\mathscr{B}_i^3 U=0 &\mathrm{~~in~~}\quad \Sn,\\
	\mathscr{B}_j^3 U=0 &\mathrm{~~in~~}\quad \Sn.
	\end{cases}
	\end{align}

	\begin{definition}[Biharmonic Green function]
Let $\bar f \in C^\infty(\overline{\B^{n+1}})$ and $i,j \in \{0,1,2,3\},  i+j \neq 3$. A function $\bar G^{(i,j)}:(\xi,\eta)\in \overline{\B^{n+1}}\times \overline{\B^{n+1}} \mapsto \R\cup \{\infty\}$ for \eqref{eqn:def_Green_fcn}
	is called a Green function, if $\bar G^{(i,j)}$ satisfies
	\begin{itemize}
	\item $\xi \mapsto \bar G^{(i,j)}(\xi,\eta)-\Gamma(\xi-\eta) \in C^\infty(\overline{\B^{n+1}})$ for all $\eta \in \B^{n+1}$;
	\item $\Delta_{\xi}^2 \left(\bar G^{(i,j)}(\xi,\eta)-\Gamma(\xi-\eta)\right)=0$ for all $(\xi,\eta) \in \B^{n+1}\times \B^{n+1}$;
	\item $\mathscr{B}_i^3 \bar G^{(i,j)}=0, \mathscr{B}_j^3 \bar G^{(i,j)}=0$ on $\Sn$.
	\end{itemize}
	\end{definition}
	
	Formally, the above \emph{Green function}  enables us to express a unique solution to \eqref{eqn:def_Green_fcn} as
	$$U(\xi)=\bar G^{(i,j)}\ast \bar f(\xi)=\int_{\B^{n+1}} \bar G^{(i,j)}(\xi,\eta) \bar f(\eta) \ud \eta.$$
	We notice that $\bar G^{(0,1)}$ coincides with the Boggio formula ( e.g., Boggio \cite{Boggio} or Gazzola-Grunau-Sweers \cite[Lemma 2.27]{GGS}). We next combine this and the biharmonic Poisson kernel to derive explicit formulae for other Green functions. 
		
				\subsection{Noncritical dimensions}
Recall the fundamental solution of $\Delta^2$ in $\R^{n+1}$:
\begin{align*}
	\Gamma(X-Y)=\Gamma(|X-Y|)=\frac{1}{2(n-1)(n-3)|\S^n|}|X-Y|^{3-n}.
\end{align*}

	For brevity, we let $\xi^\ast=|\xi|^{-2}\xi$ for $0\neq \xi \in \overline{\B^{n+1}}$. Notice that
	$$|\xi|^2|\xi^\ast-\eta|^2=|\xi-\eta|^2+(1-|\xi|^2)(1-|\eta|^2)=|\eta|^2|\eta^\ast-\xi|^2$$
	can be extended to define for all $\xi,\eta \in \B^{n+1}$.
	\begin{thm}\label{Thm:Green fcn_ball}
	For $n\neq 3$, the representation formulae of Green functiona $\bar G^{(i,j)}$ to \eqref{eqn:def_Green_fcn} are
	\begin{align*}
	\bar G^{(0,1)}(\xi,\eta)=&\Gamma(\xi-\eta)+\frac{1}{2}\Gamma(|\xi| \cdot |\xi^\ast-\eta|)\left[(n-3)\frac{|\xi-\eta|^2}{|\xi|^2 \cdot |\xi^\ast-\eta|^2}-(n-1)\right];\\
\bar G^{(2,3)}(\xi,\eta)=&\Gamma(\xi-\eta)-\frac{1}{2}\Gamma(|\xi| \cdot |\xi^\ast-\eta|)\left[(n-3)\frac{|\xi-\eta|^2}{|\xi|^2 \cdot |\xi^\ast-\eta|^2}-(n-1)\right];\\
\bar G^{(0,2)}(\xi,\eta)=&\Gamma(\xi-\eta)-\Gamma(|\xi|\cdot |\xi^{*}-\eta|);\\
\bar G^{(1,3)}(\xi,\eta)=&\Gamma(\xi-\eta)+\Gamma(|\xi|\cdot |\xi^{*}-\eta|);
\end{align*}
for all $\xi, \eta \in \B^{n+1}, \xi \neq \eta$, where $\xi^\ast=|\xi|^{-2}\xi$ for $0\neq \xi \in \overline{\B^4}$. Moreover, there holds
$$0\leq \bar G^{(0,1)}(\xi,\eta)\leq \bar G^{(0,2)}(\xi,\eta)\leq \bar G^{(1,3)}(\xi,\eta)\leq \bar G^{(2,3)}(\xi,\eta), \quad \forall~ \xi,\eta \in \overline{\B^{n+1}},\xi\neq \eta,$$
with equality if and only if either $\xi \in \Sn$ or $\eta \in \Sn$.
	\end{thm}

\noindent \textbf{Proof of Theorems \ref{Green-fcn:noncritical} and \ref{Thm:Green fcn_ball}.} Fix $X \in \R_+^{n+1}$. Then for all $Y \in \pa \R_+^{n+1}$, direct computations (cf. calculations in the proof of Theorem \ref{Thm:Poisson kernel_half-space}) show
\begin{align}\label{relation:Poisson-kernel_fund-sol}
\begin{split}
\mathscr{B}^3_0(\Gamma(X-Y))=&\frac{1}{2}P^3_3(x-y,t);\\
	\mathscr{B}^3_1(\Gamma(X-Y))=&-\frac{t}{2(n-1)|\S^n|(|x-y|^2+t^2)^{\frac{n-1}{2}}}=\frac{1}{2}P^3_2(x-y,t);\\
	\mathscr{B}^3_2(\Gamma(X-Y))=&\frac{1}{|\S^n|}\frac{t^2}{(t^2+|x-y|^2)^{\frac{n+1}{2}}}=-\frac{1}{2}P^3_1(x-y,t);\\
	\mathscr{B}^3_3(\Gamma(X-Y))=&-\frac{(n+1)}{|\S^n|}\frac{t^3}{(t^2+|x-y|^2)^{\frac{n+3}{2}}}=-\frac{1}{2}P^3_0(x-y,t).
	\end{split}
\end{align}
Here each $\mathscr{B}^3_k, k \in \{1,2,3\}$ acts on $\Gamma$ with respect to the variable $Y$.

By definition of Green functions on $\R^{n+1}_+$ and $\B^{n+1}$, we invoke the conformal invariance of $P_4^g$ to show
\begin{equation}\label{conf-change:Green_fcns}
\bar G^{(i,j)}(\xi,\eta)=G^{(i,j)}(X,Y)\left(\frac{2}{|\eta+\mathbf{e}_{n+1}|^2}\right)^{\frac{n-3}{2}}\left(\frac{2}{|\xi+\mathbf{e}_{n+1}|^2}\right)^{\frac{n-3}{2}}.
\end{equation}
For brevity, we let
$$H^{(i,j)}(X,Y)=G^{(i,j)}(X,Y)-\Gamma(X-Y)$$

Then $H^{(i,j)}$ solves
	\begin{align*}
\begin{cases}
\displaystyle~~~\Delta^2 H^{(i,j)}=0 \qquad &\mathrm{in}\qquad \R_+^{n+1},\\
\displaystyle \mathscr{B}^3_i(H^{(i,j)})= -\mathscr{B}^3_i(\Gamma) \qquad &\mathrm{on}\qquad \pa \R_+^{n+1},\\
\displaystyle \mathscr{B}^3_j(H^{(i,j)})=-\mathscr{B}^3_j(\Gamma)\qquad &\mathrm{on}\qquad \pa \R_+^{n+1}.
\end{cases}
\end{align*}

For $(i,j)=(0,2)$, we find 
\begin{align*}
	H^{(0,2)}(X,Y)=-\Gamma(\bar{X}-Y).
\end{align*}
This gives 
\begin{align*}
	G^{(0,2)}(X,Y)=\Gamma(X-Y)-\Gamma(\bar{X}-Y).
\end{align*}
For $(i,j)=(1,3)$, similarly we have
\begin{align*}
	G^{(1,3)}(X,Y)=\Gamma(X-Y)+\Gamma(\bar{X}-Y).
\end{align*}
For  $(i,j)=(0,1),(2,3)$, the Boggio formula (see, Boggio \cite{Boggio} or Gazzola-Grunau-Sweers \cite[Lemma 2.27]{GGS}) gives
\begin{align*}
	G^{(0,1)}(X,Y)=&\frac{1}{4|\S^n|}|X-Y|^{3-n}\int_{1}^{\frac{|\bar X-Y|}{|X-Y|}}(s^2-1)s^{-n}\ud s\\
	=&\Gamma(X-Y)+\frac{1}{2}\Gamma( |\bar X-Y|)\left[(n-3)\frac{|X-Y|^2}{ |\bar X-Y|^2}-(n-1)\right].
\end{align*}

Set $$f_{2}(y)=-\mathscr{B}^3_2(\Gamma(X-Y))=\frac{1}{2}P^3_1(x-y,t)$$
and
\begin{align*}
f_{3}(y)=-\mathscr{B}^3_3(\Gamma(X-Y))=\frac{1}{2}P^3_0(x-y,t).
\end{align*}
Thanks to Theorem \ref{Thm:Poisson kernel_half-space} and \eqref{relation:Poisson-kernel_fund-sol}, a clever observation is
\begin{align*}
H^{(2,3)}(X,Y)=&\int_{\R^n} P^3_2(y-z,s)f_{2}(z)\ud z+\int_{\R^n}P^3_3(y-z,s)f_{3}(z)\ud z\\
=&\frac{1}{2}\int_{\R^n}P^3_2(y-z,s)P^3_1(x-z,t)\ud z+\frac{1}{2}\int_{\R^n}P^3_3(y-z,s)P^3_0(x-z,t)\ud z\\
=&\int_{\R^n}P^3_0(x-z,t)\mathscr{B}^3_0\big( \Gamma(y-z,s)\big)\ud z+\int_{\R^n}P^3_1(x-z,t)\mathscr{B}^3_1\big( \Gamma(y-z,s)\big)\ud z\\
=&-H^{(0,1)}(X,Y).
\end{align*}
Then it follows 
\begin{align*}
G^{(2,3)}(X,Y)=&\Gamma(X-Y)-\frac{1}{2}\Gamma(|\bar X-Y| )\left[(n-3)\frac{|X-Y|^2}{ |\bar X-Y|^2}-(n-1)\right].
\end{align*}

Observe that
$$|\bar X-Y|^2=|X-Y|^2+4st, \qquad \forall~ X,Y \in \overline{\R_+^{n+1}}.$$
Then we have
\begin{align*}
 G^{(2,3)}(X,Y)\geq \Gamma(X-Y)+\Gamma(|\bar X-Y|)
=& G^{(1,3)}(X,Y)\geq G^{(0,2)}(X,Y) \geq 0,\\
[ G^{(0,2)}-G^{(0,1)}](X,Y)=&\frac{n-3}{2}\Gamma(|\bar X-Y|)\left(1-\frac{|X-Y|^2}{|\bar X-Y|^2}\right)\geq 0.
\end{align*}
This finish the proof of Theorem \ref{Green-fcn:noncritical}.

\medskip

Let $\xi^\ast= |\xi|^{-2} \xi$ whenever $\xi \neq 0$. Noitice that $F(\bar X)=\xi^\ast$.
Using 
$$|X+\mathbf{e}_{n+1}| |\xi+\mathbf{e}_{n+1}|=2, \quad |\bar{X}+\mathbf{e}_{n+1}| |\xi^\ast+\mathbf{e}_{n+1}|=2,$$
we obtain
\begin{align*}
	|\xi|^2=\frac{|\bar{X}+\mathbf{e}_{n+1}|^2}{|X+\mathbf{e}_{n+1}|^2}=\frac{ |\xi+\mathbf{e}_{n+1}|^2}{|\xi^\ast+\mathbf{e}_{n+1}|^2}.
\end{align*}
Also, the following elementary identity
\begin{align*}
	2|\xi^{*}-\eta|=|\xi^\ast+\mathbf{e}_{n+1}||\eta+\mathbf{e}_{n+1}||\bar X-Y|
\end{align*}
gives
\begin{align}\label{formula:dist-fcns_two_models2}
	2|\xi||\xi^{*}-\eta|=|\xi+\mathbf{e}_{n+1}||\eta+\mathbf{e}_{n+1}||\bar{X}-Y|.
\end{align}

Therefore, using \eqref{conf-change:Green_fcns}, \eqref{formula:dist-fcns_two_models} and \eqref{formula:dist-fcns_two_models2} the explicit representation formulae $\bar G^{(i,j)}$ for $\B^{n+1}$ follow from those of $ G^{(i,j)}$ for $\R_+^{n+1}$.

\medskip

Notice that
$$|\xi|^2 |\xi^\ast-\eta|^2=|\xi-\eta|^2+(1-|\xi|^2)(1-|\eta|^2), \qquad \forall~ \xi,\eta \in \overline{\B^{n+1}}.$$
Then 
\begin{align*}
\bar G^{(2,3)}(\xi,\eta)\geq& \Gamma(\xi-\eta)+\Gamma(|\xi|\cdot |\xi^{*}-\eta|)\\
=&\bar G^{(1,3)}(\xi,\eta)\geq\bar G^{(0,2)}(\xi,\eta) \geq 0.
\end{align*}
On the other hand, we have
\begin{align*}
[\bar G^{(0,2)}-\bar G^{(0,1)}](\xi,\eta)=\frac{n-3}{2}\Gamma(|\xi| \cdot |\xi^\ast-\eta|)\left(1-\frac{|\xi-\eta|^2}{|\xi|^2 \cdot |\xi^\ast-\eta|^2}\right)\geq 0.
\end{align*}
Hence the second assertion follows. This completes the proof of Theorem \ref{Thm:Green fcn_ball}.
\hfill $\Box$

\subsection{Critical dimension}

The purpose is to derive explicit formulae of biharmonic Poisson kernels for $\R_+^4$ and $\B^4$.

Recall that the fundamental solution for $\Delta^2$ on $\R_+^4$ is 
\begin{align*}
\Gamma(X-Y)=-\frac{1}{4|\S^3|}\log|X-Y|.
\end{align*}

	\begin{thm}\label{Thm:Green fcn_ball_dim-4}
	The representation formula of Green function $\bar G^{(i,j)}$ to \eqref{eqn:def_Green_fcn} is
	\begin{align*}
		\bar{G}^{(0,1)}(\xi,\eta)=&\Gamma(\xi-\eta)-\Gamma(|\xi|\cdot |\xi^{*}-\eta|)+\frac{1}{8|\S^3|}\left(\frac{|\xi-\eta|^2}{|\xi|^2|\xi^{*}-\eta|^2}-1\right);\\
		\bar{G}^{(2,3)}(\xi,\eta)=&\Gamma(\xi-\eta)+\Gamma(|\xi| \cdot |\xi^{*}-\eta|)-\frac{1}{8|\S^3|}\left(\frac{|\xi-\eta|^2}{|\xi|^2|\xi^{*}-\eta|^2}-1\right)+C;\\
			\bar{G}^{(0,2)}(\xi,\eta)=&\Gamma(\xi-\eta)-\Gamma(|\xi|\cdot |\xi^{*}-\eta|);\\
			\bar{G}^{(1,3)}(\xi,\eta)=&\Gamma(\xi-\eta)+\Gamma(|\xi|\cdot |\xi^{*}-\eta|)
			+C;
	\end{align*}
	for $\xi,\eta \in \overline{\B^4}, \xi \neq \eta$ and $C \in \R$.
\end{thm}

\noindent \textbf{Proof of Theorems \ref{Thm:Green fcn_half-space_dim-4} and \ref{Thm:Green fcn_ball_dim-4}.}
By definition of biharmonic Green function, we have
\begin{align*}
	H^{(i,j)}(X,Y)=G^{(i,j)}(X,Y)-\Gamma(X-Y)
\end{align*}
satisfy
	\begin{align*}
\begin{cases}
\displaystyle\Delta^2 H^{(i,j)}=0 \qquad &\mathrm{~~in~~}\qquad ~~\R^{4}_+,\\
\displaystyle \mathscr{B}^3_i(H^{(i,j)})= -\mathscr{B}^3_i(\Gamma) \qquad &\mathrm{~~on~~}\qquad \partial\R^4_+,\\
\displaystyle \mathscr{B}^3_j(H^{(i,j)})=-\mathscr{B}^3_j(\Gamma)\qquad &\mathrm{~~on~~}\qquad \partial\R^4_+.
\end{cases}
\end{align*}

For $(i,j)=(0,2)$, we find
\begin{align*}
	H^{(0,2)}(X,Y)=-\Gamma(\bar{X}-Y)
\end{align*}
whence,
\begin{align*}
	G^{(0,2)}(X,Y)=\Gamma(X-Y)-\Gamma(\bar{X}-Y).
\end{align*}

By \eqref{Poisson_Fund-sol} and symmetry we find
\begin{align*}
G^{(1,3)}(X,Y)=\Gamma(X-Y)+\Gamma(|\bar X -Y|)+C, \quad C \in \R.
\end{align*}

For $(i,j)=(0,1), (2,3)$, again the Boggio formula gives
\begin{align*}
	G^{(0,1)}(X,Y)=&\frac{1}{4|\S^3|}\int_{1}^{\frac{|\bar X-Y|}{|X-Y|}}(s^2-1)s^{-3}\ud s\\
	=&\Gamma(X-Y)-\Gamma(|\bar X-Y|)+\frac{1}{8|\S^3|}\left(\frac{|X-Y|^2}{|\bar X-Y|^2}-1\right).
\end{align*}
Notice that
\begin{align*}
	H^{(2,3)}(X,Y)=&\frac{1}{2}\int_{\R^3}P^3_2(y-z,s)P^3_1(x-z,t)\ud z+\frac{1}{2}\int_{\R^3}(P^3_3(y-z,s)+c)P^3_0(x-z,t)\ud z\\
	=&\frac{1}{2}\int_{\R^3}P^3_2(y-z,s)P^3_1(x-z,t)\ud z+\frac{1}{2}\int_{\R^3}P^3_3(y-z,s)P^3_0(x-z,t)\ud z+\frac{C}{2}.
\end{align*}
Then
\begin{align*}
	H^{(2,3)}(X,Y)=-H^{(0,1)}(X,Y)+\frac{C}{2}.
\end{align*}
Consequently, we have
\begin{align*}
	G^{(2,3)}(X,Y)=&\Gamma(X-Y)+\Gamma(|\bar X-Y|)-\frac{1}{8|\S^3|}\left[\frac{|X-Y|^2}{|\bar X-Y|^2}-1\right]+C.
\end{align*}
This finishes the proof of Theorem \ref{Thm:Green fcn_half-space_dim-4}.


Next, we derive the Green function on $\B^4$ associated with the unique smooth solution of
\begin{align*}
\begin{cases}
\displaystyle\Delta^2 U=f \qquad \quad &\mathrm{in}\qquad \B^4,\\
\displaystyle \mathscr{B}^3_i U=0 \qquad &\mathrm{on}\qquad \S^3,\\
\displaystyle \mathscr{B}^3_j U=0 \qquad &\mathrm{on}\qquad \S^3,
\end{cases}
\end{align*}
where $\bar f \in C^\infty(\overline{\B^4})$. We claim that
\begin{align*}
	\bar{G}^{(i,j)}(\xi,\eta)=G^{(i,j)}(X,Y).
\end{align*}
To this end, by conformal covariance \eqref{conf_transformation:Paneitz} and \eqref{conf_transformation:bdry_operators} of $P_4^g$ and $\mathscr{B}_k^3$ we have
\begin{align*}
	\Delta^2 (U \circ F) =e^{4\hat U_0}P_4^{e^{2\hat U_0}|\ud X|^2}(U\circ F)=e^{4\hat U_0}[P_4^{|\ud \xi|^2}(U)]\circ F=e^{4\hat U_0}\bar f\circ F
	\end{align*}
	and for $k \in \{0,1,2,3\}$,
	\begin{align*}
	(\mathscr{B}_k^3)_{|\ud X|^2}(U\circ F)=e^{k \hat U_0}(\mathscr{B}_k^3)_{e^{2\hat U_0}|\ud X|^2}(U\circ F)=e^{k \hat U_0}[(\mathscr{B}_k^3)_{|\ud \xi|^2} U] \circ F=0.
	\end{align*}
Hence, we are able to apply the Green formula on $\R^4_+$ for $U \circ F$ to obtain
\begin{align*}
U(\xi)=U\circ F (X)=&\int_{\R_+^4}G^{(i,j)}(X,Y) e^{4\hat U_0(Y)}\bar f\circ F(Y) \ud Y\\
=&\int_{\B^4} G^{(i,j)}(F^{-1}(\xi),F^{-1}(\eta)) \bar f(\eta) \ud \eta,
\end{align*}
which in turn gives
$$\bar{G}^{(i,j)}(\xi,\eta)=G^{(i,j)}(F^{-1}(\xi),F^{-1}(\eta)).$$
This completes the proof of Theorem \ref{Thm:Green fcn_ball_dim-4}.
\hfill $\Box$

\medskip
Combining the biharmonic Poisson kernel and  Green function, we are in a position to present a Green formula for an nonhomogeneous biharmonic BVP on $\B^{n+1}$.
	\begin{cor}[Biharmonic Poisson and Green integral]\label{cor:Green_formula}
	For $i,j \in \{0,1,2,3\}$ and $i+j \neq 3$, suppose $\bar f \in C^\infty(\overline{\B^{n+1}}), \bar f_i, \bar f_j \in C^\infty(\Sn)$. Then the unique solution of 
	\begin{align*}
	\begin{cases}
	\Delta^2 U=\bar f\qquad &\mathrm{~~in~~}\quad \B^{n+1},\\
	\mathscr{B}_i^3 U=\bar f_i &\mathrm{~~in~~}\quad \Sn,\\
	\mathscr{B}_j^3 U=\bar f_j &\mathrm{~~in~~}\quad \Sn,
	\end{cases}
	\end{align*}
	has the form of
	$$U(\xi)=[\bar G^{(i,j)}\ast \bar f+\bar P_i^3\ast \bar f_i+\bar P_j^3 \ast \bar f_j](\xi).$$
	\end{cor}
	As an application of the Green formula, we are able to obtain the following comparison principles for biharmonic BVPs on $\B^{n+1}$.
	\begin{prop}\label{Prop:Comparison_principle}
	Suppose $n \geq 4$ and  $U$ is a smooth function in $\overline{\B^{n+1}}$ satisfying
	\begin{align*}
	\begin{cases}
	\displaystyle\quad \Delta^2 U\geq 0 \qquad &\mathrm{~~in~~}\qquad \B^{n+1},\\
	\displaystyle \mathscr{B}^3_i(U)\geq 0 \qquad &\mathrm{~~on~~}\qquad \S^n,\\
	\displaystyle \mathscr{B}^3_j(U)\leq 0\qquad &\mathrm{~~on~~}\qquad \S^n,
	\end{cases}
	\end{align*}
	for $i\in \{0,3\}, j \in \{1,2\}$, then $U\geq 0$ in $\B^{n+1}$.
\end{prop}
\begin{pf}
This is a direct consequence of Corollary \ref{cor:Green_formula},  Theorems \ref{Thm:Poisson kernel_ball} and \ref{Thm:Green fcn_ball}.
\end{pf}

\section{ Classification theorems}\label{Sect5}
Our approach  builds on the intrinsic third-order GJMS operator $P_3=\frac{1}{2}\mathscr{B}_3^3$ in our setting, and  the basics of the geometric bubble and biharmonic Poisson integral. 

An important property of the geometric bubble $U_{x_0,\ve}$ is that for $n \geq 4$ and $i,j \in \{1,2,3\}, i<j$,  $U_{x_0,\ve}$ satisfies the BVP \eqref{Intro main eqn} with critical growth that $p_i=p_i^\ast, p_j=p_j^\ast$. 

There are two ways to see this. One is through a straightforward computation. The other is in a geometric way of exhausting all possible geometric bubbles: A conformal transformation on $\overline{\B^{n+1}}$ can be expressed as
$$\zeta:=\psi_{a}(\xi)=\frac{\xi-a-|\xi|^2 a+2(a \cdot \xi)a-|a|^2 \xi}{1-2 a \cdot \xi+|a|^2|\xi|^2}, \quad \forall~ \xi \in \B^{n+1}$$
for some $a \in \B^{n+1}$ (see e.g., \cite[Chapter 1]{Hua}), and there holds
$$\psi_a^\ast(|\ud \zeta|^2)=\left(\frac{1-|a|^2}{|\xi|^2 |a|^2-2a \cdot \xi+1}\right)^2 |\ud \xi|^2.$$
Indeed, one can verify that
$$F^\ast (\psi_a^\ast(|\ud \zeta|^2))=U_{x_0,\ve}^{\frac{4}{n-3}}(X) |\ud X|^2$$
through the relation
$$a=F(x_0,\ve) \qquad \mathrm{by~~regarding~~} (x_0,\ve)\in \R_+^{n+1}.$$
Equivalently, the BVP \eqref{Intro main eqn} with critical growth can interpreted as the pullback of the metric $(\B^{n+1},\psi_a^\ast(|\ud \zeta|^2))$ via the map $F$ has constant curvatures $(T_i^3)_{\Sn}, i\in \{1,2,3\}$.

As before,  let $P_k(x)$ denote a polynomial with degree $\leq k$ for $k \in \N$.

\subsection{Proof of Theorem \ref{Main thm 1}}

	\begin{lem}\label{Boundary data lem}
			Suppose  $u$  is a nonnegative smooth function in $\overline{\R_+^{n+1}}$ satisfying
		\begin{align*}
		\begin{cases}
		\displaystyle \Delta^2 u=0\qquad\qquad\qquad\qquad\, &\mathrm{in}\qquad \R_{+}^{n+1},\\
		\displaystyle \mathscr{B}_3^3(u)=\mathbb{T}^3_3 u^{p_3}\qquad\,\,\,\,\,\,\,&\mathrm{on}\quad\,\,\,\, \partial\R^{n+1}_{+},
		\end{cases}
		\end{align*}
		for $0<p_3\leq p_3^\ast=\frac{n+3}{n-3}$.
	Assume that
	\begin{align*}
		u(x,0)=O(|x|^{-c}),\qquad \partial_t u(x,0)=O(|x|^{-d})	\qquad \mathrm{~~as~~} |x| \to \infty
		\end{align*}
	for some constants $c>3/p_3$ and $d>1$. Then we have
	\begin{enumerate}[(i)]
		\item If $p_3=p_3^\ast$, then there exist $\ve>0$ and $x_0\in \R^n$ such that 
		\begin{align*}
		u(x,0)=\left(\frac{2\ve}{\ve^2+|x-x_0|^2}\right)^{\frac{n-3}{2}}.
		\end{align*}
		\item If $0<p_3<p_3^\ast$, then $u(x,0)=0$.
		\end{enumerate}
	\end{lem}
	\begin{pf}
		Notice that $f_1(x):=-\partial_tu(x,0)=O(|x|^{-d})$ and $f_3(x):=\mathbb{T}^3_3 u^{p_3}(x,0)=O(|x|^{-p_3c})$ as $|x| \to \infty$. We let
		\begin{align*}
			v(x,t)=-\frac{2}{|\Sn|}\int_{\Rn}\frac{t^2 f_1(y)}{(t^2+|x-y|^2)^{\frac{n+1}{2}}} \ud y+\frac{1}{(n-3)(n-1)|\S^n|}\int_{\R^n}\frac{f_3(y)}{(t^2+|x-y|^2)^{\frac{n-3}{2}}}.
		\end{align*}
		By Lemmas \ref{singular integral estimate lem 1} and \ref{singular integral estimate lem 2}, as $|X| \to \infty$ we have
		\begin{align}\label{bound:v}
			|v(x,t)|\lesssim t\cdot \begin{cases}
			|X|^{-\min\{d, n\}} &\mathrm{if~~} d\neq n\\
			|X|^{-n}\log |X| &\mathrm{if~~} d= n
			\end{cases}
			+ \begin{cases}
			|X|^{3-\min\{p_3c,n\}} &\mathrm{if~~} p_3 c \neq n\\
			|X|^{3-n} \log |X| &\mathrm{if~~} p_3 c = n
			\end{cases}
			\to 0.
		\end{align}

		By Theorem \ref{Thm:Poisson kernel_half-space}, $w:=u-v$ solves
		\begin{align*}
		\begin{cases}
		\displaystyle  ~\Delta^2 w=0  &\mathrm{~~in~~}\qquad~~ \R_{+}^{n+1},\\
		\displaystyle 	~~~ \partial_{t}w=0 &\mathrm{~~on~~}\qquad \partial\R^{n+1}_{+},\\
		\displaystyle\partial_{t} \Delta w=0\qquad\qquad &\mathrm{~~on~~}\qquad \partial\R^{n+1}_{+}.
		\end{cases}
		\end{align*}
		Moreover, since $u$ is nonnegative, by \eqref{bound:v} we have
		\begin{align*}
	w(x,t)\geq -C|X|^{-\delta} \qquad  \mathrm{~~as~~}|X|\to \infty,
		\end{align*}
	where $\delta=\min\{\min\{d, n\}-1,\min\{p_3c,n\}-3\}/2>0$, and 
		\begin{align*}
			 w(x,0)=& u(x,0)-\frac{1}{(n-3)(n-1)|\S^n|}\int_{\R^n}\frac{f_3(y)}{|x-y|^{n-3}}\ud y\\
			 =&O(|x|^{\frac{3-\min\{p_3c,n\}}{2}})+O(|x|^{-c})\to 0, \qquad \quad \mathrm{~~as~~} |x|\to \infty.
		\end{align*}
	Then it follows from Lemma \ref{Uniqueness lem} that $w(x,t) =c_2't^2$ for some nonnegative constant $c_2'$. This means
		\begin{align}\label{formula:u_global}
			u(x,t)=\frac{n+1}{4 |\S^n|}\int_{\R^n}\frac{u^{p_3}(y,0)}{(t^2+|x-y|^2)^{\frac{n-3}{2}}}\ud y+\frac{2}{|\Sn|}\int_{\Rn}\frac{t^2 \pa_t u(y,0)}{(t^2+|x-y|^2)^{\frac{n+1}{2}}} \ud y+c_2't^2.
		\end{align}
		Letting $t\to 0$ we obtain
			\begin{align*}
		u(x,0)=\frac{n+1}{4|\S^n|}\int_{\R^n}\frac{u^{p_3}(y,0)}{|x-y|^{n-3}}\ud y.
		\end{align*}
			The classification theorem (e.g., Y. Y. Li \cite{Li} and Chen-Li-Ou \cite{Chen&Li&Ou}) for  the integral equation above gives 
		\begin{align*}
		u(x,0)=\left(\frac{2\ve}{\ve^2+|x-x_0|^2}\right)^{\frac{n-3}{2}},\quad \ve \in \R_+, x_0 \in \Rn \qquad \qquad \mathrm{~~for~~} p_3=p_3^\ast
		\end{align*}
		and
		\begin{align*}
		u(x,0)=0 \qquad \qquad \mathrm{~~for~~} 0<p_3<p_3^\ast.
		\end{align*}
	\end{pf}

	\noindent \textbf{Proof of Theorem \ref{Main thm 1}.}
		For $c>\max\{3/p_3,1/p_1\}$, $u(x,0)=O(|x|^{-c})$ and $\partial_tu(x,0)=-\frac{n-3}{2}u^{p_1}(x,0)=O(|x|^{-p_1c})$ as $|x| \to \infty$.

			\emph{(1)} If $p_3=p_3^\ast=\frac{n+3}{n-3}$, then by Lemma \ref{Boundary data lem} we have 
		\begin{align*}
			u(x,0)=\left(\frac{2\ve}{\ve^2+|x-x_0|^2}\right)^{\frac{n-3}{2}}.
		\end{align*}
		
		 Now $u$ solves
			\begin{align}\label{equ tpye 1}
		\begin{cases}
		\displaystyle \Delta^2 u=0 &\mathrm{~~in~~}\quad~~ \R_{+}^{n+1},\\
		\displaystyle ~~\partial_{t}u=-\frac{n-3}{2}\Big(\frac{2\ve}{\ve^2+|x-x_0|^2}\Big)^{\frac{(n-3)p_1}{2}} &\mathrm{~~on~~}\quad \partial\R^{n+1}_{+},\\
		\displaystyle \mathscr{B}_3^3 u=\frac{(n^2-1)(n-3)}{4}\Big(\frac{2\ve}{\ve^2+|x-x_0|^2}\Big)^{\frac{n+3}{2}} \qquad &\mathrm{~~on~~}\quad\partial\R^{n+1}_{+}.
		\end{cases}
		\end{align}
	By \eqref{formula:u_global} we know
	\begin{align}\label{u:formula_I}
	u(x,t)=M_1(x,t)+c_2't^2,
	\end{align}
	where $c_2'$ is a nonnegative constant and
	\begin{align*}
	M_1(x,t)=&\frac{n+1}{4|\S^n|}\int_{\R^n}\frac{1}{(t^2+|x-y|^2)^{\frac{n-3}{2}}}\Big(\frac{2\ve}{\ve^2+|y-x_0|^2}\Big)^{\frac{n+3}{2}}\ud y\\
	&-\frac{n-3}{|\S^n|}\int_{\R^n}\frac{t^2}{(t^2+|x-y|^2)^{\frac{n+1}{2}}}\Big(\frac{2\ve}{\ve^2+|y-x_0|^2}\Big)^{\frac{(n-3)p_1}{2}} \ud y.
	\end{align*}

		Especially, assume $p_1=p_1^\ast=\frac{n-1}{n-3}$.  On one hand, by Lemmas  \ref{singular integral estimate lem 1} and \ref{singular integral estimate lem 2} we have
	\begin{equation}\label{Decay_rate_M}
	M_1(x,t)=O(|X|^{3-n}) \qquad \mathrm{~~as~~} |X| \to \infty.
	\end{equation}

		On the other hand, we consider 
		\begin{align*}
			w(x,t)=u(x,t)-U_{x_0,\ve}(x,t),
		\end{align*}
		then $w$ satisfies 
		\begin{align}\label{thm 1 formula a}
		\begin{cases}
		\displaystyle ~\Delta^2 w=0\qquad\qquad\qquad &\mathrm{in}\qquad \R_{+}^{n+1},\\
		\displaystyle 	~~~ \partial_{t}w=0&\mathrm{on}\quad\,\,\,\, \partial\R^{n+1}_{+},\\
		\displaystyle\partial_{t} \Delta w=0&\mathrm{on}\quad\,\,\,\, \partial\R^{n+1}_{+},
		\end{cases}
		\end{align}
		and  for $|X|,|x| \gg 1$,
		\begin{align*}
			w(X)\geq -C|X|^{3-n} \qquad \mathrm{and}\qquad w(x,0)=O(|x|^{-\min\{c,n-3\}}).
		\end{align*}
		Hence, it follows from Lemma \ref{Uniqueness lem} that  the solution of \eqref{equ tpye 1} has the form of 
		\begin{align}\label{u:formula_II}
			u(x,t)=U_{x_0,\ve}(x,t)+c_2t^2.
		\end{align}
		
		Therefore, we combine \eqref{u:formula_I}, \eqref{Decay_rate_M} and \eqref{u:formula_II} to conclude that $M_1=U_{x_0,\ve}$ and $c_2'=c_2$. 
	
	\medskip
	
	\emph{(2)} 	If $p_3<p_3^\ast$, then again by Lemma \ref{Boundary data lem} we obtain  $u(x,0)=0$. In this case, $u$ also solves the same homogeneous biharmonic BVP \eqref{thm 1 formula a}, thanks to Lemma \ref{Uniqueness lem}, we conclude that $u(x,t)=c_2t^2$ for every nonnegative constant $c_2$. 
\hfill $\Box$

	\subsection{Proof of Theorem \ref{Main thm 2}}
		
From the viewpoint of the \emph{Lopatinski-Shapiro} or \emph{complement condition}, the conformal boundary operator pair $(\mathscr{B}_1^3,\mathscr{B}_2^3)$ is exceptional. Nevertheless, some different technique needs to be developed for treating this type of $T$-curvature equations. 
	\medskip
	
		\noindent \textbf{Proof of Theorem \ref{Main thm 2}.}
	For $cp_2>2$,  $f_0(x)=u(x,0)=O(|x|^{-c})$ and $f_2(x)=\mathbb{T}^3_2 u^{p_2}=O(|x|^{-p_2c})$ as $|x| \to \infty$. Let
	\begin{align*}
		v(x,t)=\frac{2(n+1)}{|\S^n|}\int_{\R^n}\frac{t^3 f_0(y)}{(t^2+|x-y|^2)^{\frac{n+3}{2}}} \ud y-\frac{1}{(n-1)|\S^n|}\int_{\R^n}\frac{tf_2(y)}{(t^2+|x-y|^2)^{\frac{n-1}{2}}} \ud y
	\end{align*}
	with
	\begin{align}\label{Thm 2 formula a}
	|v(x,t)|\lesssim
	\begin{cases}
	|X|^{-\min\{c, n\}} &\mathrm{if~~} c \neq n\\
	|X|^{-n}\log |X| &\mathrm{if~~} c= n
	\end{cases}
	+
	 \begin{cases}
	|X|^{2-\min\{p_2c,n\}} &\mathrm{if~~} p_2c \neq n\\
	|X|^{2-n}\log |X| &\mathrm{if~~} p_2c= n
	\end{cases}
	\to 0
	\end{align}
	as $|X| \to \infty$ by virtue of Lemmas \ref{singular integral estimate lem 1} and  \ref{singular integral estimate lem 2}.
	
		By the assumption and Theorem \ref{Thm:Poisson kernel_half-space}, $w:=u-v$ satisfies
	\begin{align*}
	\begin{cases}
	\displaystyle \Delta^2 w=0\qquad\qquad\qquad &\mathrm{in}\qquad \R_{+}^{n+1},\\
	\displaystyle  \quad~	w=0&\mathrm{on}\quad\,\,\,\, \partial\R^{n+1}_{+},\\
	\displaystyle ~~\Delta w=0&\mathrm{on}\quad\,\,\,\, \partial\R^{n+1}_{+}.
	\end{cases}
	\end{align*}
	Moreover, using the nonnegativity of $u$ and \eqref{Thm 2 formula a} we have
	\begin{align*}
	w(x,t) \geq -C|X|^{-\delta}\qquad\mathrm{as}\qquad |X|\to \infty,
	\end{align*}
	where $\delta=\min\left\{\min\{p_2c,n\}-2,\min\{c, n\}\right\}/2>0$.
	Then we apply Lemma \ref{Uniqueness lem 2} to obtain 
	\begin{align}\label{Proof of thm 2 integral}
	u(x,t)=&\frac{2(n+1)}{|\S^n|}\int_{\R^n}\frac{t^3u(y,0)}{(t^2+|x-y|^2)^{\frac{n+3}{2}}}\ud y-\frac{n-3}{2|\S^n|}\int_{\R^n}\frac{tu^{p_2}(y,0)}{(t^2+|x-y|^2)^{\frac{n-1}{2}}} \ud y\no\\
	&+tP_2(x)+c_3 t^3,
	\end{align}
	where  $c_3$ is a nonnegative constant and 
	$\liminf_{|x|\to\infty} P_2(x)\geq 0.$
	
From Lemma \ref{lem:singular_integral_II} that
\begin{align*}
	\lim_{t\to 0}\partial _t\int_{\R^n}\frac{t^3 u(y,0)}{(t^2+|x-y|^2)^{\frac{n+3}{2}}}\ud y=0,
	\end{align*}
together with the first boundary condition we obtain
	\begin{align}\label{Proof of thm 2 formula a}
	\frac{n-3}{2}u^{p_1}(x,0)=-\partial_tu(x,0)=\frac{n-3}{2|\S^n|}\int_{\R^n}\frac{u^{p_2}(y,0)}{|x-y|^{n-1}}\ud y-P_2(x).
	\end{align}
Using the assumption that $u(x,0)=o(1)$  as $|x| \to \infty$ and Lemma \ref{singular integral estimate lem 1} we have
	\begin{align*}
		P_2(x)=o(1) \qquad \mathrm{as}\qquad |x|\to+\infty,
	\end{align*}
	so $P_2(x)=0$.

	Write $k(x)=u^{p_1}(x,0)$,  the integral equation \eqref{Proof of thm 2 formula a} becomes
	\begin{align*}
	k(x)=\frac{1}{|\S^n|}\int_{\R^n}\frac{k^{\frac{p_2}{p_1}}(y)}{|x-y|^{n-1}}\ud y.
	\end{align*}
	
		\emph{(1)} If $\frac{p_2}{p_1}=\frac{n+1}{n-1}$, then  the classification theorem (e.g.,  \cite{Chen&Li&Ou} and \cite{Li}) yields
	\begin{align*}
	k(x)=\left(\frac{2\ve}{\ve^2+|x-x_0|^2}\right)^{\frac{n-1}{2}}, \qquad \ve \in \R_+,~ x_0 \in \Rn.
	\end{align*}
	Meanwhile, the second boundary condition becomes
	\begin{align*}
		\mathscr{B}^3_2 u=\frac{(n-1)(n-3)}{2}\left(\frac{2\ve}{\ve^2+|x-x_0|^2}\right)^{\frac{n+1}{2}}.
	\end{align*}
 Going back to \eqref{Proof of thm 2 integral} we conclude that
\begin{align*}
u(x,t)=M_2(x,t)+c_3 t^3,
\end{align*}
where
	\begin{align*}
M_2(x,t)
=&\frac{2(n+1)}{|\S^n|}\int_{\R^n}\frac{t^3}{(t^2+|x-y|^2)^{\frac{n+3}{2}}}\left(\frac{2\ve}{\ve^2+|y-x_0|^2}\right)^{\frac{n-1}{2p_1}}\ud y\nonumber\\
&-\frac{n-3}{2|\S^n|}\int_{\R^n}\frac{t}{(t^2+|x-y|^2)^{\frac{n-1}{2}}} \left(\frac{2\ve}{\ve^2+|y-x_0|^2}\right)^{\frac{n+1}{2}}\ud y.
\end{align*}

 Especially, when $p_1=p_1^\ast=\frac{n-1}{n-3}$, we consider
 \begin{align*}
 w(x,t)=u(x,t)-U_{x_0,\ve}(x,t).
 \end{align*}
Then $w$ solves
 \begin{align}\label{Vanishing boundary type 2}
 \begin{cases}
 \displaystyle \Delta^2 w=0\qquad\qquad\qquad &\mathrm{in}\qquad \R_{+}^{n+1},\\
 \displaystyle 	\quad ~w=0&\mathrm{on}\quad\,\,\,\, \partial\R^{n+1}_{+},\\
 \displaystyle  ~~\Delta w=0 &\mathrm{on}\quad\,\,\,\, \partial\R^{n+1}_{+},
 \end{cases}
 \end{align}
with some $C \in \R_+$ such that
 \begin{align*}
 w(X)\geq -C|X|^{3-n} \qquad \mathrm{as~~} |X| \to \infty.
 \end{align*}
 Moreover, by Lemmas \ref{singular integral estimate lem 1} and \ref{singular integral estimate lem 2} we have
 $$M_2(x,t)=O(|X|^{-\min\{\frac{n-1}{p_1},n\}})+O(|X|^{2-n})\to 0 \qquad \mathrm{as~~} |X| \to \infty.$$
 Hence, by Lemma \ref{Uniqueness lem 2} we obtain
	\begin{align*}
	u(x,t)=U_{x_0,\ve}(x,t)+c_3t^3.
	\end{align*}

		\emph{(2)} If $\frac{p_2}{p_1}<\frac{n-1}{n-3}$, then the uniqueness theorem in  Y. Y. Li \cite{Li} yields $u(x,0)=0$ in $\Rn$, so $u$ satisfies problem \eqref{Vanishing boundary type 2}. This together with Lemma \ref{Uniqueness lem 2} yields $u(x,t)=c_3t^3$ in $\R_+^{n+1}$ for every nonnegative constant $c_3$. 
		\hfill $\Box$

	\subsection{Proof of Theorem \ref{Main thm 3}}
	
	\noindent \textbf{Proof of Theorem \ref{Main thm 3}.}
	Recall that $u$ solves
	\begin{align}\label{Type III equ}
\begin{cases}
\displaystyle ~~~\Delta^2 u=0  &\mathrm{~~in~~}\quad~ \R_{+}^{n+1},\\
\displaystyle \mathscr{B}_2^3(u)=\mathbb{T}^3_2 u^{p_2}  &\mathrm{~~on~~}\quad \partial\R^{n+1}_{+},\\
\displaystyle \mathscr{B}_3^3(u)=\mathbb{T}^3_3 u^{p_3} \qquad \qquad&\mathrm{~~on~~}\quad \partial\R^{n+1}_{+}.
\end{cases}
\end{align}
 As $cp_2>2$, a similar argument as in Theorem \ref{Main thm 2} yields
\begin{align}\label{formula:u_III}
u(x,t)=&\frac{2(n+1)}{|\Sn|}\int_{\R^n}\frac{t^3u(y,0)}{(t^2+|x-y|^2)^{\frac{n+3}{2}}}\ud y-\frac{n-3}{2|\Sn|}\int_{\R^n}\frac{tu^{p_2}(y,0)}{(t^2+|x-y|^2)^{\frac{n-1}{2}}} \ud y\no\\
	&+tP_2(x)+c_3 t^3,
\end{align}
where $c_3$ is a nonnegative constant and 
\begin{align}\label{Sec 2.3 thm formula a1}
	\liminf_{|x|\to\infty} P_2(x)\geq 0.
\end{align}
Moreover, by \eqref{formula:u_III} and Lemma \ref{lem:singular_integral_II} we have
\begin{align*}
-\partial_tu(x,0)=\frac{n-3}{2|\S^n|}\int_{\R^n}\frac{u^{p_2}(y,0)}{|x-y|^{n-1}}\ud y-P_2(x).
\end{align*}

We claim that there is some nonnegative constant $c$ such that 
\begin{align*}\label{Sec 2.3 claim}
	P_2(x)\equiv c \qquad \mathrm{~~and~~}\qquad c_3=0.
\end{align*}

 To this end, \eqref{Sec 2.3 thm formula a1} yields that $\overline{\Delta}P_2(x)$ must be a nonnegative constant. With $\tilde c_3:=c_3+\frac{1}{2}\overline{\Delta}P_2(x)$ we introduce
 \begin{align}\label{modified_u}
&\tilde{u}(x,t)=u(x,t)-P_2(x)t-c_3 t^3+\tilde c_3 t^3\no\\
=&\frac{2(n+1)}{|\S^n|}\int_{\R^n}\frac{t^3u(y,0)}{(t^2+|x-y|^2)^{\frac{n+3}{2}}}\ud y-\frac{n-3}{2|\S^n|}\int_{\R^n}\frac{tu^{p_2}(y,0)}{(t^2+|x-y|^2)^{\frac{n-1}{2}}} \ud y+\tilde c_3 t^3.
\end{align}
Clearly, $\tilde u$ is biharmonic in $\R_+^{n+1}$ and $\tilde{u}(x,0)=u(x,0)$.
Recall that $\mathscr{B}^3_2=\partial_t^2-\overline{\Delta}$ and $\mathscr{B}^3_3=\partial_t\Delta +2\overline{\Delta}\partial_t$, then $\tilde u$ and $u$ satisfy the same boundary conditions, more explicitly,
\begin{align*}\label{Sec 2.3 thm formula b0}
	\mathscr{B}^3_3(\tilde{u})=\mathscr{B}^3_3(u)+6(\tilde c_3-c_3)-3\overline{\Delta}P_2(x)=\mathscr{B}^3_3(u)
\end{align*}
and  $\mathscr{B}^3_2(\tilde{u})=\mathscr{B}^3_2(u)$ on $\pa \R_+^{n+1}$.
This means that $\tilde u$ also satisfies  problem \eqref{Type III equ},  though it might change sign in $\R_+^{n+1}$.

On the other hand, notice that there is some $C \in \R_+$ such that $\tilde{u}(x,t)\geq -C|X|^{3}$ as $|X| \to \infty$. By Lemma \ref{singular integral estimate lem 1} we have
\begin{align*}
	-\partial_t	\tilde{u}(x,0)=&\frac{n-3}{2|\S^n|}\int_{\R^n}\frac{\tilde{u}^{p_2}(y,0)}{|x-y|^{n-1}}\ud y\\
	=&\begin{cases}
	O(|x|^{1-\min\{p_2c,n\}}) &\qquad \mathrm{if~~} p_2c \neq n\\
	O(|x|^{1-n}\log |x|)&\qquad \mathrm{if~~} p_2c= n
	\end{cases}\\
	=&O(|x|^{-d}) \qquad\qquad \mathrm{~~as~~} |x| \to \infty,
\end{align*}
where $d=\min\{c p_2-1,(n-1)/2\}>1$.
As $c>\max\{2/p_2,3/p_3\}$ by assumption, we apply  Remark \ref{Rem:Uniqueness} to obtain 
\begin{align}\label{modified:u_II}
&\tilde{u}(x,t)\no\\
=&\frac{2}{|\S^n|}\int_{\R^n}\frac{t^2 \partial_t\tilde{u}(y,0)}{(t^2+|x-y|^2)^{\frac{n+1}{2}}} \ud y+\frac{n+1}{4 |\S^n|}\int_{\R^n}\frac{u^{p_3}(y,0)}{(t^2+|x-y|^2)^{\frac{n-3}{2}}}\ud y+t^2P_1(x).
\end{align} 

Combining \eqref{modified_u} and \eqref{modified:u_II}, we can apply Lemmas \ref{singular integral estimate lem 1} and \ref{singular integral estimate lem 2} to conclude that
\begin{align*}
	\tilde c_3t^3-t^2P_1(x)=o(1) \qquad \mathrm{~~as~~} |X| \to \infty.
\end{align*}
This implies $\tilde c_3=0$ and $P_1(x)\equiv 0$. Thus, \eqref{modified:u_II} becomes
\begin{align}\label{Sec 2.3 thm formula b}
\tilde{u}(x,t)=\frac{2}{|\S^n|}\int_{\R^n}\frac{t^2 \partial_t\tilde{u}(y,0)}{(t^2+|x-y|^2)^{\frac{n+1}{2}}} \ud y+\frac{n+1}{4 |\S^n|}\int_{\R^n}\frac{u^{p_3}(y,0)}{(t^2+|x-y|^2)^{\frac{n-3}{2}}}\ud y.
\end{align}

Keep in mind that 
\begin{align*}
	0=\tilde c_3=c_3+\frac{1}{2}\overline{\Delta}P_2(x),
\end{align*}
so, $c_3=0$ and $\overline{\Delta}P_2(x)=0$. This together with \eqref{Sec 2.3 thm formula a1} implies $P_2(x)\equiv c_1\geq0$. Hence, we prove the above claim.
\medskip

Going back to \eqref{formula:u_III} and  \eqref{modified_u} we have
\begin{align}\label{formula_u_final}
u(x,t)=&\frac{2(n+1)}{|\S^n|}\int_{\R^n}\frac{t^3u(y,0)}{(t^2+|x-y|^2)^{\frac{n+3}{2}}}\ud y-\frac{n-3}{2|\S^n|}\int_{\R^n}\frac{tu^{p_2}(y,0)}{(t^2+|x-y|^2)^{\frac{n-1}{2}}} \ud y\no\\
	&+c_1 t
	\end{align}
	and there is some $C\in \R_+$ such that
\begin{align}\label{Sec 2.3 thm formula e}
	\tilde{u}(x,t)=u(x,t)-c_1 t\geq -C|X| \qquad \mathrm{as~~} |X| \to \infty.
\end{align}
Furthermore, letting $t \to 0$ in  \eqref{Sec 2.3 thm formula b}, by Lemma \ref{singular integral estimate lem 1} we arrive at
\begin{align*}
u(x,0)=\tilde{u}(x,0)=\frac{n+1}{4|\S^n|}\int_{\R^n}\frac{u^{p_3}(y,0)}{|x-y|^{n-3}}\ud y.
\end{align*}

\emph{(1)}  If $p_3=\frac{n+3}{n-3}$ and $p_2>\frac{2}{n-3}$, then the classification theorem (e.g.,  \cite{Li} and \cite{Chen&Li&Ou}) gives
\begin{align*}
u(x,0)=\left(\frac{2\ve}{\ve^2+|x-x_0|^2}\right)^{\frac{n-3}{2}}, \quad x_0 \in \Rn, \ve \in \R_+.
\end{align*}
This together with \eqref{formula_u_final} yields
$$u(x,t)=M_3(x,t)+c_1t,$$
where
\begin{align*}
M_3(x,t)
			=&\frac{2(n+1)}{|\S^n|}\int_{\R^n}\frac{t^3}{(t^2+|x-y|^2)^{\frac{n+3}{2}}} \left(\frac{2\ve}{\ve^2+|y-x_0|^2}\right)^{\frac{n-3}{2}}\ud y\\
			&-\frac{n-3}{2|\S^n|}\int_{\R^n}\frac{t}{(t^2+|x-y|^2)^{\frac{n-1}{2}}} \left(\frac{2\ve}{\ve^2+|y-x_0|^2}\right)^{\frac{(n-3)p_2}{2}}\ud y.
			\end{align*}

In particular, when $p_2=p_2^\ast$,  notice that there is some $C \in \R_+$ such that  $(\tilde u-U_{x_0,\ve})(X)\geq -C |X|$ as $|X| \to \infty$. Then, we apply Remark \ref{Uniqueness lem 2 rem 2}  to conclude that $M_3=U_{x_0,\ve}$ and then
\begin{align*}
	u(x,t)=U_{x_0,\ve}(x,t)+c_1t.
\end{align*}

	\emph{(2)}  If $p_3<\frac{n+3}{n-3}$,  then the uniqueness theorem (e.g., \cite{Li}) yields $u(x,0)=0$. This together with \eqref{formula_u_final}  follows that $u(x,t)=c_1t$ in $\R_+^{n+1}$ for every nonnegative constant $c_1$.
\hfill $\Box$

\section{Applications to biharmonic boundary value problems on balls}\label{Sect6}

Although we focus on the conformal geometry nature for biharmonic BVPs, we would  like to mention that from a physical viewpoint, the boundary conditions $(\mathscr{B}_0^3,\mathscr{B}_1^3)$ and $(\mathscr{B}_0^3,\mathscr{B}_2^3)$ have a close connection with the clamped plate and hinged plate model (cf. \cite{GGS}), respectively. However, the study of biharmonic equation coupled with other boundary conditions is far less developed.

\subsection{Solutions with single boundary point singularity}

\noindent \textbf{Proof of Theorem \ref{Thm:bdry_singularity}.}
It suffices to give the proof for $i=1$ and $j=3$, since the others are the same in the spirit thanks to Theorems \ref{Main thm 2} and \ref{Main thm 3}.

Let $U$ be a nonnegative solution of \eqref{bdry_value_problem_singular pt}. As before, we introduce
$$u(X)=\left(\frac{2}{|X+\mathbf{e}_{n+1}|^2}\right)^{\frac{n-3}{2}} (U\circ F)(X).$$
Then $u$ solves
\begin{align*}
\begin{cases}
\displaystyle \Delta^2 u=0\qquad\qquad\qquad\qquad\, &\mathrm{in}\qquad \R_+^{n+1},\\
\displaystyle \mathscr{B}^3_1(u)=\mathbb{T}^3_1u^{\frac{n-1}{n-3}} &\mathrm{on}\quad\,\,\,\, \partial\R_+^{n+1},\\
\displaystyle \mathscr{B}^3_3(u)=\mathbb{T}^3_3 u^{\frac{n+3}{n-3}}&\mathrm{on}\quad\,\,\,\, \partial\R_+^{n+1}.
\end{cases}
\end{align*}
By definition of the conformal map $F$, we know that $|X+\mathbf{e}_{n+1}||\xi+\mathbf{e}_{n+1}|=2$ and $F(\pa \R_+^{n+1})=\pa \B^{n+1}$. By the assumption \eqref{assump:isolated_singularity} for $U$, we obtain 
\begin{align*}
u(x,0)=\left(\frac{2}{1+|x|^2}\right)^{\frac{n-3}{2}}U\circ F(x,0)\lesssim |x|^{3-n}|X+\mathbf{e}_{n+1}|^{\alpha}\lesssim |x|^{-(n-3-\alpha)}, \quad \mathrm{as~~} |x| \to \infty
\end{align*}\label{sol:zero_second}
with $n-3-\alpha>\frac{3(n-3)}{n+3}$. By Theorem \ref{Main thm 1}, we have
\begin{align*}
u(x,t)=U_{x_0,\ve}(x,t)+c_2 t^2
\end{align*}
for some $x_0\in \Rn, \ve \in \R_+$ and nonnegative constant $c_2$. 

If we let $\xi_0=F(x_0,\ve) \in \B^{n+1}$, then it is not hard to verify that
$$F^\ast(U_{\xi_0}(\xi)^{\frac{4}{n-3}} |\ud \xi|^2)=U_{x_0,\ve}(X)^{\frac{4}{n-3}}|\ud X|^2,$$
where
$$U_{\xi_0}(\xi)=\left(\frac{1-|\xi_0|^2}{|\xi|^2 |\xi_0|^2-2 \xi_0 \cdot \xi+1}\right)^{\frac{n-3}{2}}, \qquad \xi \in \B^{n+1}.$$
In other words, 
$$U_{\xi_0}(\xi)=\left(\frac{2}{|\xi+\mathbf{e}_{n+1}|^2}\right)^{\frac{n-3}{2}}U_{x_0,\ve}\circ F^{-1}(\xi).$$
Using 
$$t\circ F^{-1}=\frac{1-|\xi|^2}{|\xi+\mathbf{e}_{n+1}|^2}$$
and \eqref{sol:zero_second} we obtain
\begin{align*}
U(\xi)=&\left(\frac{2}{|\xi+\mathbf{e}_{n+1}|^2}\right)^{\frac{n-3}{2}} u \circ F^{-1}(\xi)\\
=&U_{\xi_0}(\xi)+c_2 \left(\frac{2}{|\xi+\mathbf{e}_{n+1}|^2}\right)^{\frac{n-3}{2}}\left(\frac{1-|\xi|^2}{|\xi+\mathbf{e}_{n+1}|^2}\right)^2.
\end{align*}
This follows the desired assertion.
\hfill $\Box$\\

\subsection{Sharp geometric inequalities}

On a smooth compact manifold with boundary, the maximizing problem of the isometric ratio over a scalar-flat conformal class was studied in \cite{Hang-Wang-Yan1,Hang-Wang-Yan2}. With comparison principles in Proposition \ref{Prop:Comparison_principle}, we are motivated to study of the maximizers of the isometric ratio   and  another sharp geometric inequality over some $Q$-curvature flat conformal class of the flat metric  $[|\ud \xi|^2]$ on $\overline{\B^{n+1}}$.

More general, for $n\geq 4$ we define
\begin{align*}
\mathcal{M}_1=&\{g\in [|\ud \xi|^2]; Q_g\leq 0 \mathrm{~~in~~} \B^{n+1},~ (T_2^3)_g\geq 0 \mathrm{~~on~~} \Sn\};\\
\mathcal{M}_2=&\{g\in [|\ud \xi|^2]; Q_g\leq 0 \mathrm{~~in~~} \B^{n+1},~  h_g\geq 0 \mathrm{~~on~~} \Sn\}
\end{align*}
and
\begin{align*}
\mathring{\mathcal{M}}:=&\{g\in [|\ud \xi|^2]; Q_g= 0 \mathrm{~~in~~} \B^{n+1},~ (T_2^3)_g=0 \mathrm{~~on~~} \Sn\}\\
=&\{g\in [|\ud \xi|^2]; Q_g=0 \mathrm{~~in~~} \B^{n+1},~  h_g= 0 \mathrm{~~on~~} \Sn\}
\end{align*}
by virtue of Proposition \ref{Prop:ex-intrinsic_GJMS}.

\begin{lem}\label{lem:lifting}
	Let  $g=U^{\frac{4}{n-3}} |\ud \xi|^2$ be a smooth metric in $\overline{\B^{n+1}}$. Suppose  $g \in \mathcal{M}_1\cup \mathcal{M}_2$. Then we have 
	\begin{enumerate}[(1)]
	\item there exists a conformal metric  $\tilde g=\tilde U^\frac{4}{n-3} |\ud \xi|^2 \in \mathring{\mathcal{M}}$ such that  $\tilde{U}\geq U$ in $\B^{n+1}$ with $\tilde{U}=U$ on $\Sn$;
	\item there exists a conformal metric  $\tilde g=\tilde U^\frac{4}{n-3} |\ud \xi|^2 \in \mathring{\mathcal{M}}$ such that  $\tilde{U}\geq U$ in $\B^{n+1}$ with $\mathscr{B}^3_3(U)=\mathscr{B}^3_3(\tilde U)$ on $\Sn$.
	\end{enumerate}
\end{lem}

\begin{pf}
	For \emph{(i)}, we introduce an auxiliary function by	
	\begin{align*}
	\tilde{U}(\xi)=&\frac{(n+1)\left(1-|\xi|^2\right)^3}{4|\Sn|}\int_{\Sn}\frac{U(\eta)}{|\xi-\eta|^{n+3}}\ud V_{\Sn}(\eta).
	\end{align*}
	By Theorem \ref{Thm:Poisson kernel_ball} we have 
	\begin{align*}
	\begin{cases}
	\displaystyle  \quad \Delta^2 U\leq 0=\Delta^2 \tilde U \qquad &\mathrm{in}\qquad \B^{n+1},\\
	\displaystyle  \qquad~ U=\tilde U \qquad &\mathrm{on}\qquad \S^n,\\
	\displaystyle \mathscr{B}^3_i(U)\geq 0= \mathscr{B}^3_i(\tilde U)\qquad &\mathrm{on}\qquad \S^n,
	\end{cases}
	\end{align*}
	for $i \in \{1,2\}$.
	Then the desired assertion follows by Proposition \ref{Prop:Comparison_principle}.
	
	For \emph{(ii)}, we consider
	\begin{align*}
	\tilde{U}(\xi)=
	\frac{1}{(n-1)(n-3)|\S^n|}\int_{\S^n} \frac{\mathscr{B}^3_3(U) (\eta)}{|\xi-\eta|^{n-3}}  \ud V_{\S^n}(\eta).
	\end{align*}
	By Theorem \ref{Thm:Poisson kernel_ball} we have  
	\begin{align*}
	\begin{cases}
	\displaystyle \quad \Delta^2 U\leq 0=\Delta^2 \tilde U \qquad &\mathrm{in}\qquad \B^{n+1},\\
	\displaystyle \mathscr{B}^3_i(U)\geq 0=\mathscr{B}^3_i(\tilde U) \qquad &\mathrm{on}\qquad \Sn,\\
	\displaystyle \mathscr{B}^3_3(U)=\mathscr{B}^3_3(\tilde U)\qquad &\mathrm{on}\qquad \Sn,
	\end{cases}
	\end{align*}
	for $i \in \{1,2\}$. Then the desired assertion follows from Proposition \ref{Prop:Comparison_principle}.
\end{pf}

In the following, we consider the isoperimetric ratio over $\mathcal{M}_1\cup \mathcal{M}_2$:
\begin{align}\label{Isoperimetric variation problem}
\sup_{g\in\mathcal{M}_1\cup \mathcal{M}_2}\frac{|\B^{n+1}|_g}{|\S^{n}|_g^{\frac{n+1}{n}}}=\sup_{g\in \mathring{\mathcal{M}}}\frac{|\B^{n+1}|_g}{|\S^{n}|_g^{\frac{n+1}{n}}}
\end{align}
by virtue of Lemma \ref{lem:lifting}. Given each $g=U^{\frac{4}{n-3}}|\ud \xi|^2 \in \mathring{\mathcal{M}}$, notice that
\begin{align*}
\int_{\S^n}|(T_3^3)_g|^{\frac{2n}{n+3}}\ud V_g=\left(\frac{2}{n-3}\right)^{\frac{2n}{n+3}}	\int_{\S^n}|\mathscr{B}^3_{3}(u)|^{\frac{2n}{n+3}}\ud V_{\S^n}.
\end{align*}
then it follows from Theorem \ref{Thm:Poisson kernel_ball} that $(T_3^3)_g$ does not vanish identically. So, again by Lemma \ref{lem:lifting} it is reasonable to study
\begin{align}\label{Curvature integral variation problem}
\sup_{g\in\mathcal{M}_1\cup\mathcal{M}_2}	\frac{|\B^{n+1}|^{\frac{n-3}{2(n+1)}}_g}{\|(T_3^3)_g\|_{L^{\frac{2n}{n+3}}(\S^n,g)}}=	\sup_{g\in\mathring{\mathcal{M}}}	\frac{|\B^{n+1}|^{\frac{n-3}{2(n+1)}}_g}{\|(T_3^3)_g\|_{L^{\frac{2n}{n+3}}(\S^n,g)}}.
\end{align}

\begin{thm} Let $n \geq 4$. Then for all $g \in \mathcal{M}_1\cup \mathcal{M}_2$, 
\begin{align*}
|\B^{n+1}|_g \leq d_n |\Sn|_g^{\frac{n+1}{n}} \quad \mathrm{~~and~~} \quad
|\B^{n+1}|^{\frac{n-3}{2(n+1)}}_g \leq e_n \|(T_3^3)_g\|_{L^{\frac{2n}{n+3}}(\Sn,g)},
\end{align*}
where 
$$d_n=|\S^n|^{-\frac{1}{n}}\int_{0}^{1}r^n\big(1+\frac{n-3}{4}(1-r^2)\big)^{\frac{2(n+1)}{n-3}}\ud r$$
and
$$e_n=|\S^n|^{-\frac{7n+3}{2n(n+1)}}\frac{4}{(n^2-1)(n-3)}\left[\int_{0}^{1}r^n\big(1+\frac{n-3}{4}(1-r^2)\big)^{\frac{2(n+1)}{n-3}}\ud r\right]^{\frac{n-3}{2(n+1)}}.$$
Moreover,  the above two equalities hold if and only if  $g=c U_a^{4/(n-3)} |\ud \xi|^2 \in \mathring{\mathcal{M}}$ with $c \in \R_+$ and
\begin{equation}\label{extremal_fcns_on_balls}
U_a(\xi)=\Big(\frac{1-|a|^2}{|a|^2 |\xi|^2-2a \cdot \xi+1}\Big)^{\frac{n-3}{2}}+\frac{n-3}{4}(1-|\xi|^2)\Big(\frac{1-|a|^2}{|a|^2 |\xi|^2-2a \cdot \xi+1}\Big)^{\frac{n-1}{2}}
\end{equation}
for every $a \in \B^{n+1}$.
 \end{thm}
\begin{pf}
Given a conformal metric $g=U^{\frac{4}{n-3}}|\ud \xi|^2\in \mathring{\mathcal{M}}$, if we let
\begin{align*}
u(x)=\left(\frac{2}{1+|x|^2}\right)^{\frac{n-3}{2}}U\circ F(x,0),
\end{align*}
then a similar argument in the proof of Theorem \ref{Thm:bdry_singularity} yields
\begin{align*}
\left(\frac{2}{(1+t)^2+|x|^2}\right)^{\frac{n-3}{2}}U\circ F(x,t)=\frac{2(n+1)}{|\S^n|}\int_{\R^n}\frac{t^3 u(y)}{(t^2+|x-y|^2)^{\frac{n+3}{2}}}\ud y.
\end{align*}
To be consistent with notation in Gluck \cite{Gluck}, we define
$$E_{-2,3}(u)(x,t):=\left(\frac{2}{(1+t)^2+|x|^2}\right)^{\frac{n-3}{2}}U\circ F(x,t).$$
 So the extremal problem \eqref{Isoperimetric variation problem} is equivalent to 
\begin{align}\label{HLS type 1 ineq}
d_n:=\sup_{\substack{0<u \in C^\infty(\Rn)\\ u=O(|x|^{3-n})\mathrm{~as~} |x| \to \infty}}\frac{\|E_{-2,3}(u)\|_{L^{\frac{2(n+1)}{n-3}}(\R^{n+1}_{+})}}{\|u\|_{L^{\frac{2n}{n-3}}(\R^n)}}.
\end{align}
Gluck \cite[Theorem 1.1]{Gluck} states that  the optimal constant $d_n$ for the extremal problem \eqref{HLS type 1 ineq} is attained by 
\begin{align*}
		c(|x-x_0|^2+d^2)^{\frac{3-n}{2}}, \qquad c,d \in \R_+, x_0 \in \Rn,
		\end{align*}
equivalently, modulo of a positive constant, the extremal metric is a conformal metric $u_a^{4/(n-3)}g_{\Sn}$ induced by some conformal transformation on $\Sn$, which can be expressed as
\begin{equation}\label{conf-factor:trans_sphere}
u_a(\xi)=\Big(\frac{1-|a|^2}{|a|^2 |-2a \cdot \xi+1}\Big)^{\frac{n-3}{2}} \qquad \mathrm{on~~}\quad \Sn
\end{equation}
for every $a \in \B^{n+1}$. Meanwhile, the uniqueness Theorem \ref{Thm:Poisson kernel_ball} for  both $(\mathscr{B}_0^3,\mathscr{B}_1^3)$ and $(\mathscr{B}_0^3,\mathscr{B}_2^3)$ shows that
\begin{align*}
U_a(\xi)=	\frac{(n+1)\left(1-|\xi|^2\right)^3}{4|\S^n|}\int_{\Sn}\frac{u_a(\eta)}{|\xi-\eta|^{n+3}}\ud V_{\Sn}(\eta)
\end{align*}
such that $U_a^{\frac{4}{n-3}}|\ud \xi|^2 \in  \mathring{\mathcal{M}}$ is an extremal metric  with $U_a=u_a$ on $\Sn$.
Indeed, the trick in authors \cite[Section 6]{Chen-Zhang} can be applied to derive an explicit representation formula of
$$U_a(\xi)=\Big(\frac{1-|a|^2}{|a|^2 |\xi|^2-2a \cdot \xi+1}\Big)^{\frac{n-3}{2}}+\frac{n-3}{4}(1-|\xi|^2)\Big(\frac{1-|a|^2}{|a|^2 |\xi|^2-2a \cdot \xi+1}\Big)^{\frac{n-1}{2}}.$$
Moreover, we prefer to compute the exact value of $d_n$ in a geometric way: Take $1+\frac{n-3}{4}(1-|\xi|^2)$ as an extremal function, and then
\begin{align*}
d_n=&\frac{\int_{\B^{n+1}}\left(1+\frac{n-3}{4}(1-|\xi|^2)\right)^{\frac{2(n+1)}{n-3}}\ud \xi}{	|\S^{n}|^{\frac{n+1}{n}}}\\
=&|\S^n|^{-\frac{1}{n}}\int_{0}^{1}r^n\left(1+\frac{n-3}{4}(1-r^2)\right)^{\frac{2(n+1)}{n-3}}\ud r.
\end{align*}

Similarly, given a conformal metric $g=U^{\frac{4}{n-3}}|\ud \xi|^2\in \mathring{\mathcal{M}}$, if we let
\begin{align}\label{transfer-fcns:Dou-Zhu}
u(x)=\left(\frac{2}{1+|x|^2}\right)^{\frac{n+3}{2}}\mathscr{B}^3_3(U)\circ F(x,0),
\end{align}
then
\begin{align*}
&\left(\frac{2}{(1+t)^2+|x|^2}\right)^{\frac{n-3}{2}}U\circ F(x,t)\\
=&\frac{1}{(n-3)(n-1)|\S^n|}\int_{\R^n}\frac{u(y)}{(t^2+|x-y|^2)^{\frac{n-3}{2}}}\ud y.
\end{align*}
To be consistent with the notation in Dou-Zhu \cite{Dou-Zhu}, we define
$$E_4(u)(x,t):=\left(\frac{2}{(1+t)^2+|x|^2}\right)^{\frac{n-3}{2}}U\circ F(x,t).$$
So the extremal problem \eqref{Curvature integral variation problem} is equivalent to
\begin{align}\label{HLS type 2 ineq}
\sup_{\substack{0 \not \equiv u\in C^\infty(\Rn)\\ u=O(|x|^{3-n})\mathrm{~as~} |x| \to \infty}}\frac{\|E_{4}(u)\|_{L^{\frac{2(n+1)}{n-3}}(\R^{n+1}_{+})}}{\|u\|_{L^{\frac{2n}{n+3}}(\R^n)}}.
\end{align}

  Dou-Zhu \cite[Theorem 1.4]{Dou-Zhu} states that every maximizer for \eqref{HLS type 2 ineq} has the form of
$$c(|x-x_0|^2+d^2)^{-\frac{n+3}{2}}, \qquad c,d \in \R_+, x_0 \in \Rn.$$
Combining the above discussion and the authors \cite[Section 6]{Chen-Zhang}, we conclude from \eqref{transfer-fcns:Dou-Zhu} that modulo a positive constant,
$$2P_3 U=\mathscr{B}^3_3(U)=\frac{n(n-2)(n-4)}{4} u_a^{\frac{n+3}{n-3}} \quad \mathrm{~~on~~} \Sn,$$
where $u_a$ is a conformal factor induced by a conformal transformation on $\Sn$ as in \eqref{conf-factor:trans_sphere}. Since $\mathrm{ker} P_3=\{0\}$ implies $U=u_a$ on $\Sn$, the remaining part of proof is identical to that of \emph{(1)}.
\end{pf}

\section{Perspectives}\label{Sect7}

Similar to the boundary Yamabe problem (see for example Chen-Sun \cite{Chen-Sun} and Escobar \cite{Escobar1,Escobar2}),  it is natural to ask \emph{whether the round metric on  a geodesic ball $\Sigma_r^{n+1}$  with radius $r \in (0,\pi)$  in $\S^{n+1}$ is the unique conformal metric (up to a conformal diffeomorphism) of constant $Q$-curvature and constant boundary $T$-curvature under certain geometric constraints on the pair of $T$-curvature constants.}  

Consider  $\Sigma_r^{n+1}$ equipped with  the round metric $g_c:=g_{\S^{n+1}}$. Using the precise formulae of $T$-curvatures in \cite[p.296]{Case} and $Q$-curvature in \cite{Chen-Zhang}, we directly compute that $Q_{g_{\mathrm{c}}}=Q_{\S^{n+1}}=\frac{(n+1)((n+1)^2-4)}{8}$ and 
\begin{align}\label{T-curvatures_geodesic-ball}
	(T^3_1)_{g_{\mathrm{c}}}=&h_{g_{\mathrm{c}}}=\cot r, \no\\
	(T^3_2)_{g_{\mathrm{c}}}=&(n-1)\big(\frac{1}{2}+h_{g_{\mathrm{c}}}^2\big),\qquad (T^3_3)_{g_{\mathrm{c}}}=\frac{n^2-1}{2}h_{g_{\mathrm{c}}}\big(\frac{3}{2}+h_{g_{\mathrm{c}}}^2\big).
\end{align}

Notice that $\B^{n+1}$ is conformal equivalent to $\Sigma_r^{n+1}$, see for example, the first author and Sun \cite[pp.8-9]{Chen-Sun}. Based on \eqref{T-curvatures_geodesic-ball}, we can transfer the above problem into  the following conjecture.
\begin{conjecture}
    \label{Conjecture}
	For $n \geq 4$ and $i \in \{1,2\}$, let $U\in C^4(\overline{\B^{n+1}})$ be a positive smooth solution to
	\begin{align}\label{Conj:geometric_PDE_ball}
	\begin{cases}
	\displaystyle\quad \Delta^2 U= \frac{n-3}{2}Q_{\S^{n+1}}U^{p^{*}}  &\mathrm{~~in~~}\quad~~ \B^{n+1},\\
	\displaystyle \mathscr{B}^3_i(U)=c_i \mathbb{T}^3_i U^{p^{*}_i}  \qquad \qquad&\mathrm{~~on~~}\quad~~\S^{n},\\
	\displaystyle \mathscr{B}_3^3(U)=c_3 \mathbb{T}^3_3 U^{p^{*}_3} &\mathrm{~~on~~}\quad~~ \S^{n},
	\end{cases}
	\end{align} 
	with $c_k \in \R, k \in \{1,2,3\}$ and
\begin{align*}
	\begin{cases}
	\displaystyle c_3=\tfrac{3}{2}c_1+c_1^3 \qquad&\mathrm{for}\qquad i=1,\\
	\displaystyle c_2\geq \tfrac{1}{2},\quad c_3^2=(c_2-\tfrac{1}{2})(1+c_2) \qquad&\mathrm{for}\qquad i=2.
	\end{cases}
\end{align*}
Then $U$ is of the form
\begin{align}\label{sols:general}
U(\xi)=\left(\frac{2\ve (1-|\xi_0|^2)}{(\ve^2+1)(|\xi_0|^2|\xi|^2-2\xi_0 \cdot \xi+1)+(1-|\xi_0|^2)(|\xi|^2-1)}\right)^{\frac{n-3}{2}}
\end{align}
for some $\xi_0 \in \B^{n+1}$ and
\begin{align*}
c_1=\cot r, \quad c_1=\frac{\ve^2-1}{2\ve}\qquad\mathrm{~~for~~}\qquad i=1,\\
c_2=\frac{\ve^4+1}{4\ve^2}\qquad\mathrm{~~for~~}\qquad i=2.
\end{align*}
\end{conjecture} 

This conjecture can be regarded as an extension of Corollary \ref{Geom_Rigidity} to the constant $Q$-curvature and constant $T$-curvature equations on $\B^{n+1}$. Moreover, we can prove a uniqueness theorem  of positive solutions to a general 4-th order ODE including the BVP \eqref{Conj:geometric_PDE_ball} for radial solutions. See Theorem \ref{Uniqueness theorem}, whose proof is left to Appendix \ref{Append:B}. On the other hand, it is direct to verify that the radial function $\left(\frac{2\ve}{\ve^2+|\xi|^2}\right)^{\frac{n-3}{2}}$ corresponding to $\xi_0=0$ in \eqref{sols:general} satisfies \eqref{Conj:geometric_PDE_ball}. Hence we confirm the above conjecture for all radial solutions.

We hope that the biharmonic Poisson kernel and Green function will come into play in the resolution of the above conjecture.

\appendix
\section{Appendix: Proof of Proposition \ref{Prop:ex-intrinsic_GJMS}}\label{Append:A}

The combination of Theorem \ref{Thm:Poisson kernel_ball} and \cite[Theorem 1(3)]{Chen-Shi} demonstrates that the Green functions of fractional GJMS operators $P_1$ and $P_3$ are 
$$G_{P_1}(\xi,\eta)=\frac{2}{(n-1)|\S^n|}|\xi-\eta|^{1-n}, \quad G_{P_3}(\xi,\eta)=2\bar P_3^3(\xi,\eta), \qquad \xi,\eta \in \Sn.$$

\noindent \textbf{Proof of Proposition \ref{Prop:ex-intrinsic_GJMS}.}
				By \eqref{conf_change_bdry_operators} we set $\bar{f}_i=\mathscr{B}^3_iU$ and
				$$f_i(x)=\left(\frac{2}{1+|x|^2}\right)^{\frac{n+2i-3}{2}}\bar f_i\circ F\left((x,0)\right)$$
				for $i \in \{0,1,2,3\}$. 
				
				When $n\neq 3$, by Theorem \ref{Thm:Poisson kernel_ball} we have
				\begin{align*}
					U(\xi)=-&\frac{(1-|\xi|^2)^2}{2|\Sn|}\int_{\Sn} \frac{\bar{f}_1 (\eta)}{|\xi-\eta|^{n+1}}  \ud V_{\Sn}(\eta)\nonumber\\
					+&\frac{1}{(n-1)(n-3)|\S^n|}\int_{\S^n} \frac{\bar{f}_3 (\eta)}{|\xi-\eta|^{n-3}}  \ud V_{\Sn}(\eta), \qquad \xi \in \B^{n+1}.
				\end{align*}
				When $n=4$, by Theorem \ref{Thm:biharmonic_Poisson_kernerl_dim-4} we have
					\begin{align*}
&U(\xi)-\frac{1}{|\S^3|}\int_{\S^3}U   \ud V_{{\S^3}}\\
=&\int_{\S^3} \bar{P}^3_1(\xi-\eta)\bar{f}_1(\eta) \ud V_{{\S^3}}(\eta)-\frac{1}{2|\S^3|}\int_{\S^3}\log |\xi-\eta|\bar{f}_3(\eta)\ud V_{\S^3}(\eta), \quad \forall~\xi \in \overline{\B^4}.
\end{align*}

	Using the elementary formula \eqref{formula:dist-fcns_two_models} we obtain
				\begin{align*}
					\frac{(1-|\xi|^2)^2}{2|\Sn|}\int_{\Sn} \frac{\bar{f}_1 (\eta)}{|\xi-\eta|^{n+1}}  \ud V_{\Sn}(\eta)=\frac{2}{|\Sn|}\left(\frac{(1+t)^2+|x|^2}{2}\right)^{\frac{n-3}{2}}\int_{\Rn}\frac{t^2f_1(y) \ud y}{(t^2+|x-y|^2)^{\frac{n+1}{2}}} .
				\end{align*}
			Since $|\xi| \to 1$ means $t \to 0$, we have
				\begin{align*}
					\lim_{|\xi|\to 1}\frac{(1-|\xi|^2)^2}{2|\Sn|}\int_{\Sn} \frac{\bar{f}_1 (\eta)}{|\xi-\eta|^{n+1}}  \ud V_{\Sn}(\eta)=0.
				\end{align*}
				Consequently, we arrive at 
				\begin{align*}
			\forall~\xi \in \Sn,\quad		U(\xi)=\bar P_3^3 \ast \bar f_3(\xi)=\frac{1}{2}G_{P_3}\ast \bar{f}_3(\xi), \qquad \mathrm{when~~} n\neq 3
									\end{align*}
									and
									\begin{align*}
			\forall~\xi \in \Sn,\quad		U(\xi)-\frac{1}{|\S^3|}\int_{\S^3}U \ud V_{\S^3}=\bar P_3^3 \ast \bar f_3(\xi)=\frac{1}{2}G_{P_3}\ast \bar{f}_3(\xi), \qquad \mathrm{when~~} n=3.
									\end{align*}
							Consequently, we immediately obtain
									\begin{align*}
										2P_3 U=\bar{f}_3=\mathscr{B}^3_3U \qquad \mathrm{~~on~~} \Sn.
									\end{align*}

									Similarly, again by Theorem \ref{Thm:Poisson kernel_ball} we have
									\begin{align*}
										U(\xi)=&\frac{(n+1)\left(1-|\xi|^2\right)^3}{4|\Sn|}\int_{\Sn}\frac{\bar{f}_0(\eta)}{|\xi-\eta|^{n+3}}\ud V_{\Sn}(\eta)\nonumber\\
										&-\frac{(1-|\xi|^2)}{2(n-1)|\Sn|}\int_{\Sn}  \frac{\bar{f}_2(\eta)}{|\xi-\eta|^{n-1}} \ud V_{\Sn}(\eta) , \qquad \xi \in \B^{n+1}.
									\end{align*}
									Connecting a function $u$ with $U$ as in \eqref{fcns:ball_half-space} we employ \eqref{formula:dist-fcns_two_models} and conformal changes to show
									\begin{align*}
										u(x,t)=&\frac{2(n+1)}{|\S^n|}\int_{\R^n}\frac{t^3f_0(y)}{(t^2+|x-y|^2)^{\frac{n+3}{2}}}\ud y
										-\frac{1}{(n-1)|\S^n|}\int_{\R^n}\frac{tf_2(y)}{(t^2+|x-y|^2)^{\frac{n-1}{2}}} \ud y\\
										=&[P^3_0\ast f_0+P^3_2\ast f_2](x,t).
									\end{align*}
									Then it is not hard to see that
									\begin{align*}
										\lim_{t\to 0}\mathscr{B}^3_1(u)(x,t)=\frac{1}{(n-1)|\S^n|}\int_{\R^n}\frac{f_2(y)}{|x-y|^{n-1}}\ud y.
									\end{align*}
								This together with \eqref{conf_change_bdry_operators} and \eqref{formula:dist-fcns_two_models} in turn implies that $\forall~\xi \in \Sn$,
									\begin{align*}
										\mathscr{B}^3_1U(\xi)=&\left[\lim_{t \to 0}\left(\frac{(1+t)^2+|x|^2}{2}\right)^\frac{n-1}{2} \mathscr{B}^3_1(u)\right]\circ F^{-1}(\xi)\\
										=&\left[\frac{1}{(n-1)|\S^n|}\left(\frac{1+|x|^2}{2}\right)^\frac{n-1}{2}\int_{\R^n}\frac{f_2(y)}{|x-y|^{n-1}}\ud y\right]\circ F^{-1}(\xi)\\
										=&\frac{1}{(n-1)|\S^n|}\int_{\S^n}  \frac{\bar{f}_2(\eta)}{|\xi-\eta|^{n-1}} \ud V_{\Sn}(\eta)=\frac{1}{2}G_{P_1}\ast \bar{f}_2,
									\end{align*}
									which leads to
									\begin{align*}
										2P_1\mathscr{B}^3_1 U=\mathscr{B}^3_2U \qquad \mathrm{~~on~~}\quad \Sn.
									\end{align*}
			\hfill $\Box$
			
			\section{Appendix: Uniqueness of positive solutions to fourth-order ODE }\label{Append:B}

The following uniqueness theorem of the BVP for a fourth-order ODE  is of independent interest.

\begin{thm}\label{Uniqueness theorem}
Let $n \geq 4$ and $i \in \{1,2\}$. Assume that there is a positive radial solution $U \in C^\infty(\overline{\B^{n+1}})$ to 
\begin{align}\label{BVP_U}
\begin{cases}
\displaystyle \quad\Delta^2 U=U^{p^{*}}\qquad\qquad\qquad\qquad\, &\mathrm{in}\qquad \B^{n+1},\\
\displaystyle ~\mathscr{B}^3_i(U)=C_i\qquad\,\,\,\,\,\,\,&\mathrm{on}\quad\,\,\,\, \S^{n},\\
\displaystyle ~\mathscr{B}^3_3(U)=C_3\qquad\,\,\,\,\,\,\,&\mathrm{on}\quad\,\,\,\, \S^{n}.
\end{cases}
\end{align}
 Here, $C_k\in \R,k \in \{1,2,3\}$.
Then $U$ is unique.
\end{thm}
\begin{pf}
Denote by $U=C_0>0$ on $\Sn$.
Let $r=|\xi| \in (0,1]$ and $U=U(r)$. Under the change of variable $t=-\log r \in [0,\infty)$, the conformal metric is
\begin{align*}
U(r)^{\frac{4}{n-3}}|\ud \xi|^2=U(e^{-t})^{\frac{4}{n-3}}e^{-2t} (\ud t^2+g_{\Sn}):=V(t)^{\frac{4}{n-3}} g_{\mathrm{cyl}},
\end{align*}
where $ g_{\mathrm{cyl}}=\ud t^2+g_{\Sn}$ is the product metric on the cylinder $[0,\infty)\times \Sn$ and
\begin{align}\label{change:fcns}
	V(t)=e^{\frac{(3-n)t}{2}}U(e^{-t}) ,~~t\in [0,\infty) \qquad\Longleftrightarrow\qquad U(r)=r^{\frac{3-n}{2}}V(-\log r),~~ r\in (0,1].
\end{align}
Then the equation $\Delta^2U=U^{p^{*}}$ is geometrically interpreted as
\begin{align*}
1=P_4^{U(r)^{\frac{4}{n-3}}|\ud \xi|^2}(1)
=&U(r)^{-\frac{n+5}{n-3}}P_4^{|\ud \xi|^2} (U)
=U(e^{-t})^{-\frac{n+5}{n-3}} P_4^{e^{-2t}g_{\mathrm{cyl}}} (U(e^{-t}))\\
=&U(e^{-t})^{-\frac{n+5}{n-3}} e^{\frac{n+5}{2}t} P_4^{g_{\mathrm{cyl}}} (e^{\frac{3-n}{2}t}U(e^{-t}))\\
=&V(t)^{-\frac{n+5}{n-3}}P_4^{g_{\mathrm{cyl}}}(V),
\end{align*}
equivalently, 
$$P_4^{g_{\mathrm{cyl}}}(V)=V^{\frac{n+5}{n-3}}.$$
Notice that
$$Q_{g_{\mathrm{cyl}}}=\frac{(n+1)^2(n-3)}{8}, \quad P_4^{g_{\mathrm{cyl}}}(V)=\left(\pa_t^4-\frac{(n-1)^2+4}{2}\pa_t^2+\frac{n-3}{2}Q_{g_{\mathrm{cyl}}}\right)V.$$
So, $V=V(t)$ satisfies
\begin{align}\label{eqn:V}
	V^{(4)}-\frac{n^2-2n+5}{2}V''+\left(\frac{(n+1)(n-3)}{4}\right)^2V=V^{p^{*}}, \quad t\in[0,+\infty).
\end{align}
Similarly, by the conformal invariance of $\mathscr{B}^3_k, k \in \{1,2,3\}$ and the fact that $U(1)=V(0)$, the same argument yields
$$\mathscr{B}^3_k (V) =C_k \qquad \mathrm{~~on~~}\quad \{0\}\times \Sn.$$

Since the boundary $\{0\}\times \Sn$ is totally geodesic with respect to $g_{\mathrm{cyl}}$, it is easy to compute
\begin{align*}
\mathscr{B}^3_1(V)=&-\pa_t V;\\
\mathscr{B}^3_2(V)=&\Big(\pa_t^2+\frac{(n+1)(n-3)}{4}\Big)V;\\
\mathscr{B}^3_3(V)=&\Big(\pa_t^3-\frac{3n^2-6n+7}{4}\pa_t\Big)V;
\end{align*}
using the explicit formulae $\mathscr{B}^3_k, k \in \{1,2,3\}$ on manifolds with boundary; see \cite[pp.295-296]{Case}.

    Hence, collecting these facts together we conclude that
\begin{align}\label{BVP:V}
\begin{cases}
\displaystyle V^{(4)}-\frac{n^2-2n+5}{2}V''+\left(\frac{(n+1)(n-3)}{4}\right)^2V(t)=V^{p^{*}} & \mathrm{~~in~~} [0,+\infty),\\
\displaystyle -V'(0)=C_1, \quad V'''(0)-\frac{3n^2-6n+7}{4} V'(0)=C_3 &\mathrm{~~for~~} i=1,\\
 \displaystyle V''(0)+\frac{(n-3)(n+1)}{4}V(0)=C_2,\quad  V'''(0)-\frac{3n^2-6n+7}{4} V'(0)=C_3  & \mathrm{~~for~~} i=2,
\end{cases}
\end{align}
and $V(0)=C_0$.

Alternatively, problem \eqref{BVP:V} of $V$ can be derived by direct calculations. We leave it as an exercise for the readers.

 \medskip

 It reduces to showing the uniqueness of positive solutions to \eqref{BVP:V}. 
 
 In the following, we adapt the strategy of Frank-K\"onig \cite[Proposition 4]{Frank&Konig} or van den Berg \cite[Theorem 1]{van den Berg Jan Bouwe} to the uniqueness result of \eqref{BVP:V}. However, the 4-order ODE studied in \cite{van den Berg Jan Bouwe} or \cite{Frank&Konig} is defined in $\R$, and their proof  hinges on the symmetry of $\R$, so it is impossible to  apply their results directly to our setting.
 
 Suppose not, there exists another positive smooth solution $W=W(t)$ to \eqref{BVP:V} with $W(0)=C_0$, which arises from another distinct positive radially symmetric solution to \eqref{BVP_U} through the change of function \eqref{change:fcns}.
Rewrite equation \eqref{eqn:V} of $V$ as
$$(\frac{\ud^2}{\ud t^2}-\mu)(\frac{\ud^2}{\ud t^2}-\lambda)V=V^{p^\ast}$$
with
\begin{align*}
\lambda=\big(\frac{n-3}{2}\big)^2,\quad \mu=\big(\frac{n+1}{2}\big)^2.
\end{align*}

We introduce an auxiliary function by
\begin{align*}
\varphi(t)=(V-W)''(t)-\lambda (V(t)-W(t)).
\end{align*}
Notice that $V(0)=W(0)$.

For $i=1$, the boundary conditions give
\begin{align*}
V'(0)=W'(0),~~ V'''(0)=W'''(0).
\end{align*}
Without loss of generality, we assume $V''(0)>W''(0)$.
This implies
\begin{align}\label{initial_cond_1}
\varphi(0)>0, \quad \varphi'(0)=0.
\end{align}
For $i=2$,  the boundary conditions yield that
\begin{align*}
 V''(0)=&W''(0),~~ \\
V'''(0)-\frac{3n^2-6n+7}{4}V'(0)=&W'''(0)-\frac{3n^2-6n+7}{4}W'(0).\end{align*}
Without loss of generality, we assume $V'(0)>W'(0)$.
This yields
\begin{align}\label{initial_cond_2}
\varphi(0)=0,\quad 
\varphi'(0)=(V-W)'''(0)-\lambda (V-W)'(0)
=\frac{n^2-1}{2}(V-W)^{'}(0)>0.
\end{align}

In both cases, we deduce that  $V(t)>W(t)$ near zero. Then we define 
\begin{align*}
\sigma:=\sup\{s; V(t)>W(t), t\in (0,s)\} \quad \mathrm{~~and~~} \quad \sigma\in \R_+\cup\{+\infty\}.
\end{align*}
Thus we deduce from the equations of $V$ and  $W$ that 
\begin{align*}
V^{p^{*}}-W^{p^{*}}=&\left(\frac{\ud^2}{\ud t^2}-\mu\right)\left(\frac{\ud^2}{\ud t^2}-\lambda\right)(V-W)
=\varphi''-\mu\varphi.
\end{align*}
This together with the definition of $\sigma$ yields 
\begin{align}\label{ineq:auxiliary_fcn}
\varphi''(t)-\mu\varphi(t)>0 \qquad\mathrm{~~in~~}\quad (0,\sigma).
\end{align}

Next we claim that 
\begin{align*}
\varphi(t)>0 \qquad\mathrm{in}\qquad (0,\sigma).
\end{align*}

To this end, it follows from either \eqref{initial_cond_1} or \eqref{initial_cond_2} that $\varphi(t)>0$ near zero.  Then we define 
\begin{align*}
    t_0:=\sup\{T\in (0,\sigma); \varphi(t)>0, t\in (0,T)\}.
\end{align*} 
It suffices to show $t_0=\sigma$. Otherwise,
if $t_0<\sigma$, then $\varphi(t_0)=0$ by definition of $t_0$. Then we multiply both sides of \eqref{ineq:auxiliary_fcn} by $\varphi$ and integrate by parts over $(0,t_0)$ to show
\begin{align*}
-\int_{0}^{t_0}(\varphi'(s))^2\ud s=\int_{0}^{t_0}\varphi''(s)\varphi(s)\ud s>\mu\int_{0}^{t}\varphi^2(s)\ud s,
\end{align*}
where we also used either \eqref{initial_cond_1} or \eqref{initial_cond_2}. A contradiction! This finishes the proof of the above claim. 

\medskip
The above claim asserts that
\begin{align*}
(V-W)''-\lambda(V-W)>0 \qquad\mathrm{in}\qquad (0,\sigma).\end{align*}
Set $\Psi=V-W$ and then $\varphi=\Psi''-\lambda \Psi$. Then the above inequality becomes 
\begin{equation}\label{ineq:Psi}
\Psi''-\lambda \Psi>0 \qquad \mathrm{~~in~~}\quad  (0,\sigma),
\end{equation}
equipped with the boundary conditions that
\begin{align*}
\Psi(0)=\Psi'(0)=0 \qquad \mathrm{for}\qquad i=1,\\
\Psi(0)=0, \Psi'(0)>0 \qquad \mathrm{for}\qquad i=2.
\end{align*}
By definition of $\sigma$, we know that $\Psi(t)>0$ in $(0,\sigma)$. For any fixed $t \in (0,\sigma)$, we multiply both sides of \eqref{ineq:Psi} by $\Psi$ and integrate by parts to show
\begin{align*}
\Psi'(t)\Psi(t)>\int_{0}^{t}(\Psi'(s))^2\ud s+\lambda\int_{0}^{t}(\Psi(s))^2\ud s.
\end{align*}
This leads to 
\begin{align*}
\Psi(t)>0,\quad\Psi'(t)>0, \quad \Psi''(t)>0\qquad\mathrm{~~in~~}\quad (0,\sigma).
\end{align*}
Therefore, we conclude that $\sigma=+\infty$. Otherwise, by definition of $\sigma$ we have $\Psi(\sigma)=0$, and arrive at a contradiction.

\medskip

For $\sigma=+\infty$, we repeat the same argument to show 
\begin{align*}
\Psi(t)>0,\qquad\Psi'(t)>0, \qquad \Psi''(t)>0\qquad\mathrm{~~in~~}\quad (0,+\infty).
\end{align*}
This implies $\lim_{t \to +\infty} \Psi(t)=+\infty$.
On the other hand, we deduce from the change of function in \eqref{change:fcns}  that $V,W$ are bounded in $[0,\infty)$. A contradiction! 
\end{pf}

  		\bibliographystyle{unsrt}

\begin{thebibliography}{40}
		
		\bibitem{ADN}
		S. Agmon, A. Douglis and  L. Nirenberg, Estimates near the boundary for solutions of elliptic partial differential equations satisfying general boundary conditions. I. \textit{Comm. Pure Appl. Math.} \textbf{12} (1959), 623-727.
		
		\bibitem{BGW}
		E. Berchio, F. Gazzola and T. Weth, 
Radial symmetry of positive solutions to nonlinear polyharmonic Dirichlet problems.
\textit{J. Reine Angew. Math.} \textbf{620} (2008), 165-183.

\bibitem{van den Berg Jan Bouwe}  
          J. B. van den Berg, The phase-plane picture for a class of fourth-order conservative differential equations. \textit{J. Differential Equations} \textbf{161} (2000), no. 1, 110-153.
		
		\bibitem{Boggio}
			T. Boggio, Sulle funzioni di Green d'ordine $m$. \textit{Rend. Circ. Mat. Palermo}, \textbf{20}, 97-135, 1905.
			
			\bibitem{Branson-Gover}
			 T. Branson and A. Gover, Conformally invariant non-local operators.
\textit{Pacific J. Math.} \textbf{201} (2001), no. 1, 19-60.

           \bibitem{BFGM}
D. Buoso, C. Falc\'o, M. Gonz\'alez and M. Miranda, Bulk-boundary eigenvalues for bilaplacian problems.
\textit{Discrete Contin. Dyn. Syst.} \textbf{43} (2023), no. 3-4, 1175-1200.
			
							
			\bibitem{Caffarelli&Gidas&Spruck} L. Caffarelli, B. Gidas and J. Spruck,
			 Asymptotic symmetry and local behavior of semilinear elliptic equations with critical Sobolev growth.
			 \textit{Comm. Pure Appl. Math.} \textbf{42} (1989), no.3, 271--297.
						
			\bibitem{CLMNY}
			J. Case, Y. Lin, S. McKeown, C. Ndiaye and P. Yang, A fourth-order Cherrier-Escobar problem with prescribed corner behavior on the half-ball. \href{https://arxiv.org/abs/2401.00937}{arXiv:2401.00937}.

			
			\bibitem{Case}
			J. Case, Boundary operators associated with the Paneitz operator.
			\textit{Indiana Univ. Math. J.} \textbf{67}  (2018), no. 1, 293-327.
			
			\bibitem{Chang-Qing}
  S.-Y. A. Chang and J. Qing, The zeta functional determinants on manifolds with boundary. I. The formula, \textit{ J. Funct. Anal.} \textbf{147} (1997), no. 2, 327-362.
  
                      \bibitem{CSF}
  M. Chipot, I. Shafrir and M. Fila, 
On the solutions to some elliptic equations with nonlinear Neumann boundary conditions.
\textit{Adv. Differential Equations} \textbf{1} (1996), no.1, 91-110.
			
						
			\bibitem{Chen&Li&Ou} W. Chen, C. Li and B. Ou, 
				Classification of solutions for an integral equation.
			\textit{Comm. Pure Appl. Math.} \textbf{59}  (2006), no. 3, 330-343.
			
			\bibitem{Chen-Shi}
			X. Chen and Y. Shi, Green functions for GJMS operators on spheres, Gegenbauer polynomials and rigidity theorems. \href{https://arxiv.org/abs/2401.02087}{arXiv:2401.02087}.
			
			\bibitem{Chen-Sun} 
			X. Chen and L. Sun, 
				Existence of conformal metrics with constant scalar curvature and constant boundary mean curvature on compact manifolds.
					\textit{Commun. Contemp. Math.}\textbf{21} (2019), no. 3, 1850021, 51 pp.

			
			\bibitem{Chen-Zhang}
			X. Chen and S. Zhang, A sharp Sobolev trace inequality of order four on three-balls.  \href{https://arxiv.org/abs/2403.00380}{arXiv:2403.00380}. 
			
  
                 \bibitem{Dou-Zhu} J. Dou and M. Zhu, 
		Sharp Hardy-Littlewood-Sobolev inequality on the upper half space.
		\textit{Int. Math. Res. Not. IMRN}   (2015), no. 3, 651-687.
			
		\bibitem{Escobar1}
 J. Escobar, Sharp constant in a Sobolev trace inequality, \textit{Indiana Univ. Math. J.} \textbf{37} (1988), no. 3, 687-698.
			\bibitem{Escobar2}
			J. Escobar, Uniqueness theorems on conformal deformation of metrics, Sobolev inequalities, and an eigenvalue estimate. \textit{Comm. Pure Appl. Math.} \textbf{43} (1990), no. 7, 857-883.
			
			\bibitem{Frank&Konig}  	R. Frank and T. K\"onig, 
					Classification of positive singular solutions to a nonlinear biharmonic equation with critical exponent.
					\textit{Anal. PDE} \textbf{12} (2019), no. 4, 1101-1113.
							
		\bibitem{GGS}
 F. Gazzola, H. Grunau and G. Sweers, Polyharmonic boundary value problems. Positivity preserving and nonlinear higher order elliptic equations in bounded domains. \textit{Lecture Notes in Mathematics}, 1991. Springer-Verlag, Berlin, 2010. xviii+423 pp. 
 
                \bibitem{GNN}
			B. Gidas, W. Ni and L. Nirenberg, Symmetry and related properties via the maximum principle.
\textit{Comm. Math. Phys.} \textit{68} (1979), no. 3, 209-243.

		\bibitem{Gluck} M. Gluck, 
		Subcritical approach to conformally invariant extension operators on the upper half space.
		\textit{J. Funct. Anal.} \textbf{278} (2020), no. 1, 108082, 46 pp.
		
		\bibitem{GS}
			M. Gonz\'alez and M. S\'aez, 
			Eigenvalue bounds for the Paneitz operator and its associated third-order boundary operator on locally conformally flat manifolds. \textit{Ann. Sc. Norm. Super. Pisa Cl. Sci.} \textbf{(5) 24} (2023), no. 4, 2047-2081.
			
			 \bibitem{Gover-Peterson}
			A. Gover and L. Peterson, Conformal boundary operators,  $T$-curvatures, and conformal fractional Laplacians of odd order. \textit{Pacific J. Math.} \textbf{311} (2021), no. 2, 277-328.

		
		\bibitem{Grant}
			D. Grant, A conformally invariant third order Neumann-type operator for hypersurfaces, Master's thesis, The University of Auckland, 2003.
			
			\bibitem{Hang-Wang-Yan2}
		F. Hang, X. Wang and X. Yan, An integral equation in conformal geometry.
\textit{Ann. I. H. Poincar\'{e} Anal. Non lin\'{e}aires} \textbf{26} (2009), no. 1, 1-21. 
			
			\bibitem{Hang-Wang-Yan1} F. Hang, X. Wang and X. Yan, 
	Sharp integral inequalities for harmonic functions.
		\textit{Comm. Pure Appl. Math.} \textbf{61} (2008), no. 1, 54-95. 
		
		\bibitem{Hua}
		L.-K. Hua, Starting with the unit circle, Springer-Verlag, New York-Berlin, 1981, xi+179 pp.
		
		
			
			\bibitem{Juhl}
                         A. Juhl, Families of conformally covariant differential operators,  $Q$-curvature and holography. \textit{Progr. Math.,} \textit{275} Birkh\"auser Verlag, Basel, 2009, xiv+488 pp.
			
			\bibitem{Li} Y. Y. Li,
			Remark on some conformally invariant integral equations: the method of moving spheres.
			\textit{J. Eur. Math. Soc. (JEMS)} \textbf{6} (2004), no. 2, 153-180.
			
			\bibitem{Li-Zhang}
Y. Y. Li and L. Zhang, Liouville-type theorems and Harnack-type inequalities for semilinear elliptic equations. \textit{Journal d'Analyse Mathematique} \textbf{90} (2003), 27-87.
			
				\bibitem{Li-Zhu}
				Y. Y.  Li and  M. Zhu,
			Uniqueness theorems through the method of moving spheres.
			\textit{Duke Math. J.} \textbf{80} (1995), no.2, 383-417.
			
			
			\bibitem{Martinazzi}
		 L. Martinazzi, Classification of solutions to the higher order Liouville's equation on  $\R^{2m}$.
\textit{Math. Z.} \textbf{263} (2009), no. 2, 307-329.

                   \bibitem{Obata}
M. Obata, The conjectures of conformal transformations of Riemannian manifolds. \textit{Bull. Amer. Math. Soc.} \textbf{77} 1971 265-270.

 \bibitem{Stafford}
	                R. Stafford, Tractor calculus and invariants for conformal sub-manifolds, Master's thesis, The University of Auckland, 2006.
			
			
				\bibitem{Sun-Xiong} 
				L. Sun and J. Xiong,
			Classification theorems for solutions of higher order boundary conformally invariant problems, I.
			\textit{J. Funct. Anal. } \textbf{271} (2016), no. 12, 3727–3764.
			
			
			
			
			
			
		\end{thebibliography}

	\end{document}